\documentclass[10pt, a4paper]{article}
\usepackage[a4paper, margin=1.1in]{geometry}

\title{Abelian surfaces over $\Fq(t)$ with large Tate--Shavarevich groups}
\author{\Large Martin \textsc{Azon}}
\date{}

\usepackage[utf8]{inputenc}
\usepackage{amsmath}
\usepackage{amsthm}
\usepackage{amssymb}
\usepackage{amsfonts}
\usepackage{graphicx}
\usepackage{tikz-cd}
\usepackage{enumerate}
\usepackage{hyperref}
\usepackage{xcolor}
\usepackage{pbox,tkz-graph}
\usetikzlibrary{positioning}
\usetikzlibrary{decorations.pathmorphing}
\usepackage{float}
\usepackage{caption}
\usepackage[colorinlistoftodos]{todonotes}
\usepackage{mathrsfs}
\usepackage{multicol}
\usepackage{comment}

\hypersetup{colorlinks, linkcolor={red!70!black}, citecolor={orange!90!yellow},}

\usetikzlibrary{decorations.pathmorphing, arrows.meta,bending,positioning,patterns,decorations.pathreplacing}

\newtheorem{lemma}{Lemma}[section]
\newtheorem{theorem}[lemma]{Theorem}
\newtheorem{proposition}[lemma]{Proposition}
\newtheorem{corollary}[lemma]{Corollary}

\newtheorem{remark}[lemma]{Remark}

\theoremstyle{definition}
\newtheorem{definition}[lemma]{Definition}

\let\originalleft\left
\let\originalright\right
\renewcommand{\left}{\mathopen{}\mathclose\bgroup\originalleft}
\renewcommand{\right}{\aftergroup\egroup\originalright}

\renewcommand{\Re}{\mathsf{Re}}
\newcommand{\cf}{\textit{cf.}}
\newcommand{\ie}{\textit{i.e.}}

\newcommand{\ord}{\mathrm{ord}}
\newcommand{\divi}{\, | \, }
\newcommand{\st}{^{\times}}
\newcommand{\prun}{\operatorname{pr}_{1}}

\renewcommand{\mod}{\operatorname{\, mod}}
\newcommand{\lcm}{\operatorname{lcm}}
\newcommand{\charact}{\operatorname{char}}
\newcommand{\rank}{\operatorname{rk}}
\newcommand{\rkan}{\operatorname{rk}_{\mathrm{an}}(S)}
\newcommand{\image}{\mathrm{im}}
\newcommand{\dist}{\operatorname{dist}}

\newcommand{\abs}[1]{\left\arrowvert #1 \right\arrowvert }

\newcommand{\Mpl}{M_{\K}}
\newcommand{\dv}{d_{v}}
\newcommand{\vbar}{\overline{v}}
\newcommand{\infbar}{\overline{\infty}}

\newcommand{\ensemble}[2]{\left\lbrace #1 \ \left\vert \  #2  \right. \right\rbrace}

\newcommand{\K}{\mathbb{K}}
\newcommand{\C}{\mathbb{C}}
\newcommand{\Z}{\mathbb{Z}}
\newcommand{\Q}{\mathbb{Q}}
\newcommand{\R}{\mathbb{R}}
\newcommand{\Ff}{\mathbb{F}}
\newcommand{\Qbar}{\overline{\Q}}
\newcommand{\Zbar}{\overline{\Z}}
\newcommand{\Qlbar}{\overline{\Q_{\ell}}}
\newcommand{\Fpbar}{\overline{\Fp}}
\newcommand{\Fqbar}{\overline{\Fq}}
\newcommand{\Fvbar}{\overline{\Fv}}
\newcommand{\Ksep}{\K^{\mathrm{sep}}}

\newcommand{\Fp}{\Ff_p}
\newcommand{\Fq}{\Ff_q}
\newcommand{\Fqa}{\Ff_{q^a}}
\newcommand{\Ia}{\mathbb{I}_{a}}
\newcommand{\Pqa}{\mathbf{P}_{q}(a)}
\newcommand{\Oqa}{\mathbf{O}_{q}(a)}
\newcommand{\Oqamax}{\mathbf{O}_{q}^{\mathrm{max}}(a)}
\newcommand{\Pqamax}{\mathbf{P}_{q}^{\mathrm{max}}(a)}
\newcommand{\qa}{q^a}

\newcommand{\Fv}{\Ff_v}

\newcommand{\Fo}{\Ff_o}

\newcommand{\pfrac}{\mathfrak{p}}
\newcommand{\valp}{\ord_{\pfrac}}
\newcommand{\Frv}{\text{Fr}_v}
\newcommand{\Fr}{\operatorname{Fr}_{q}}
\newcommand{\Tr}[2]{\operatorname{Tr}_{#1 / #2}}
\newcommand{\Nrm}[2]{\operatorname{N}_{#1 / #2}}
\newcommand{\Kloos}[3]{\operatorname{Kl}_{#1}\left(#2 \, , #3 \right)}
\newcommand{\Gauss}[3]{\operatorname{G}_{#1}\left(#2 \, , #3 \right)}
\newcommand{\kl}[3]{\kappa_{1, \, #1}^{ }\left(#2 \, , #3 \right)}
\newcommand{\klbis}[3]{\kappa_{2, \, #1}^{ }\left(#2 \, , #3 \right)}

\newcommand{\gamo}{\gamma(o)}
\newcommand{\kapv}{\kappa_{1}^{ }(v)}
\newcommand{\kapvbis}{\kappa_{2}^{ }(v)}

\newcommand{\Klo}{\operatorname{Kl}(o)}
\newcommand{\kapo}{\kappa_{1}^{ }(o)}
\newcommand{\kapobis}{\kappa_{2}^{ }(o)}

\newcommand{\Fqt}{\Fq(t)}

\newcommand{\sep}[1]{#1^{\mathrm{sep}}}
\newcommand{\unr}[1]{#1^{\mathrm{unr}}}

\newcommand{\Ov}{\mathcal{O}_{v}}
\newcommand{\Ql}{\Q_{\ell}}
\newcommand{\Zl}{\Z_{\ell}}
\newcommand{\Gal}{\operatorname{Gal}}
\newcommand{\TlS}{\operatorname{T}_{\ell}(S)}
\newcommand{\VlS}{\operatorname{V}_{\ell}(S)}
\newcommand{\cvbar}{c_{\overline{v}}(S)}
\newcommand{\Oinf}{\mathcal{O}_{\infty}}

\newcommand{\ArtSch}{\wp_a}
\newcommand{\As}{\ArtSch(t)}

\newcommand{\Wrst}{\mathcal{W}}
\newcommand{\Cmod}{\mathcal{C}}
\newcommand{\Kmod}{\mathcal{K}}
\newcommand{\Smod}{\mathcal{S}}
\newcommand{\Xmod}{\mathcal{X}}
\newcommand{\Emod}{\mathcal{E}}
\newcommand{\Xhat}{\widetilde{\mathcal{X}}}
\newcommand{\Zsurf}{\mathcal{Z}}
\newcommand{\Zsurfa}{\mathcal{Z}_{a}}
\newcommand{\res}{\mathsf{res}}

\newcommand{\Sa}{S_{a}}
\newcommand{\Uzero}{\mathcal{U}_{0}}
\newcommand{\Uinft}{\mathcal{U}_{\infty}}
\newcommand{\Pproj}[2]{\mathbb{P}_{#2}^{#1}}
\newcommand{\Spec}{\operatorname{Spec}}
\newcommand{\Pp}{\Pproj{1}{\Fq}}
\newcommand{\Ppo}{\Pproj{1}{ }}

\newcommand{\Sner}{\mathcal{S}^{0}}
\newcommand{\Snerv}{\mathcal{S}_{v}^{0}}

\newcommand{\Ga}{\mathbb{G}_a}
\newcommand{\Gm}{\mathbb{G}_m}
\newcommand{\Pic}{\operatorname{Pic}}
\newcommand{\NS}{\operatorname{NS}}
\newcommand{\ShaS}{\Sh(S / K)}
\newcommand{\Reg}{\operatorname{Reg}}

\newcommand{\tor}[1]{#1_{\mathrm{tor}}}

\newcommand{\GXa}{G_{X,a}}

\newcommand{\Div}{\operatorname{Div}}
\newcommand{\Jac}{\operatorname{Jac}}
\newcommand{\Het}{\operatorname{H}_{\acute{e}t}^{i}}
\renewcommand{\H}{\operatorname{H}}
\newcommand{\HS}{\operatorname{H}^{1}(S)}
\newcommand{\chpol}[2]{\operatorname{det}\left(1 - \Fr T \, \left | \, \H^{#1}(#2) \right. \right)}
\newcommand{\chpolIn}[2]{\operatorname{det}\left(1 - \Frv T^{\dv} \left | \, \H^{#1}(#2)^{I_v} \right. \right)}
\newcommand{\chpolv}[2]{\operatorname{det}\left(1 - \Frv \, T^{\dv} \left | \, \H^{#1}(#2) \right. \right) }
\newcommand{\chpolg}[1]{\operatorname{det}\left(1 - \Fr T \, \left | \, #1 \right. \right)}
\newcommand{\Trpol}[2]{\operatorname{Tr}\left( #1 \, \left | \, #2 \right. \right)}
\newcommand{\spval}{L^{*}(S)}
\newcommand{\spvala}{L^{*}(\Sa)}

\newcommand{\e}{\mathrm{e}}
\newcommand{\pig}{\pi_{g}}
\newcommand{\height}{\textit{ht}}
\newcommand{\normun}[1]{\left \lVert #1 \right \rVert_{1}}
\newcommand{\norminf}[1]{\left \lVert #1 \right \rVert_{\infty}}
\newcommand{\eps}{\varepsilon}
\newcommand{\epspm}{\eps_{\Ff}^{ }(\chi, \psi)}
\newcommand{\epspmbis}{\eps_{\Ff}^{ }(\chi, \psi')}
\newcommand{\thet}{\theta}
\newcommand{\thetpm}{\thet_{\Ff}^{ }(\psi, \alpha)}
\newcommand{\GamFn}{\Gamma_{\Ff, n}^{ }}
\newcommand{\epso}{\varepsilon_{o}}
\newcommand{\epsoprim}{\varepsilon_{o'}}
\newcommand{\theto}{\theta_{o}}
\newcommand{\thetv}{\theta_{v}}

\newcommand{\intrv}{[0, \pi ]}
\newcommand{\interv}{[0, 2 \pi)}
\newcommand{\intv}{\left[0, \frac{\pi}{2r}\right]}

\renewcommand{\d}{\mathrm{d}}

\newcommand{\nua}{\nu_{a}}

\newcommand{\dnua}{\d \nua}

\newcommand{\dlthetv}{\delta \lbrace \theta_v \rbrace}
\newcommand{\nust}{\nu_{\mathrm{ST}}}
\newcommand{\dnust}{\d \nust}
\newcommand{\dx}{\d x}

\newcommand{\Discrep}{\mathcal{D}_{q}^{*}(a)}

\newcommand{\F}{\mathcal{F}}
\newcommand{\Klsh}{\mathcal{K}l}
\newcommand{\GL}{\operatorname{GL}}

\newcommand{\SL}{\operatorname{SL}}
\newcommand{\SU}{\operatorname{SU}}
\newcommand{\nat}{^{\natural}}
\newcommand{\semis}{^{\operatorname{s.s.}}}
\newcommand{\Diag}{\operatorname{Diag}}
\newcommand{\Trace}{\operatorname{Tr}}

\newcommand{\Disc}{\mathcal{D}_{r, \lambda}}
\newcommand{\Disca}{\mathcal{D}_{r_{a}, \lambda_{a}}}

\newcommand{\beteta}{\beta_{\eta}^{r, \lambda}}
\newcommand{\Wlam}{W^{r, \lambda}}
\newcommand{\Wlama}{W^{r_{a}, \, \lambda_{a}}}

\newcommand{\Weta}{W_{\eta}^{r, \lambda}}
\newcommand{\Wetaa}{W_{\eta}^{r_a, \lambda_a}}
\newcommand{\Wprim}{(\Weta)'}
\newcommand{\Wpri}{W_{\eta}'}
\newcommand{\nrmbet}{\norminf{\beta_{0}'}}
\newcommand{\Delteta}{\Delta_{\eta}^{r, \lambda}}

\newcommand{\BS}{\mathfrak{Bs}}

\newcommand{\etalchar}[1]{$^{#1}$}

\newcommand{\twin}{\mathfrak{t}}
\newcommand{\Rroots}{\mathcal{R}}

\newcommand{\refGauss}[1]{\hyperref[#1]{\textsf{(Ga\,\ref*{#1})}}}
\newcommand{\refKloos}[1]{\hyperref[#1]{\textsf{(Kl\,\ref*{#1})}}}

\newcommand{\map}[5]{
	#1:\left \lbrace 
	\begin{array}{ccc}
		#2 & \longrightarrow & #3 \\
		#4 & \longmapsto & #5 \\
	\end{array}\right.
}

\newcommand{\ratmap}[5]{
	#1:\left \lbrace 
	\begin{array}{ccc}
		#2 & \dashrightarrow & #3 \\
		#4 & \longmapsto & #5 \\
	\end{array}\right.
}

\newcommand*{\EnsembleQuotient}[2]%
{\ensuremath{%
		#1/\!\raisebox{-.65ex}{\ensuremath{\mathcal{#2}}}}}

\numberwithin{equation}{section}

\DeclareFontFamily{U}{wncy}{}
\DeclareFontShape{U}{wncy}{m}{n}{<->wncyr10}{}
\DeclareSymbolFont{mcy}{U}{wncy}{m}{n}
\DeclareMathSymbol{\Sh}{\mathord}{mcy}{"58}

\begin{document}

\maketitle

\begin{abstract}
	We produce an explicit sequence $\left(S_a \right)_{a \geq 1}$ of abelian surfaces over the rational function field $\Fq(t)$ whose Tate--Shafarevich groups are finite and large. More precisely, we establish the estimate $\abs{\Sh (S_a)} = H(S_a)^{1 + o(1)}$ as $a \rightarrow \infty$, where $H(S_a)$ denotes the exponential height of $S_a$. Our method is to prove that each $S_a$ satisfies the BSD conjecture, analyse the geometry and arithmetic of its Néron model and give an explicit expression for its $L$-function in terms of Gauss and Kloosterman sums. By studying the relative distribution of the angles associated to these character sums, we estimate the size of the central value of $L(S_a, T)$, hence the order of $\Sh(S_a)$.
\end{abstract}

\section{Introduction}~

\indent Consider a finite field $\Fq$ of characteristic $p \geq 7$ and let $\K := \Fq(t)$ be the rational function field over $\Fq$. Let $A$ be an abelian variety defined over $\K$, the Tate--Shafarevich group $\Sh(A)$ is the group of isomorphism classes of principal homogeneous spaces of $A$ over $\K$ that become trivial in all the local completions of $\K$ (see \cite{Hindry07} for more details on this definition). It is an arithmetic object that measures the extent to which the local-global principle fails for $A$. Many questions regarding the Tate--Shafarevich group are still open: its finiteness, while proven in various cases, remains conjectural in general. A natural question to ask is: assuming finiteness of $\Sh(A)$, what can be said about its size? Can we compare it to other arithmetic invariants of $A$? In \cite{GoldfeldSzpiro} and \cite{Rajan}, Goldfeld and Szpiro and, independently, Rajan, addressed the case of elliptic curves and gave upper bounds on the order of the Tate--Shafarevich group in terms of the degree of the conductor and the exponential height. More recently, in \cite{HindryPacheco}, Hindry and Pacheco relate the order of $\Sh(A)$ to the exponential height $H(A)$, which measures the arithmetic complexity of $A$.

\begin{theorem}[Hindry-Pacheco \cite{HindryPacheco}]
	Let $d \geq 1$ be an integer and $\eps > 0$ be a real number. There exists a real constant $c_{\eps, d}>0$ such that for any abelian variety $A / \K$ of dimension $d$ with finite Tate--Shafarevich group, the following estimate holds:
	\begin{equation*}
		\abs{\Sh(A)} \leq c_{\eps, d} \, H(A)^{1 + \eps}.
	\end{equation*}
\end{theorem}

Kato and Trihan proved in \cite{KatoTrihan} that the finiteness of $\Sh(A)$ implies the Birch and Swinnerton-Dyer conjecture for $A$ (hereafter abbreviated as BSD). Under the hypothesis of the theorem, Hindry and Pacheco use the BSD formula to relate $\Sh(A)$ to the special value $L^{*}(A)$ of the $L$-function of $A$. By bounding various arithmetic invariants of $A$ and employing analytic techniques, they manage to estimate the size of $L^{*}(A)$, hence deriving their theorem.

The main result of this article is that in dimension $2$, this bound can essentially not be improved, meaning that the exponent $1$ cannot be replaced by a smaller quantity:

\begin{theorem}[= \ref{thm: Sha compared to H}]\label{thm: Main thm introd}
	There exists a sequence $\left(S_a \right)_{a \geq 1}$ of non-isotrivial abelian surfaces over $\K$, all having finite Tate--Shafarevich groups, satisfying the following property: for any $\eps > 0$, there exist positive constants $c_{\eps}', c_{\eps}'' > 0$ such that for all $a \geq 1$, 
	\begin{equation*}
		c_{\eps}' \, H(S_a)^{1 - \eps} \leq \abs{\Sh(S_a)} \leq c_{\eps}'' \,  H(S_a)^{1 + \eps}.
	\end{equation*} 
\end{theorem}

Our proof of this theorem is constructive and effective: we explicitly construct such a sequence of abelian surfaces. We build $S_a$ as the Jacobian of a hyperelliptic genus $2$ curve in a tower of Artin--Schreier extensions. Following a construction by Pries and Ulmer (see \cite{PriesUlmer}), we define, for any $a \geq 1$, $C_a / \K$ to be the smooth projective curve described by the hyperelliptic equation
\begin{equation*}\label{eq: Hyp model introd}
	C_{a} : y^2 = \left(x^3 + \As - 2 \right) \left( x^3 + \As + 2 \right), \quad \text{where } \As := t^{q^a} - t \in \Fqt.
\end{equation*} 
We let $S_a := \Jac(C_a)$ be the Jacobian of $C_a$. One of the main goals of this article is to describe several arithmetic invariants of $S_a$ and their variation as $a$ grows.

Recently, there has been a growing interest in the explicit arithmetic of Jacobians of hyperelliptic curves. In \cite{D2M2}, the authors introduce the notion of \textit{cluster picture}, a combinatorial object that helps understand the local arithmetic of a hyperelliptic curve and its Jacobian. In Section~\ref{sect: Studying C and S}, we apply this theory in our context to describe various invariants of $S_a$. We explicitly obtain:

\begin{proposition}[= \ref{cor: Conductor} + \ref{prop: Tamagawa} + \ref{prop: Height}]
	For any $a \geq 1$: 
	\begin{itemize}
		\item The local exponents of the conductor are $n_{v} = 2$ if $v \divi \As^2 - 4$, $n_{\infty} = 4$ and $n_{v} = 1$ for all other places $v$.
		\item The Tamagawa numbers of $S_a$ satisfy $c_{v}(S_a) = 1$ if $v \neq \infty$ and $c_{\infty}(S_a) \in \lbrace 1, 3, 9 \rbrace$.
		\item The differential height of $S_a$ is given by $h(S_a) = q^a + 1$.
	\end{itemize}
\end{proposition}

Since $S_a$ is a Jacobian, its arithmetic is intimately related to the geometry of the underlying curve $C_a$, which is well-described \textit{via} its minimal regular model $\Cmod_{a} \rightarrow \Pp$. The family $\left(S_a \right)_{a \geq 1}$ constitutes an \textit{Artin--Schreier} family of abelian varieties (such families have been extensively studied in \cite{PriesUlmer}). In our context, one can analyse the geometry of $\Cmod_{a}$ and prove that, as a surface over $\Fq$, it is dominated by a product of two curves. Building upon the work of Shioda and Tate, we deduce from the general framework of \cite{PriesUlmer}:

\begin{theorem}[= \ref{thm: BSD}]\label{prop: Invariants intro}
	For any $a \geq 1$, the abelian surface $S_a$ satisfies the BSD conjecture:
	\begin{itemize}
		\item The algebraic and analytic ranks of $S_a$ coincide, i.e. $\rank \left(S_{a}(\K) \right) = \ord_{T = q^{-1}} L(S_a, T)$.
		
		\item The group $\Sh(S_a)$ is finite.
		
		\item The BSD formula holds:
		\begin{equation*}
			\spvala = \frac{\abs{\Sh(S_a)} \Reg(S_a) \, q^2 \, \prod_{v} c_v(S_a)}{\abs{\tor{S_{a}(\K)}}^{2} H(S_a)},
		\end{equation*}
		where $\Reg(S_a)$ is the Néron--Tate regulator of $S_a$.
	\end{itemize}
\end{theorem}

The BSD conjecture paves the way for using analytic techniques to compute the rank of $S_a(\K)$ and better understand other arithmetic data of $S_a$, such as $\abs{\Sh(S_a)}$ and $\Reg(S_a)$. The $L$-function $L(S_a, T)$ has a cohomological interpretation: it is the characteristic polynomial of the Frobenius on certain étale cohomology groups. Using the explicit construction of our hyperelliptic curves $C_a$, we compute in Section~\ref{sect: L func} the $L$-function of $S_a$ as follows. By studying the reduction types of the Néron model of $S_a$, we first establish a connection between $L(S_a, T)$ and the zeta function $Z(\Cmod_{a}, T)$ of the surface $\Cmod_{a} / \Fq$. As we stated above, there exist two smooth projective curves $X_a$, $Y_a$ defined over $\Fq$ such that $\Cmod_{a}$ is dominated by $X_a \times Y_a$. We study in deeper detail the geometry of $X_a \times Y_a \dashrightarrow \Cmod_{a}$ and describe $Z(\Cmod_{a}, T)$ in terms of the cohomology groups of $X_a$ and $Y_a$. This way we obtain a completely explicit expression for $L(S_a, T)$ which is computable in practice.

\begin{theorem}[= \ref{thm: L func formula}]\label{thm: L func introd}
	For any $a \geq 1$, the $L$-function of $S_a$ is given by the formula
	\begin{equation*}
		L(S_a, T) = \prod_{o \, \in \Oqa} \left( 1 - \gamo \Klo \,  T^{\abs{o}} + \gamo^2 \, q^{\abs{o}} \, T^{2 \abs{o}} \right) \in \Z[T],
	\end{equation*} 
	where $\Oqa$ is a finite set of orbits, and for each $o \in \Oqa$, $\gamo$ is a Gauss sums and $\Klo$ is a Kloosterman sum. Precise definitions are given in section~\ref{sect: Orbits}.
\end{theorem}

Various authors have previously computed $L$-functions of abelian varieties over function fields. For example, in \cite{GriDWit} and \cite{ExplJacLeg}, the authors use the Grothendieck--Lefschetz trace formula to compute the local factors of the $L$-function through point-counting arguments and direct manipulation of character sums, and then recover a global formula. The method employed here is inspired by a theorem of Shioda (\cite{Shioda}), who treated the case of elliptic curves , and \cite{GriUlm,AGTT} where the authors provide cohomological computations of $L$-functions.

Theorem~\ref{thm: L func introd} describes $L(S_a, T)$ in terms of Gauss and Kloosterman sums, which are rather well-understood: in Section \ref{sect: Orbits} we recall several properties of these sums (their size, $p$-adic valuation, ...) that are classical (see \cite{LidlNiederreiter, Cohen}). Moreover, it is possible to associate angles to Gauss and Kloosterman sums: for any $o \in \Oqa$, there are $\epso \in \interv$ and $\theto \in \intrv$ such that $\gamo = q^{\abs{o} / 2} \, \e^{i \epso}$ and $\Klo = q^{\abs{o} / 2} \, \cos(\theto)$ (see Sections~\ref{sect: Special value section} and~\ref{sect: Equid sect} for details). Using the explicit expression for $L(S_a, T)$ and manipulating trigonometric formulae we establish an explicit expression for the special value $\spvala$ in terms of Gauss and Kloosterman angles. Using the properties of the character sums mentioned above, we deduce the special value and the rank of the Mordell-Weil group.

\begin{proposition}[= \ref{thm: Special value thm} + \ref{cor: Alg rank 0}]
	For any $a \geq 1$, we have
	\begin{equation*}
		L(S_a, q^{-1}) = \prod_{o \, \in \Oqa} 4 \, \abs{\sin \left( \frac{\epso + \theto}{2}\right) \, \sin \left( \frac{\epso - \theto}{2}\right)} \neq 0.
	\end{equation*}
	We deduce that, $\spvala = L(S_a, q^{-1})$, hence $\rank(S_{a}(\K)) = 0$ and $\Reg(S_a) = 1$.
\end{proposition}

We then focus on the archimedean size of the special value $\spvala$. The regulator of $S_a$ being trivial, we use the BSD formula to estimate the order of $\Sh(S_a)$. This point of view has revealed to be very fruitful: the interested reader is referred to \cite{HindryPacheco, Gri19, GriDWit} for some examples where the authors approach the arithmetic of abelian varieties \textit{via} an analytic study of their $L$-function.

Estimating the size of the special value involves understanding the distribution of Gauss and Kloosterman angles. First, we use Katz's work \cite{Katz88} to establish an equidistribution result concerning Kloosterman angles. More precisely, if one lets $\nua$ be the average of Dirac measures supported at Kloosterman angles, then the sequence $(\nua)_{a \geq 1}$ converges weak--$*$ to the Sato--Tate measure $\nust:= \frac{2}{\pi}\sin^{2}\theta \, \d\theta$. Moreover, following \cite{Gri19}, this equidistribution result is made effective and we give an upper bound on the error term. After that we study the minimal distance between Gauss and Kloosterman angles, topic that was partly addressed in \cite{Gri19, GriDWit}. In Section \ref{sect: Distance angles} we generalise the results of \textit{loc. cit.} and establish a Diophantine statement asserting that Kloosterman angles ``avoid" Gauss angles and cannot be arbitrarily close to them. 

With these results at hand, employing analytic techniques we estimate the size of the special value in terms of the exponential height. 

\begin{theorem}[= \ref{thm: Ratio spval height thm}]\label{thm: Ratio spval height introd}
	For any $a \geq 1$, we have
	\begin{equation*}
		\frac{\log \spvala}{\log H(S_a)} = O \left( \frac{1}{a}\right).
	\end{equation*}
\end{theorem}

Theorem \ref{thm: Main thm introd} follows from Theorem \ref{thm: Ratio spval height introd} using the BSD formula and the estimates on the arithmetic invariants of $S_a$ obtained in Proposition \ref{prop: Invariants intro}. We conclude by deriving an analogue of the Brauer--Siegel theorem for the abelian surfaces $S_a$.  Our family provides a new example of non-isotrivial higher-dimensional abelian varieties whose Brauer--Siegel ratio $\BS(S_a):=\log(\abs{\Sh(S_a)} \Reg(S_a))/ \log H(S_a)$ tends to $1$ as $a \rightarrow \infty$. This provides some more evidence towards the phenomenon (already observed in \cite{HindryPacheco, Ulmer19} and related work) that the Brauer--Siegel ratio of an abelian variety over $\K$ should ``often" be near $1$ when the exponential height is large.


\section[The arithmetic of the abelian surfaces Sa]{The arithmetic of the abelian surfaces $S_a$}\label{sect: Studying C and S}

In this section we introduce the curve $C_a / \K$, its Jacobian $S_a / \K$ and we describe different models of those. In particular, we use the minimal regular model of $C_a$ to compute the conductor of $S_a$, and the minimal regular simple normal crossings model of $C_a$ to describe the differential height $h(S_a)$ and the Tamagawa numbers $c_v(S_a)$.

We fix some notation first: let $p \geq 7$ be a prime number, $\Fq$ be a finite field of characteristic $p$ and let $\K := \Fq(t)$ be the function field of $\Pp$. We write $\Ksep$ for a separable closure of $\K$ and $\Mpl$ for the set of places of $\K$. Recall that $\Mpl$ is in bijection with the set of closed points of $\Pp$. For any $v \in \Mpl$, we denote the completion of $\K$ at $v$ by $\K_v$, its ring of integers by $\Ov$, and a uniformizer of $\Ov$ by $\pi_{v}$. Let $\Fv$ be the residue field at $v$, and $\dv := \left[\Fv : \Fq\right]$ be the degree of $v$. Denote by $\unr{\K_{v}}$ the maximal unramified extension of $\K_{v}$ and $\unr{\Ov}$ for its ring of integers: its residue field is an algebraic closure of $\Fv$. Moreover, we let $\infty$ denote the place at infinity of $\K$, corresponding to the valuation defined on $\Fq[t]$ by $P \mapsto - \deg P$, where $\deg P$ denotes the degree of the polynomial $P$. Equivalently, $\infty$ corresponds to the $t^{-1}$-adic valuation on $\K$.

\subsection[Constructing the curves Ca and their Jacobians]{Constructing the curves $C_a$ and their Jacobians}\label{sect: Constructing ab surf}

We construct the abelian surfaces $S_a$ as the Jacobians of certain genus $2$ curves, in such a way that they satisfy the BSD conjecture (see \cite{HindryPacheco} or Section \ref{sect: Def Lfunc and BSD} for a precise statement of the BSD conjecture). Building upon the work of Shioda and Tate, the following result concerning Jacobians of curves over $\K$ holds:

\begin{theorem}\label{thm: DPC implies BSD}
	Consider a smooth, proper, geometrically irreducible surface $\Zsurf$ defined over $\Fq$. Assume that $\Zsurf$ is equipped with a generically smooth morphism $\Zsurf \rightarrow \Pp$ and denote by $C / \K$ its generic fiber. If the surface $\Zsurf / \Fq$ is dominated by a product of curves, then the Jacobian of $C / \K$ satisfies the BSD conjecture.
\end{theorem}

\begin{proof}[Skecth of the proof]
	This follows from the equivalence of the BSD conjecture with Tate's conjecture over finite fields, and the fact that products of curves satisfy Tate's conjecture. We refer the reader to \cite[Theorems 6.3.1, 8.1.2 and Lemma 8.2.1]{Ulmer14}.
\end{proof}

In Theorem~\ref{thm: BSD} we will state a precise result in our context and give more details about its proof. The goal of this section is to exploit this phenomenon to construct the curves $C_a$ and their Jacobians $S_a$ presented in the introduction. The discussion below is a particular case of a more general framework, described by Pries and Ulmer in \cite{PriesUlmer}. In that article, the authors give an explicit recipe to construct families of abelian varieties in Artin--Schreier extensions, such that all the terms in the family satisfy the BSD conjecture.\\

Fix an integer $a \geq 1$, consider the morphism $\ArtSch : \Pp \rightarrow \Pp$ defined on $\mathbb{A}_{\Fq}^{1}$ by  $\ t \mapsto t^{\qa} - t$, which we call the \textit{Artin--Schreier} morphism. Let $X_{a}$ and $Y_{a}$ be the smooth projective curves over $\Fq$ described on affine subsets of $\Pproj{2}{\Fq}$ by the equations  
\begin{equation}\label{eq: Models X Y}
	X_{a} : u^3 = \ArtSch(t_1) \quad \text{and} \quad Y_{a}: v + \frac{1}{v}= \ArtSch(t_2).
\end{equation}
These are instances of so-called \textit{Artin--Schreier} curves: we will study in detail their geometry in Section~\ref{sect: Geom Cmod}. Let $\Zsurfa$ be the projective irreducible surface over $\Fq$ defined on an affine open subset by 
\begin{equation}\label{eq: Model Zsurfa}
	\Zsurfa :  v + \frac{1}{v} - u^3 = \As.
\end{equation}
The morphism $\Zsurfa \rightarrow \Pp$ extending $(u, v, t) \mapsto t$ turns $\Zsurfa$ into a fibered surface over $\Pp$. Using the Jacobian criterion one checks that the generic fiber of $\Zsurfa \rightarrow \Pp$ is a smooth curve over $\K$, which we denote by $C_a$. The arithmetic surface $\Zsurfa$ is a model of the curve $C_a$ over $\K$ and equation \eqref{eq: Model Zsurfa}, viewed as an element of $\K(u, v)$, defines a dense open affine subset of $C_a / \K$. The surface $\Zsurfa$ is not regular, but it is birational to a regular proper irreducible surface fibered over $\Pp$, called the minimal regular model of $C_a$ (we will construct this and other models of $C_a$ in Section~\ref{sect: Models}). 

The rational map 
\begin{equation*}
	\ratmap{\vartheta_{0}}{X_a \times Y_a}{\Zsurf_a}{(u, t_1, v, t_2)}{(u, v, t_2 - t_1),}
\end{equation*}
is well-defined because of the relation $\ArtSch(t_2 - t_1) = \ArtSch(t_2) - \ArtSch(t_1) $. We will show in Proposition~\ref{prop: DPC} that the map $\vartheta_{0}$ is dominant. Since $\Zsurfa$ is birational to the minimal regular model of $C_a$, Theorem~\ref{thm: DPC implies BSD} states that $\Jac(C_a)$ satisfies the BSD conjecture (see Theorem~\ref{thm: BSD} for a precise statement). 

In the terminology of \cite{PriesUlmer}, the curves $C_a$ correspond to the type $(3, 1+1)$, our choice of rational functions being $f(u) = u^3$ (which has a pole of order $3$ at infinity) and $g(v) = v + v^{-1}$ (which has two poles of order $1$). Using Proposition $3.1.5$ of the mentioned article, we see that $C_a$ has genus $2$, and is thus hyperelliptic. The hyperelliptic involution, read off from equation \eqref{eq: Model Zsurfa}, sends $(u, v)$ to $(u, 1/v)$ and a hyperelliptic equation of $C_a$ is:
\begin{equation}\label{eq: Hyp model}
	C_{a} : y^2 = \left(x^3 + \As - 2 \right) \left( x^3 + \As + 2 \right).
\end{equation} 
The change of variables is explicitly given by $(u, v) \mapsto (x, y) = \left(u, 2v-u^3 - \As\right).$
\vspace{0.5em}

\begin{definition}
	For any $a \geq 1$, let $C_a$ be the genus $2$ hyperelliptic curve defined over $\K$ by equation \eqref{eq: Hyp model}, and define $S_a := \Jac(C_a)$ to be its Jacobian, which is an abelian surface over $\K$.
\end{definition}

\begin{remark}
	For any $a \geq 1$, let $\K_{a} := \K[z] / (\ArtSch(z) - t)$, then $\K_{a} / \K$ is a field extension of degree $q^a$, which corresponds to the covering $\ArtSch : \Pp \rightarrow \Pp$. We say that $\K_a / \K$ is an \textit{Artin--Schreier} extension. One can also view $C_a$ as the base change from $\K$ to $\K_a$ of the hyperelliptic curve $C_0 / \K$ defined by $y^2 = (x^3 + t - 2 )( x^3 + t + 2)$. In the same way $S_a$ is isomorphic to the base change $\Jac (C_0) \times_{\K} \K_a$. We say that $(S_{a})_{a \geq 1}$ is an Artin--Schreier family of abelian surfaces.
\end{remark}

\subsection[Models of Ca and Sa]{Models of $C_a$ and $S_a$}\label{sect: Models}

In this section we discuss models of $C_a$ and $S_a$ over $\Pp$. Since each $C_a$ is a hyperelliptic curve, there are many explicit tools that help us understand the geometry and the arithmetic of its Jacobian $S_a$. We refer to \cite[§9 \& 10]{Liubook} and to \cite[§9]{BoschLutkebohmerRaynaud} for basics about models of curves and Néron models of Jacobians respectively. From now on, we fix the integer $a \geq 1$ and, unless stated otherwise, we omit, the subscript $a$ and write, for example, $C$ instead of $C_a$.

\vspace{1em}

First, we recall from \cite{Liu23} the definition of a Weierstrass model associated to a hyperelliptic equation of a hyperelliptic curve over a PID.

\begin{definition}
	Let $A$ be a principal ideal domain, $K := \operatorname{Frac}(A)$ and $C / K$ be a hyperelliptic curve of genus $g$ defined by an integral Weierstrass equation 
	\begin{equation}\label{eq: Hyp eq Wrst model}
		y^2 = F(x), \qquad F(x) \in A[x].
	\end{equation}
	We define the Weierstrass model of $C$ over $\Spec(A)$ associated to equation \eqref{eq: Hyp eq Wrst model} to be the glueing of the affine schemes 
	\begin{equation*}
		\Spec \left( A[x, y] \, / \, (y^2 - F(x)) \right) \quad \text{ and } \quad \Spec \left(A[z, w] \, / \, (w^2 - F(z^{-1})) \right)
	\end{equation*}
	along the identification $(z, w) = (1/x, \, y/ x^{g+1})$. We define the discriminant of such a model as $\Delta := 2^{4g} \, \operatorname{disc}(F)$. We say that the hyperelliptic equation \eqref{eq: Hyp eq Wrst model} (or the associated Weierstrass model) is minimal at a place $\mathfrak{P}$ of $A$ if the valuation $\ord_{\mathfrak{P}}(\Delta)$ is minimal among all other integral hyperelliptic equations for $C$. 
\end{definition}

We now construct a Weierstrass model of $C$ over $\Pp$. Let $\Uzero := \Spec(\Fq[t])$ and $\Uinft := \Spec(\Fq[s])$ be two affine lines covering $\Pp$, which we glue along the morphism $t \mapsto s = t^{-1}$. Over $\Uinft$, replacing $t = s^{-1}$ in \eqref{eq: Hyp model} gives a non-integral equation of $C$. In order to construct an integral model of $C$ over $\Uinft$, let $b \in \lbrace 1, 2\rbrace$ be such that $3 \divi \qa + b$. Setting the change of variables $(x_1, y_1, s) := (t^{-(q^a +b)/3}x, ^{-(q^a +b)}y, s^{-1})$, equation \eqref{eq: Hyp model} becomes
\begin{equation}\label{eq: Integral eq infty}
	y_{1}^2 = \left(x_{1}^3 + s^b - s^{q^a+b} - 2s^{q^a + b} \right) \left(x_{1}^3 + s^b - s^{q^a+b} + 2s^{q^a + b} \right).
\end{equation}
This is an integral equation of $C$ over $\Oinf$. As we will see in Lemma \ref{lem: Wrst is minimal}, when $b=1$, \ie \ $q^a \equiv 2 \mod 3$, this equation is minimal. However, when $b = 2$, \ie \ $q^a \equiv 1 \mod 3$, a further change of variables is required to obtain a minimal equation. Following Liu's minimization algorithm (see \cite{Liu23}), we let $(x_2, y_2) := (s^{-1} x_1, s^{-2} y_1)$. Then \eqref{eq: Integral eq infty} becomes 
\begin{equation}\label{eq: Integ eq infty 2}
	y_{2}^2 = \left(s x_{2}^3 + 1 - s^{q^a - 1} - 2 s^{q^a} \right) \left(s x_{2}^3 + 1 - s^{q^a - 1} + 2 s^{q^a} \right).
\end{equation}

\begin{definition}\label{def: Wrst model}
	Fix an integer $a \geq 1$. We let $\Wrst_{0} \rightarrow \Uzero$ be the Weierstrass model associated to the hyperelliptic equation \eqref{eq: Hyp model}. Consider $b \in \lbrace 1, 2\rbrace$ be such that $3 \divi \qa + b$ as above. We let $\Wrst_{\infty} \rightarrow \Uinft$ be the Weierstrass model associated to the hyperelliptic equation \eqref{eq: Integral eq infty} if $b= 1$, and to \eqref{eq: Integ eq infty 2} if $b=2$. We define the \textit{Weierstrass model} of $C$ over $\Pp$, denoted by $\Wrst$, as the glueing of $\Wrst_{0}$ and $\Wrst_{\infty}$ along the morphism 
	\begin{equation*}
		(x, y, t) \mapsto \begin{cases}
			(x_1, y_1, s) = (t^{-(\qa + 1)/3}x, \, t^{-(\qa + 1)}y, \, t^{-1}) & \text{ if } b = 1 \text{ \ie \ } q^a \equiv 2 \mod 3, \\
			(x_2, y_2, s) = (t^{-(\qa - 1)/3}x, \, t^{-\qa}y, \, t^{-1}) & \text{ if } b = 2 \text{ \ie \ } q^a \equiv 1 \mod 3.
		\end{cases}
	\end{equation*} 
\end{definition}

The scheme $\Wrst$ is an irreducible, proper, normal surface over $\Fq$. The projections $\Wrst_{0} \rightarrow \Uzero$ and $\Wrst_{\infty} \rightarrow \Uinft$ glue into a morphism $\Wrst \rightarrow \Pp$, whose generic fiber is the curve $C$. The discriminant of $\Wrst$ is an effective divisor on $\Pp$, which we write in multiplicative notation. Recall that we identify finite places of $\K$ with monic irreducible polynomials.

\begin{lemma}\label{lem: Wrst is minimal}
	The model $\Wrst \rightarrow \Pp$ is minimal at all places of $\K$. The discriminant $\Delta(\Wrst)$ is given by
	\begin{equation*}
		\Delta(\Wrst) = \prod_{v \divi \As^2 - 4} (v)^{2} \ (\infty)^{6q^a + 10}.
	\end{equation*}
\end{lemma}

\begin{proof}
	One checks that the discriminant of a sextic polynomial of the form $\alpha x^6 + \beta x^3 + \gamma$ is given by $\operatorname{disc}(\alpha x^6 + \beta x^3 + \gamma) = -3^6 (\alpha \gamma)^2 \left(\beta^2 - 4\alpha \gamma\right)^3$. Therefore, the discriminant of $\Wrst_{0} \rightarrow \Uzero$ is obtained in terms of the discriminant of the defining polynomial, which we compute now:
	\begin{equation*}
		\operatorname{disc}\left(\left(x^3 + \As - 2 \right) \left( x^3 + \As + 2 \right) \right) = 2^{12} 3^{6} \left(\As^2 -4\right)^2.
	\end{equation*}
	But $\As^2 -4$ is a square-free polynomial in $t$, so we deduce that $\ord_{v}\left(\Delta(\Wrst_0)\right)=  \ord_{v}\left(\Delta(\Wrst)\right) = 2$ at any $v \divi \As^2 - 4$. Now Remark 7.3 in \cite{Liu96} implies that $\Wrst_0$ (and thus $\Wrst$) is minimal at any place $v$ of $\Uzero$. Applying the minimization algorithm given in \cite[§5.5]{Liu23} one checks that \eqref{eq: Integral eq infty} (resp. \eqref{eq: Integ eq infty 2}) is a minimal equation at $\infty$ when $q^a \equiv 2 \mod 3$ (resp. $q^a \equiv 1 \mod 3$). Again, a straightforward computation shows that the discriminant of both \eqref{eq: Integral eq infty} and \eqref{eq: Integ eq infty 2} is $6q^a + 10$, hence the result.
\end{proof}

\begin{remark}\label{rmk: Defining eqt infinity}
	When $b = 2$, \ie \ $q^a \equiv 1 \mod 3$, one can modify the equation defining $\Wrst$ \eqref{eq: Integ eq infty 2} to obtain a monic polynomial on the right-hand side. Indeed, letting $\varrho := \left( 1- s^{q^a -1}\right)^2 - 4s^{2q^a}$, $\widetilde{x} := 1/x_2$ and $\widetilde{y} := y_{2} / (x_{2} \sqrt{\varrho})$ (notice that $\varrho$ is a square in $\Oinf$ by Hensel's lemma), equation \eqref{eq: Integ eq infty 2} becomes 
	\begin{equation*}
		\widetilde{y}^2 = \left(\widetilde{x}^3 + \frac{s(1-s^{q^a -1})- 2s^{q^a + 1} }{\varrho}\right)\left(\widetilde{x}^3 + \frac{s(1-s^{q^a -1})+ 2s^{q^a + 1} }{\varrho}\right).
	\end{equation*}
	Thus, independently of $q^a \mod 3$, the defining equation of $\Wrst$ at infinity is of the form 
	\begin{equation}\label{eq: Def eqt infty}
		y^2 = \left(x^3 + su - s^{q^a+1} v \right) \left(x^3 + su + s^{q^a+1} v \right) \quad \text{ for some } u,v \in \Oinf\st.
	\end{equation} 
\end{remark}

In the following, we will manipulate regular models of $C$. We denote the proper minimal regular model of $C$ by $\Cmod \rightarrow \Pp$ (see \cite[§9.3.2 \& 9.3.3]{Liubook} for its definition and main properties). This is a regular surface fibered over $\Pp$ whose generic fiber is isomorphic to $C$, and that contains no exceptional divisors. Moreover, $\Cmod$ is minimal, in the sense that any birational map $\mathcal{C'} \dashrightarrow \Cmod$ of regular surfaces fibered over $\Pp$ is a birational morphism. We denote the Simple Normal Crossings (SNC) minimal regular model of $C$ by $\Kmod \rightarrow \Pp$ (see \cite[§9.3.4]{Liubook} for details). This is a regular surface fibered over $\Pp$, whose generic fiber is isomorphic to $C$, such that all its singularities are ordinary double points. Finally, we let $\Smod \rightarrow \Pp$ be the Néron model of $S$ (see \cite[§9.5]{BoschLutkebohmerRaynaud} for details).

For any $v \in \Mpl$, we denote by $\Cmod_{v}$ the base change of $\Cmod$ to $\Spec(\Fv)$, and we call this the \textit{special fiber} or \textit{reduction} of $\Cmod$ at $v$. We define in the same way the fiber at $v$ of the models $\Wrst$, $\Kmod$ and $\Smod$ introduced above. Since $\Spec \left(\unr{\Ov}\right) \rightarrow \Spec\left(\Ov\right)$ is étale, Proposition 10.1.17 from \cite{Liubook} states that $\Cmod \times_{\Ppo} \unr{\Ov}$ is the minimal regular model of $C \times_{\K} \unr{\K_{v}}$ (the same holds if one replaces $\Cmod$ by $\Kmod$). Similarly, Néron models commute to étale base change, so $\Smod \times_{\Ppo} \unr{\Ov}$ is the Néron model of $S \times_{\K} \unr{\K_{v}}$. We denote by $\Cmod_{\vbar}$ the special fiber of $\Cmod \times_{\Ppo} \unr{\Ov}$, and similarly for $\Kmod_{\vbar}$ and $\Smod_{\vbar}$.

\subsection[Special fibers of Cmod, Kmod and Smod]{Special fibers of $\Cmod$, $\Kmod$ and $\Smod$} \label{sect: Reduction}  

In this section, we describe the fibers of $\Cmod$, $\Kmod$ and $\Smod$ at various places of $\K$ in order to compute arithmetic invariants of $S$ in Section~\ref{sect: Invariants}. To do so, we  make use of the theory of cluster pictures: we refer the reader to \cite{D2M2} for details, or to \cite{Hyperuser} for a quick survey.

We begin with the minimal regular model $\Cmod$. We identify places of $\K$ with the corresponding monic irreducible polynomial in $\Fq[t]$. If $v \nmid \Delta(\Wrst)$, then the fiber $\Wrst_{v}$ is smooth and $\Wrst_{v} \simeq \Cmod_{v}$ are isomorphic over $\Fv$. We will refer to the places dividing $\Delta(\Wrst)$ as \textit{places of bad reduction} (as we will see later, for our specific models, $\Cmod$ has bad reduction at $v \in \Mpl$ if and only if $\Smod$ has bad reduction at $v$). For simplicity, we will write $\Ppo$ instead of $\Pp$ whenever the context is clear.

\begin{proposition} \label{prop: Fibers}
	The singular fibers of the minimal regular model $\Cmod$ are described as follows:
	\begin{itemize}
		\item At any finite place of bad reduction $v$, the special fibers of $\Cmod$ and $\Wrst$ are isomorphic over $\Fv$. On an affine subset, they are described by the hyperelliptic equation:
		\vspace{-0.25em}\begin{equation*}
			\Cmod_{v} \simeq \Wrst_{v} : y^2 = x^3 \left(x^3 + \rho_{v} \right) \quad \text{where } \rho_{v} :=  \begin{cases}
				\As - 2 \mod v & \text{ if } v \divi \As+2, \\
				\As + 2 \mod v & \text{ if } v \divi \As-2. \\
			\end{cases}
		\end{equation*}
		\item At infinity, the special fiber $\Cmod_{\infty}/ \Fq$ is a tree of $6q^a + 7$ $\Ppo$'s arranged as in Figure \ref{fig: Fiber Cmod infinity}. All components are defined over $\Fq$ and numbers indicate multiplicities.
		\begin{figure}[H]
			\centering
			
                \begin{tikzpicture}[
                    lfnt/.style={font=\small},
                    rightl/.style={above right=2.25em and 0pt,lfnt},
                    leftl/.style={above left= 2.25em and -1pt,lfnt},
                    mainl/.style={xshift = -0.7em, yshift = 0.65em}]
            
                    \draw (-4.75, 0) -- (-1.5, 0) node [above left = -0.5pt and 2.4em, lfnt] {$3$};
                    \draw (-3.25, 0.25) -- (-3.25, -2) node [leftl] {$1$};
                    \draw (-2, 0.2) -- (-0.75, -1) node [above left = 0.65em and 1.6em, lfnt] {$3$};
                    \draw (-4.5, 0.25) -- (-4.5, -2) node [leftl] {$2$};
                    \draw (-5.75, -1.75) -- (-4.25, -1.75) node [above left= -1.5pt and 2em, lfnt] {$1$};
            
                    \filldraw [black] (-0.35, -1) circle (1.25pt);
                    \filldraw [black] (0, -1) circle (1.25pt);  
                    \filldraw [black] (0.35, -1) circle (1.25pt);
                    
                    \draw (4.75, 0) -- (1.5, 0) node [above right = -0.5pt and 2.4em, lfnt] {$3$};
                    \draw (3.25, 0.25) -- (3.25, -2) node [leftl] {$1$};
                    \draw (2, 0.2) -- (0.75, -1) node [above right = 0.65em and 1.6em, lfnt] {$3$};
                    \draw (4.5, 0.25) -- (4.5, -2) node [leftl] {$2$};
                    \draw (5.75, -1.75) -- (4.25, -1.75) node [above right= -1.5pt and 2em, lfnt] {$1$};

                    \draw [decorate,decoration={brace,amplitude=8pt,mirror}]
                    (-1.3, -1.2) -- (1.3, -1.2) node[midway,yshift=-1.5em]{$6q^a -1$};
            
                \end{tikzpicture}

			\vspace{0.8em}
			\caption{Fiber $\Cmod_{\infty}$.}
			\label{fig: Fiber Cmod infinity}			
		\end{figure}
	\end{itemize}
\end{proposition}

\begin{proof}
	 Fix a finite place of bad reduction $v$. From Lemma \ref{lem: Wrst is minimal}, we know that $\Wrst$ is a minimal Weierstrass model at $v$. Proposition $6$ from \cite{Liu96} states that $\Cmod \times_{\Ppo} \Ov$ is the minimal desingularization of $\Wrst \times_{\Ppo} \Ov$. Since $\operatorname{disc}(x^3 + \As \pm 2) = -2^4 3^3 \left( \As \pm 2 \right)^2$, reducing equation \eqref{eq: Hyp model} $\mod v$ gives the following affine equation for the special fiber of $\Wrst$ at $v$:
	\begin{equation*}
		\Wrst_{v} : y^2 = x^3 \left(x^3 + \rho_{v} \right).
	\end{equation*}
	The projective curve $\Wrst_{v} / \Fv$ admits a singularity at $(0, 0)$ and is smooth at all other points. We claim that $\Wrst \times_{\Ppo} \Ov$ is regular at the point corresponding to the ideal $(x, y, \pi_{v}) \subset \Ov[x, y]$. Indeed, $\Wrst \times_{\Ppo} \Ov$, localized at this ideal is
	\begin{equation*}
		\Spec \left( \Ov \, [x, y]_{(x, y)} \, / \left(y^2 - \left(x^3 + \As - 2 \right) \left( x^3 + \As + 2 \right) \right) \right).
	\end{equation*}
	The defining polynomial is not in $\left(x, y, \pi_{v} \right)^2$ as the constant term $\As^2 - 4$ is square-free, so it has $v$-adic valuation $1$. Therefore, Corollary $4.2.12$ from \cite{Liubook} implies that $\Wrst \times_{\Ppo} \Ov$ is regular at the point corresponding to the ideal $(x, y, \pi_{v})$, hence everywhere.	We deduce the existence of an isomorphism $\Wrst \times_{\Ppo} \Ov \simeq \Cmod \times_{\Ppo} \Ov$: the latter is the minimal desingularization of the former, which is already regular. In particular the special fibers $\Wrst_{v}$ and $\Cmod_{v}$ are isomorphic over $\Fv$.
	
	At infinity, we know from Remark \ref{rmk: Defining eqt infinity} that we can take \eqref{eq: Def eqt infty} as a defining equation for $\Wrst$. This corresponds to Namikawa--Ueno type $\left[III_{6q^a}\right]$ (see \cite[p.184]{NamikawaUeno}): using their table, we see that the fiber of the minimal regular model $\Cmod_{\infbar}$ is as depicted in Figure \ref{fig: Fiber Cmod infinity}. Notice that all singularities in the fiber $\Cmod_{\infbar}$ are ordinary double points: hence $\Cmod_{\infbar}$ coincides with the fiber $\Kmod_{\infbar}$. Using the theory of cluster pictures, we are going to see that all irreducible components appearing in Figure \ref{fig: Fiber Cmod infinity} are actually defined over $\Fq$.
	
	The roots of the defining polynomial in \eqref{eq: Def eqt infty} are $\eta_{i}^{\pm} := \zeta_{3}^{i} \sqrt[3]{s (u \pm s^{q^a} v)}$, where $\zeta_{3} \in \Ksep$ is a third root of unity. It is easy to check that all roots have $s$-adic valuation $1/3$. Computing the Taylor expansion in $\Oinf$ of $\left(1 \pm s^{q^a} v /u \right)^{1/3}$ we see that $\ord_{s}\left(\eta_{i}^{+} - \eta_{i}^{-}\right) = 1/3 + q^a$ for any $0 \leq i \leq 2$, and $\ord_{s}\left(\eta_{i}^{+} - \eta_{j}^{-}\right) = 1/3$ if $i \neq j$. We deduce that the cluster picture of $C$ at $\infty$ is the following:
	\begin{figure}[H]
		\centering
		\includegraphics[width=5cm]{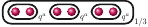}
		\caption{Cluster picture of $C$ at $\infty$.}
		\label{fig: Cluster infinity}
	\end{figure}
	\vspace{-0.5em} Using \cite{FaraggiNowell}, one can describe the special fiber $\Kmod_{\infbar} = \Cmod_{\infbar}$ and the action of Frobenius on irreducible components, in terms of the cluster picture at $\infty$ and the action of $\Gal(\K_{\infty}^{\mathrm{sep}} / \K_{\infty})$ on it. We freely adopt the notations of \cite{FaraggiNowell} until the end of the proof. The orbit $X_{\Rroots} := \lbrace \Rroots \rbrace$ contributes two central components with tails of $\Ppo$'s. Write $\twin_i := \left \lbrace \eta_{i}^{+}, \eta_{i}^{-} \right \rbrace$ for $0 \leq i \leq 2$ and let $X_{\twin} := \lbrace \twin_0, \twin_1, \twin_2 \rbrace$ be the orbit under Galois of $\twin_0$. Then $X_{\twin}$ contributes with a chain of $\Ppo$'s of length $6q^a -1$ linking the two central components associated to the orbit $X_{\Rroots}$. The quantities $\epsilon_{X_{\Rroots}}, \epsilon_{X_{\twin}} \in \lbrace \pm 1 \rbrace$ describe the action of the Frobenius on irreducible components of $\Cmod_{\infbar}$. It is easy to check that for any $i \in \lbrace 0, 1, 2 \rbrace$, we have $\Rroots^{*} = \twin_{i}^{*} = \Rroots$ which implies that $\epsilon_{X_{\Rroots}} = \epsilon_{X_{\twin}} = 1$. We then deduce from Theorem 7.21 of \cite{FaraggiNowell} that the Frobenius fixes each irreducible component of $\Cmod_{\infbar}$. Hence all these components are defined over $\Fq$, so the picture of $\Cmod_{\infbar}$ ``descends" to $\Fq$ and $\Cmod_{\infty}$ looks like Figure \ref{fig: Fiber Cmod infinity}.
\end{proof}

\begin{remark}
	The author fears that the description of the special fiber of the minimal regular SNC model given in \cite[Theorem 7.18]{FaraggiNowell} includes a typographical error. Indeed, applying the mentioned theorem, one obtains for the fiber of $\Kmod_{\infbar}$ a similar drawing as the one obtained in Figure \ref{fig: Fiber Cmod infinity}, with a subtle difference. The $\Ppo$'s of multiplicity $1$ branching the ``central" components of multiplicity $3$ are replaced by a chain of length $2$, with multiplicities $2$ and $1$ (identical to the ones already appearing in Figure \ref{fig: Fiber Cmod infinity}). The author has checked by other methods that $\Kmod_{\infbar} = \Cmod_{\infbar}$ is indeed as shown in Figure \ref{fig: Fiber Cmod infinity}, and that all irreducible components are indeed defined over $\Fq$.
\end{remark}

\begin{remark}
	The strategy from \cite{FaraggiNowell} can also be applied to describe $\Kmod_{\vbar}$ when $v$ is a finite place of bad reduction. At any such place, it is easy to check that the cluster picture of $C$ is as shown in Figure \ref{fig: Cluster finite v}. The second author of \textit{loc. cit.} has made an explicit table giving the correspondence between cluster pictures and fibers of SNC minimal regular models for genus $2$ curves. Our cluster appears in \cite[p.241]{Nowell} and gives the fiber $\Kmod_{\vbar}$ shown in Figure \ref{fig: Fiber Kmod finite v}. In it, numbers denote multiplicities and the label ``g1" points out that the component is a genus $1$ curve.
		\begin{figure}[H]
		\centering
		\begin{minipage}[t]{.5\textwidth}
			\centering
			\includegraphics[width=4cm]{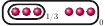}
			\caption{Cluster picture at $v|\Delta(\Wrst)$, $v \neq \infty$.}
			\label{fig: Cluster finite v}
		\end{minipage}%
		\begin{minipage}[t]{.5\textwidth}
			\centering
			
                \begin{tikzpicture}[scale = 0.8, lfnt/.style={font=\small},
                    rightl/.style={below right=0.2em and -4pt,lfnt},
                    leftl/.style={below left= 1.7em and -1pt,lfnt},
                    leftl2/.style={below left= 2.1em and 1.75em,lfnt},
                    mainl/.style={xshift = -0.7em, yshift = 0.55em}]]
                    \draw[domain=2.6953:3.55,samples=100,smooth,variable=\x] plot ({\x},{sqrt(10*(\x-3)^3 - 0.6*(\x-3) + 0.1)}) node[leftl2] {$1$} node[rightl] {g$1$};
                    \draw[domain=2.6953:3.55,samples=100,smooth,variable=\x] plot ({\x},{-sqrt(10*(\x-3)^3 - 0.6*(\x-3) + 0.1)});
                
                    \draw (-0.25, -1) -- (4.5, -1) node[mainl] {$6$};
                    \draw (0.25, -1.25) -- (0.25, 1.25) node[leftl] {$3$};
                    \draw (1.75, -1.25) -- (1.75, 1.25) node[leftl] {$2$};
                    \draw (2.6953, -0.01) -- (2.6953, 0.01);
            			
                \end{tikzpicture}

			\caption{Fiber $\Kmod_{\vbar}$ at $v|\Delta(\Wrst)$, $v \neq \infty$.}
			\label{fig: Fiber Kmod finite v}
		\end{minipage}
	\end{figure}
\end{remark}

We now turn our attention to describing the fibers of the  Néron model $\Smod \rightarrow \Pp$ of $S / \K$. For any smooth group scheme $G$ of finite type over $\Pp$ with geometrically connected generic fiber (e.g. $\Smod$), denote by $G^{0}$ the relative identity component. This is the open subscheme of $G$ obtained by removing the closed complement of the identity component in the finitely many disconnected fibers. Combining Chevalley's structure theorem for connected smooth algebraic groups (see \cite[Theorem 16]{Rosenlicht}) and the structure theorem for linear connected groups over perfect fields (see \cite[Exposé XVII]{SGA3II}) we get the following result. At any place $v \in \Mpl$,  $\Snerv$ fits in an exact sequence of $\Fv$-groups
\begin{equation*}
	 1 \rightarrow T \times U \rightarrow \Snerv \rightarrow A \rightarrow 1,
\end{equation*}
where $T$ is a torus, $U$ is a unipotent group and $A$ is an abelian variety. We refer the reader to \cite{Rosenlicht} and \cite[§9.1 - 9.5.]{BoschLutkebohmerRaynaud} for details.  The data of these three groups $A$, $T$ and $U$ is what we call the reduction type of $\Sner$ at $v$. We will talk about abelian, toric and unipotent dimension of $\Snerv$ to refer to the respective dimensions of the algebraic groups fitting in the exact sequence above. Whenever the unipotent group $U$ is trivial, \ie \ $\Sner_v$ is the extension of an abelian variety by a torus, we say that $\Smod$ has semistable reduction at $v$. \\

Let $\Pic_{\Cmod / \Ppo}$ be the Picard functor of $\Cmod / \Pp$ (\cf \ \cite{BoschLutkebohmerRaynaud}). Proposition \ref{prop: Fibers} implies that for any $v \in \Mpl$, the gcd of the multiplicities of the irreducible components of $\Cmod_{v}$ is $1$. Therefore, at any place $v \in \Mpl$, we have the following isomorphisms over $\Fv$ (see \cite[§9.5, §9.6]{BoschLutkebohmerRaynaud})
\begin{equation}\label{eq: Isom fibres modeles}
	\Snerv \simeq \left(\Pic_{\Cmod / \Ppo}^{0}\right)_{v} \simeq \Pic^{0}(\Cmod_v).
\end{equation}
Using this, we describe the reductions of $\Sner$. 

\begin{theorem}\label{thm: Red Smod}
The special fibers of the identity component of the Néron model $\Sner$ are as follows:
	\begin{itemize}
		\item At any $v \nmid \Delta(\Wrst)$, $\Cmod_{v}$ is smooth and we have $\Smod_{v} = \Sner_{v} = \Jac(\Cmod_{v})$.
		\item At any $v \divi \Delta(\Wrst)$, $v \neq \infty$, $\Snerv$ fits in the exact sequence of $\Fv$-groups
		\begin{equation*}
			1 \rightarrow \Ga \rightarrow \Snerv \rightarrow E_v \rightarrow 1,
		\end{equation*}
		where $E_v / \Fv$ is the elliptic curve given by the affine equation $E_v : y^2 = x(x^3 + \rho_{v})$, and $\rho_{v}$ is defined as in Proposition \ref{prop: Fibers}.
		\item At $\infty$, $\Sner$ has totally unipotent reduction, \ie, $\Sner_{\infty}$ is a unipotent group of dimension $2$.
	\end{itemize}
\end{theorem}

\begin{proof}
	For a finite place of bad reduction $v$, we use Section 7.5 in \cite{Liubook} to describe $\Pic^{0}(\Cmod_{v})$ in terms of the normalizations of the irreducible components of $\Cmod_{v}$. Recall that the special fiber $\Cmod_{v} \simeq \Wrst_{v}$ is reduced and irreducible. We compute its normalization by blowing up along the closed singular point $(0, 0)$. The result is the elliptic curve $E_v / \Fv$ presented in the statement of the theorem and given by
	\begin{equation}\label{eq: Def ell curve Ev}
		E_{v} : y^2 = x \left( x^3 + \rho_{v} \right).
	\end{equation}
	Theorem 7.5.19 from \cite{Liubook} states that the ``abelian part" of $\Snerv$ is the elliptic curve $E_v$. To compute the toric dimension of $\Snerv$, consider the normalization morphism 
	\begin{equation*}
		\map{\nu_v}{E_v}{\Wrst_{v}}{(x, y)}{(x, xy).}
	\end{equation*}
	Every point on $\Wrst_{v}$ has a unique preimage by $\nu_v$, even the singular point $(0, 0)$ as $\nu_{v}^{-1}((0, 0)) = \lbrace (0, 0) \rbrace$. The formula given in \cite[Proposition 7.5.18]{Liubook} implies that the toric dimension of $\Snerv$ equals $0$. Thus, the ``linear part" of $\Snerv$ is a unipotent group of dimension $1$ defined over the perfect field $\Fv$, which is therefore split and thus isomorphic to $\Ga$, giving the desired exact sequence. 
	
	At infinity, we compute first $\Sner_{\infbar}$ using \cite[Corollary 1.4]{Lorenzini90}. It states that the abelian dimension of $\Sner_{\infbar}$ equals the sum of the genus of the irreducible components of $\Cmod_{\infbar}$, and the toric dimension equals the first Betti number of the dual graph of $\Cmod_{\infbar}$. But all these irreducible components are $\Ppo$'s, and the dual graph of is a tree, which has trivial homology. We deduce that the abelian and toric dimension are both equal to $0$. Hence, $\Sner_{\infbar}$ is a unipotent group of dimension $2$, which is isomorphic to $\Sner_{\infty}$ over $\Fqbar$, hence the result.
\end{proof}

Recall that an abelian variety defined over $\K$ is said to be \textit{constant} if it is isomorphic over $\K$ to the base change to $\K$ of an abelian variety defined over $\Fq$. An abelian variety is said to be \textit{isotrivial} if it becomes constant after base change to a finite extension of $\K$. From the description of the cluster picture of $C$ at infinity we deduce:

\begin{proposition}
	The abelian surface $S$ is non-isotrivial.
\end{proposition}

\begin{proof}
	Let $F / \K$ be the extension where $S$ attains semistable reduction (see \cite[Exposé IX]{SGA7I}). The cluster picture of $C$ at $\infty$, shown in Figure \ref{fig: Cluster infinity}, is such that every cluster has at most $2$ odd children. Theorem 10.3 of \cite{D2M2} then states that the Néron model of $S_F$ has totally toric reduction at $\infty$. By Corollary 7.4.4 of \cite{BoschLutkebohmerRaynaud}, the same holds for $S_{F'}$, where $F' / F$ is any finite extension. But the base change of a constant variety remains constant, and its Néron model must have good reduction everywhere. Hence we conclude that $S$ cannot be isotrivial.
\end{proof}

\subsection[Arithmetic invariants of S]{Arithmetic invariants of $S$}\label{sect: Invariants} 
We now exploit the results from the previous section to estimate arithmetic invariants of $S$ such as its conductor, Tamagawa numbers and differential height. 

Fix, once and for all, a prime number $\ell \neq p$. Let $\TlS$ be the $\ell$-adic Tate module of $S$ and $\VlS := \TlS \otimes_{\Zl} \Ql$: this is a $4$-dimensional $\Ql$-vector space endowed with a continuous linear action of $\Gal(\sep{\K}/\K)$. At any place $v \in \Mpl$, denote the inertia subgroup at $v$ by $I_v$, and let $\VlS^{I_v}$ be the subspace fixed by $I_v$. Let $N(S) = \prod_{v \in \Mpl} v^{n_{v}(S)}$ be the conductor divisor of $S$, written in multiplicative notation. We refer the reader to \cite{Serre70} for details on the conductor of an abelian variety. By the Néron-Ogg-Shafarevich criterion we have $n_{v}(S) = 1$ at any place $v$ of good reduction.

\begin{corollary}\label{cor: Conductor}
	The conductor divisor of $S$ is given by 
	\begin{equation*}
		N(S) = \prod_{\substack{v \divi \Delta(\Wrst) \\ v \neq \infty}} (v)^2 \ (\infty)^4.
	\end{equation*}
	The degree of the conductor is $\deg N(S) = 4 \, q^a + 4$.
\end{corollary}

\begin{proof}
	Our assumption $\charact(\K) > 5$ guarantees that $\VlS$ is tamely ramified (see Corollary $2$ to Theorem $2$ in \cite[§3]{SerreTate}), so at any place $v \in \Mpl$ we have $n_v(S) := \dim V_{\ell}(S) - \dim V_{\ell}(S)^{I_v} $. Now Remark $1$ from \cite[§3]{SerreTate} states that $\dim V_{\ell}(S) - \dim V_{\ell}(S)^{I_v}$ is twice the unipotent dimension of $\Snerv$. Therefore at any finite place of bad reduction $v$, we have $n_{v}(S) = 2$ and at the place at infinity, $n_{\infty}(S) = 4$. Finally, we compute
	\begin{align*}
		\deg N(S) & = \sum_{v \divi \As - 2} 2 \, \dv +  \sum_{v \divi \As + 2} 2 \, \dv \  + \ 4  \\
		& = 2 \, \left(\deg (\As - 2) \, + \, \deg (\As + 2) \right) + 4= 4 \,q^a + 4.
	\end{align*}
\end{proof}
We turn now to the Tamagawa numbers of $S$. For $v \in \Mpl$, the local Tamagawa number of $S$ at $v$ is defined as $c_{v}(S) := \abs{\left( \Smod_{v} / \Sner_{v} \right) (\Fv)}$. 

\begin{proposition}\label{prop: Tamagawa}
	For any $v \neq \infty$, the local Tamagawa number $c_v(S)$ equals $1$. At infinity we have $c_{\infty}(S) \in \lbrace 1, 3, 9 \rbrace$.
\end{proposition}
	
\begin{proof}
	At any finite place of good reduction $v$, the special fiber $\Smod_{v}$ is connected, so we have $c_v(S) = 1$. Fix now a place of bad reduction $v$, we first compute $\cvbar := \abs{\left(\Smod_{\vbar}^{ } / \Smod_{\vbar}^{0} \right)(\Fvbar)}$. Since $\left(\Smod_{v}^{ } / \Snerv \right)(\Fv)$ is a subgroup of $\left(\Smod_{\vbar}^{ } / \Smod_{\vbar}^{0} \right)(\Fvbar)$, we have $c_{v}(S) \divi \cvbar$. The advantage of introducing this new quantity is that we can apply a result by Lorenzini to compute $\cvbar$ in terms of $\Kmod_{\vbar}$. We write the latter as a divisor $\Kmod_{\vbar} = \sum_{i = 1}^{m} r_{i} \Gamma_{i}$ and let $d_{i} := \sum_{j \neq i} \Gamma_{i} \cdot \Gamma_{j}$ (where $\Gamma_{i} \cdot \Gamma_{j}$ stands for the intersection number of $\Gamma_{i}$ and $\Gamma_{j}$). Since $\Sner_{\vbar}$ has toric dimension $0$, Corollary $1.5$ from \cite{Lorenzini90} yields
	\begin{equation}\label{eq: cvbar Lorenzini}
		\cvbar = \prod_{i=1}^{m} r_{i}^{d_{i}-2}.
	\end{equation}
	Using the description of the fibers of $\Kmod$ drawn in Figures \ref{fig: Fiber Cmod infinity} and \ref{fig: Fiber Kmod finite v}, we compute $\cvbar = 3^{-1} 2^{-1} 1^{-1} 6 = 1$ for any $v \divi \Delta(\Wrst), v \neq  \infty$ and 
	\vspace{-0.8em}\begin{equation*}
		c_{\overline{\infty}}(S)= 1^{-1} 2^{0}1^{-1} 3 \prod_{i = 1}^{6q^a -1} \left(3^{0} \right) \, 3 \cdot 1^{-1} 2^{0}1^{-1} = 9,\vspace{-0.5em}
	\end{equation*}
	hence the result.
\end{proof}

To conclude this section, we compute the differential height of $S$. Let $z : \Pp \rightarrow \Smod$ be the identity section of the Néron model $\Smod \rightarrow \Pp$ and $\Omega_{\Smod / \Ppo}^{2}$ be the sheaf of relative $2$-differentials of $\Smod$ over $\Pp$. The pullback $\omega_{S}:= z^{*}\Omega_{\Smod / \Ppo}^{2}$ is a line bundle on $\Pp$, and we define the differential height of $S$ as
\begin{equation*}
	h(S) := \deg(\omega_{S}).
\end{equation*}
We define also the exponential height of $S$ as $H(S) := q^{h(S)}$.

\begin{proposition}\label{prop: Height}
	The differential height of $S$ is given by $h(S) = q^a + 1$.
\end{proposition}

\begin{proof}
	Notice first that the polynomial in the right-hand side of \eqref{eq: Hyp model} is unitary and has even degree. Thus, $C$ has two points at infinity, which we denote by $\infty_1$, $\infty_2$; both are $\K$-rational. Each of them induces a section of $\xi : \Kmod \rightarrow \Pp$, whose image lies in the smooth locus of $\xi$. Let $\omega_{\Kmod / \Ppo}$ be the relative dualizing sheaf of $\Kmod \rightarrow \Pp$, then Proposition 7.4 of \cite{ExplJacLeg} gives the isomorphism of line bundles
	\begin{equation*}
		\omega_{S} \simeq \bigwedge^{2} \xi_{*} \omega_{\Kmod / \Ppo}.
	\end{equation*}
	Consider the two differentials $\mu_{1}:= \frac{dx}{y}$ and  $\mu_{2}:= \frac{x\,dx}{y}$ in $\H^{0}\left(\Kmod, \omega_{\Kmod / \Ppo}\right)$. The contraction of exceptional divisors $\Kmod \rightarrow \Cmod$ is a birational morphism of regular surfaces over $\Pp$. Corollary 9.2.25 in \cite{Liubook} then shows that the subgroups $\H^{0}\left(\Kmod, \omega_{\Kmod / \Ppo}\right)$ and $\H^{0}\left(\Cmod, \omega_{\Cmod / \Ppo}\right)$ of $\Omega_{C / \K}^{1}$ coincide. We thus may view $\mu_{1}, \mu_{2}$ as differentials on $\Cmod$ and work with $\omega_{\Cmod / \Ppo}$ rather than $\omega_{\Kmod / \Ppo}$. Moreover $(\mu_{1|C}, \mu_{2| C})$ forms a $\K$-basis of $\H^{0}\left(C, \Omega_{C / \K}^{1} \right)$, so $\xi_{*}(\mu_1 \wedge \mu_2)$ is a generator of $\omega_{S}$. We deduce
	\begin{equation*}
		\deg(\omega_{S}) = \deg(\xi_{*}(\mu_1 \wedge \mu_2)) = \sum_{v \in \Mpl} \ord_{v}(\xi_{*}(\mu_1 \wedge \mu_2)) \, \deg v.
	\end{equation*}
	At any place of good reduction $v$, the base change $\Cmod \times_{\Ppo} \Ov$ is smooth. Thus, after restriction, $(\mu_{1}, \mu_{2})$ forms an $\Ov$-basis of $\H^{0}\left(\Cmod \times \Ov, \omega_{\Cmod / \Ov}\right)$, implying that $\ord_{v}(\xi_{*}(\mu_1 \wedge \mu_2)) = 0$. At a finite place of bad reduction $v$, we know from Lemma \ref{lem: Wrst is minimal} that equation \eqref{eq: Hyp model} defining $C$ \eqref{eq: Hyp model} is minimal. Moreover, $\Cmod$ has Namikawa--Ueno type $\left[I_{0}-II-0\right]$ at $v$, which corresponds to Ogg's numerical type $[0]$. Propositions 4 and 7 in \cite{Liu94} then imply that the restriction of $(\mu_{1}, \mu_{2})$ is a basis of $\H^{0}\left(\Cmod \times \Ov, \omega_{\Cmod / \Ov}\right)$, hence $\ord_{v}(\xi_{*}(\mu_1 \wedge \mu_2)) = 0$. 
	
	It follows that $\deg(\omega_{S}) = \ord_{\infty}(\xi_{*}(\mu_1 \wedge \mu_2))$. To compute the contribution at $\infty$, let $b \in \lbrace 1, 2 \rbrace$ be such that $q^a + b \equiv 0 \mod 3$. We know from Lemma \ref{lem: Wrst is minimal} that, at $\infty$, the defining equation \eqref{eq: Integral eq infty} (resp. \eqref{eq: Integ eq infty 2}) if $b = 1$ (resp. $b = 2$) is minimal. The proof of Proposition \ref{prop: Fibers} shows that the Namikawa--Ueno type of $\Cmod_{\infty}$ is  $\left[ III_{6q^a}\right]$, which corresponds to Ogg's numerical type $[43]$. Write $(x_b, y_b)$ for the coordinates appearing in \eqref{eq: Integral eq infty} or \eqref{eq: Integ eq infty 2} depending on the value of $b$. Then Propositions 4 and 7 in \cite{Liu94} imply that $\left( \frac{\d x_{b}}{y_{b}}, \frac{x_{b} \, \d x_{b}}{y_{b}} \right)$ forms an $\Oinf$-basis of $\H^{0}\left(\Cmod \times \Oinf, \omega_{\Cmod / \Oinf}\right)$. Using the change of variables giving $(x, y)$ in terms of $(x_b, y_b)$ given in Definition \ref{def: Wrst model}, it is easy to check that, for both values of $b \in \lbrace 1, 2 \rbrace$, we have
	\begin{equation*}
		\mu_{1}\wedge \mu_{2} = s^{q^a + 1} \left( \frac{\d x_{b}}{y_{b}} \wedge \frac{x_{b} \, \d x_{b}}{y_{b}} \right), 
	\end{equation*}
	which finally gives the result.
\end{proof}

\vspace{-1.25em}


\section[Character sums and the geometry of Cmod]{Character sums and the geometry of $\Cmod$}\label{sect: Characters and curves}

As we saw in Section~\ref{sect: Constructing ab surf}, the minimal regular model $\Cmod$ of $C$ is dominated by a product of curves. In this section we describe the geometry of $\Cmod$ in further detail in terms of the two auxiliary curves $X$ and $Y$ defined by \eqref{eq: Models X Y} (to ease notation we drop the subscript $a$). We will study the cohomology of $X$ and $Y$: the eigenvalues of the action of the Frobenius on these spaces may be expressed in terms of certain character sums over finite fields. We begin by recalling some facts about the relevant sums.

\subsection{Background on Gauss and Kloosterman sums}\label{sect: Character sums} 
In this section, we let $p \geq 3$ be a prime number and temporarily drop the assumption $p \geq 7$. We fix a finite field $\Fq$ of characteristic $p$.
~\\

We first discuss additive and multiplicative characters on finite fields. Our main reference here is \cite[Chapter 5, §1]{LidlNiederreiter}. Fix, once and for all, a non-trivial additive character $\psi_0 : \Fp \rightarrow \Q(\zeta_p)\st$ with values in the $p$-th cyclotomic field. Let $\Ff$ be a finite field of characteristic $p$ and write $\Tr{\Ff}{\Fp}: \Ff \rightarrow \Fp$ for the relative trace map. We recall that all additive characters on $\Ff$ are obtained in the following way (see \cite[Theorem 5.7]{LidlNiederreiter}): for any $\alpha \in \Ff$, define the map 
\begin{equation*}
	\psi_{\Ff, \, \alpha}: \Ff \rightarrow \Q(\zeta_p)\st, \quad x \mapsto \psi_0 (\Tr{\Ff}{\Fp}(\alpha x)).
\end{equation*}
One easily checks that $\psi_{\Ff, \, \alpha}$ is an additive character on $\Ff$, and that it is non-trivial if and only if $\alpha \neq 0$. To ease notation, we write $\psi_{\Ff}$ instead of $\psi_{\Ff, \, 1}$. Recall that the group $\widehat{\Ff}$ of additive characters on $\Ff$ is isomorphic to $\Ff$ itself \textit{via} $\alpha \mapsto \psi_{\Ff, \, \alpha}$.

We now discuss multiplicative characters. Fix an algebraic closure $\Qbar$ of $\Q$, and denote the ring of algebraic integers by $\Zbar$. We fix a prime ideal $\pfrac$ of $\Zbar$ lying over the rational prime $p$ and let $\valp : \overline{\Q} \rightarrow \Q$ be the $\pfrac$-adic valuation on $\Qbar$, normalised so that $\valp(q) = 1$. The quotient $\Zbar/\pfrac$ is an algebraic closure of $\Fp$, denoted by $\Fpbar$. The reduction map $\Zbar \rightarrow \Zbar / \pfrac$ induces an isomorphism between the group of roots of unity in $\Qbar$ whose order is prime to $p$ and $\Fpbar\st$. Let $\boldsymbol{\chi}: \Fpbar\st \rightarrow \Qbar\st$ denote the inverse of this isomorphism. For any positive integer $n$ dividing $\abs{\Ff\st}$, define the map
\begin{equation*}
	\chi_{\Ff, n}^{ }: \Ff\st \rightarrow \Qbar\st, \quad x \mapsto \boldsymbol{\chi}(x)^{|\Ff\st| / n}.
\end{equation*}
This map is a multiplicative character on $\Ff\st$ of exact order $n$ (see \cite[§4.3.2]{Cohen}). If $\Fp \subset \Ff \subset \Ff'$ are successive finite extensions and $\Nrm{\Ff'}{\Ff} : (\Ff')\st \rightarrow \Ff\st$ denotes the relative norm, we have the equality ${\chi_{\Ff', n}^{ } = \chi_{\Ff, n}^{ } \circ \Nrm{\Ff'}{\Ff}}$.
\begin{definition}
	Let $\Ff / \Fp$ be a finite field extension, $\psi : \Ff \rightarrow \Q(\zeta_p)\st$ be an additive character on $\Ff$ and $\chi : \Ff\st \rightarrow \Qbar\st$ a multiplicative character on $\Ff\st$. Define the Gauss sum 
	\begin{equation*}
		\Gauss{\Ff}{\chi}{\psi} := - \sum_{x \in \Ff\st} \chi(x) \psi(x).
	\end{equation*}
\end{definition}

If $\chi$ has order $n$, without loss of generality we can assume that $\chi$ takes values in the cyclotomic field $\Q(\zeta_{n})$. Therefore the sum $\Gauss{\Ff}{\chi}{\psi}$ is an algebraic integer in $\Q\left(\zeta_{n p}\right)$. For any additive character $\psi$ on $\Ff$ and any multiplicative character $\chi$ on $\Ff\st$, we have the following facts:
\begin{enumerate}[(\sf {Ga}\,1)]
	\item \label{Gauss norm} If $\psi$ and $\chi$ are non-trivial then, in any complex embedding of $\Q(\zeta_{n p})$, $\abs{\Gauss{\Ff}{\chi}{\psi}} = \abs{\Ff}^{1/2}$.

	\item\label{Extract Gauss} For any $\alpha \in \Ff\st$, $\Gauss{\Ff}{\chi}{\psi_{\Ff, \, \alpha}} = \chi(\alpha)^{-1} \Gauss{\Ff}{\chi}{\psi_{\Ff}}$.
	
	\item\label{Gauss bien def} For any finite extension $\Ff' / \Ff$, any $\alpha \in \Ff^{\times}$,
	\begin{equation*}
		\Gauss{\Ff'}{\chi^{\abs{\Ff}} \circ \Nrm{\Ff'}{\Ff}^{ }}{\ \psi_{\Ff', \, \alpha^{1/ \abs{\Ff}}}} = \Gauss{\Ff'}{\chi \circ \Nrm{\Ff'}{\Ff}}{\psi_{\Ff', \, \alpha}}.
	\end{equation*}

	\item\label{Hasse Davenport} (Hasse--Davenport) For any finite extension $\Ff' / \Ff$, 
	\begin{equation*}
		\Gauss{\Ff'}{\chi \circ \Nrm{\Ff'}{\Ff}}{\psi \circ \Tr{\Ff'}{\Ff}} = \Gauss{\Ff}{\chi}{\psi}^{[\Ff' : \Ff]}.
	\end{equation*}
	
	\item\label{Stickelberger Griffon} If $\chi$ has order $n$,
	\begin{equation*}
		\frac{1}{n} \, \leq \,  \frac{\valp(\Gauss{\Ff}{\chi}{\psi})}{[\Ff : \Fq]} \, \leq \, 1 - \frac{1}{n}.
	\end{equation*}
\end{enumerate}
These are well-known facts: the reader can find proofs for~\refGauss{Gauss norm},~\refGauss{Extract Gauss}, and~\refGauss{Hasse Davenport} in \cite[Chap. V, §2]{LidlNiederreiter}, ~\refGauss{Gauss bien def} is \cite[Lemma 2.5.8]{Cohen}, and~\refGauss{Stickelberger Griffon} follows from Stickelberger theorem: see \cite[Corollary 6.11]{AGTT} with notation adapted to our setting.

\begin{definition}
	Let $\Ff / \Fp$ be a finite field extension, $\psi : \Ff \rightarrow \Q(\zeta_p)\st$ be an additive character, and $\alpha \in \Ff$. Define the Kloosterman sum 
	\begin{equation*}
		\Kloos{\Ff}{\psi}{\alpha} := - \sum_{x \in \Ff\st} \psi\left(x+\frac{\alpha}{x}\right).
	\end{equation*}
\end{definition}
It is well-known that the sum $\Kloos{\Ff}{\psi}{\alpha}$ is a totally real algebraic integer in $\Q(\zeta_p)$. For any non-trivial additive character $\psi$ on $\Ff$, we have the following facts:
\begin{enumerate}[(\sf {Kl}\,1)]
	\item\label{Kloos bien def} For any finite extension $\Ff' / \Ff$ and any $\alpha \in \Ff$, $\Kloos{\Ff'}{\psi \circ \Tr{\Ff'}{\Ff}}{\alpha^{\abs{\Ff}}} = \Kloos{\Ff'}{\psi \circ \Tr{\Ff'}{\Ff}}{\alpha}$.
 
	\item\label{Kl HasseDav} For any $\alpha \in \Ff\st$, there is a unique pair $\left\{\kl{\Ff}{\psi}{\alpha}, \klbis{\Ff}{\psi}{\alpha} \right\}$ of algebraic integers, whose product is $\abs{\Ff}$ and such that, for any finite extension $\Ff' / \Ff$,
	\begin{equation*}
		\Kloos{\Ff'}{\psi \circ \Tr{\Ff'}{\Ff}}{\alpha} = 	\kl{\Ff}{\psi}{\alpha}^{[\Ff' : \Ff]} + \klbis{\Ff}{\psi}{\alpha}^{[\Ff' : \Ff]}.
	\end{equation*}
		
	\item\label{Kloos norm} For any $\alpha \in \Ff\st$ one has $\abs{\kl{\Ff}{\psi}{\alpha}} = \abs{\klbis{\Ff}{\psi}{\alpha}} = \abs{\Ff}^{1/2}$ in any complex embedding. In particular $\abs{\Kloos{\Ff}{\psi}{\alpha}} \leq 2 \abs{\Ff}^{1/2}$ and $\kl{\Ff}{\psi}{\alpha}$, $\klbis{\Ff}{\psi}{\alpha}$, seen as complex numbers, are conjugate.
	
	\item\label{Kloos conj valuation} If $\Ff$ is a finite extension of $\Fq$, then for any $\alpha \in \Ff\st$, the $\pfrac$-adic valuations of $\kl{\Ff}{\psi}{\alpha}$ and $\klbis{\Ff}{\psi}{\alpha}$ are 
	\begin{equation*}
		 \left \lbrace \valp(\kl{\Ff}{\psi}{\alpha}), \, \valp(\klbis{\Ff}{\psi}{\alpha})  \right \rbrace = \lbrace 0, [\Ff : \Fq] \rbrace.
	\end{equation*}	
\end{enumerate}
 
The first point is straightforward;~\refKloos{Kl HasseDav} and~\refKloos{Kloos norm} are Theorems 5.43 and 5.44 in \cite{LidlNiederreiter}. Finally, a proof of~\refKloos{Kloos conj valuation} can be found in \cite[Lemma 5.1]{GriDWit}.

\subsection{Orbits}\label{sect: Orbits} 
Assume once again that $p \geq 7$. Fix an integer $a \geq 1$, consider the action of $\langle q \rangle \subset \Q\st$ on $\Fqa\st$ \textit{via} $q \cdot \alpha = \alpha^{1/q}$ for any $\alpha \in \Fqa$. Write $\Pqa$ for the set of orbits of $\Fqa\st$ under this action. $\Pqa$ is in bijection with the set of closed points of the multiplicative group $\Gm$ over $\Fq$ of degree dividing $a$. Indeed, for any $v \in \Pqa$, any $\beta \in v$, we have $v = \left \lbrace \beta, \beta^{q}, \beta^{q^2}, \, \ldots, \, \beta^{q^{\abs{v} - 1}} \right \rbrace$, and all elements in this set induce the same closed point of $\Gm$. Under the identification between closed points of $\Pp$ and places of $\K$, $\Pqa$ corresponds to the set of places $\neq 0, \infty$ whose degree divides $a$. Given an orbit $v \in \Pqa$, its size is $\abs{v} = \dv = [ \Fv : \Fq]$. \\

From now on, we will only need to manipulate multiplicative characters of order $3$. Endow the set $\Ia := \left(\Z / 3\Z \right)\st \times \Fqa\st$ with the action of $\langle q \rangle$ defined by $q \cdot (j, \alpha) = (q \, j \mod 3,\, \alpha^{1/q})$ for any $(j, \alpha) \in \Ia$. Define $\Oqa$ to be the set of orbits under this action. There is a well-defined map $\Theta : \Oqa \rightarrow \Pqa$ given as follows: for any $o \in \Oqa$ and any representative $(j, \alpha) \in o$, we let $\Theta(o)$ be the orbit of $\alpha$ in $\Pqa$. Moreover, this map $\Theta$ is surjective. The size of an orbit in $\Oqa$ can be expressed as follows. For any $o \in \Oqa$, let $v := \Theta(o) \in \Pqa$, then we have 
\begin{equation}\label{eq: Size orbit eqt}
	\abs{o} =  \lcm \left(\ord\st(q \mod 3), \abs{v} \right). 
\end{equation}
Indeed, if $(j, \alpha) \in o$ then $\abs{o}$ is the smallest integer $n$ such that $q^n \, j \equiv j \mod 3$ and $\Fq(\alpha) = \Fv \subset \Ff_{q^{n}}$. We can give a more detailed description of the orbits in $\Oqa$.
\begin{itemize}
	\item When $q \equiv 1 \mod 3$, we have $\ord\st(q \mod 3) = 1$ and  $\langle q \rangle$ acts trivially on $\left(\Z / 3\Z \right)\st$. Hence, in that case, the orbits in $\Oqa$ are vertical, in the sense that for any $o \in \Oqa$, the ``first coordinate" $j \in \left(\Z / 3\Z \right)\st$ is the same for all representatives $(j, \alpha) \in o$. There is a map $\Oqa \rightarrow \left( \Z / 3\Z \right)\st$, sending the orbit $o \in \Oqa$ to $j \in \left( \Z / 3\Z \right)\st$, where $(j, \alpha) \in o$ is any representative. We choose $1, -1 \in \Z$ as representatives of $\left( \Z / 3\Z \right)\st$ and lift the map above to  $\prun : \Oqa \rightarrow \lbrace \pm 1 \rbrace$. Moreover, $\Theta$ is $2$-to-$1$ and \eqref{eq: Size orbit eqt} implies that $\abs{o} = \abs{\Theta(o)}$ for any $o \in \Oqa$.
	
	\item When $q \equiv 2 \mod 3$, the orbits in $\Oqa$ ``zigzag". This means that, for any $o \in \Oqa$ the first coordinate of the representatives of $o$ is iteratively exchanged each time $q$ acts on one of them: if $(j, \alpha) \in o$, then $q \cdot (j, \alpha) = (- j, \alpha^{1/q})$. Letting $v := \Theta(o) \in \Pqa$, the size of $o$ depends on the parity of $\abs{v}$. If $\alpha \in v$, then the orbits of $(1, \alpha)$ and $(-1, \alpha)$ are equal if and only if $\abs{o} = \abs{v}$. This implies that 
	\begin{equation}\label{eq: Preimages by Theta}
		\abs{\Theta^{-1}(v)} = \begin{cases}
			2 & \text{ if } \abs{o} = 2 \abs{v}, \\
			1 & \text{ if } \abs{o} = \abs{v}.
		\end{cases} 
	\end{equation}
\end{itemize}
We let $\Fo$ be the extension of $\Fq$ of degree $\abs{o}$: if $(j, \alpha) \in o$ is any representative, then $\Fo$ is the smallest extension of $\Fq$ containing $\alpha$ and having a non-trivial multiplicative character of order $3$. 

For any $d \geq 1$, define $\Pi_{q}(d)$ to be the number of orbits $v \in \Pqa$ having size $\abs{v} = d$. 

\begin{lemma}\label{lem: Estimates numb orbits}
	For any $d, a \geq 1$, we have the following estimates:
	\begin{multicols}{3}
		\begin{enumerate}
			\item\label{item: i lem estimates numb orbits} $\abs{\Pi_q(d) - \frac{q^d}{d}} \leq \frac{q^{d/2}}{1-q^{-1}}.$
			\item $\frac{q^a}{a} \ll_q \abs{\Pqa} \ll_q \frac{q^a}{a}$.
			\item $\frac{q^a}{a} \ll_q \abs{\Oqa} \ll_q \frac{q^a}{a}$.
		\end{enumerate}
	\end{multicols}
\end{lemma}

\begin{proof}
	Item $i)$ is the Prime Number Theorem for $\Fqt$ (see \cite[Proposition 6.3]{Brumer} for details). We deduce $ii)$ from equality $\abs{\Pqa} = \sum_{d \divi a} \Pi_{q}(d)$ . To prove $iii)$, notice that $\Theta$ is surjective and for any $v \in \Pqa$, $\abs{\Theta^{-1}(v)} = 1$ or $2$ as stated above. We deduce that $\abs{\Pqa} \leq \abs{\Oqa}\leq 2 \abs{\Pqa}$, which combined with $ii)$ gives the desired result.
\end{proof}

Let $o \in \Oqa$, applying~\refGauss{Gauss bien def} to the finite extension $\Fo / \Fq$ shows that for any $(j, \alpha) \in o$, we have the equality $\Gauss{\Fo}{\chi_{\Fo, 3}^{q \, j}}{\psi_{\Fo, \, \alpha^{1/q}}} = \Gauss{\Fo}{\chi_{\Fo, 3}^{j}}{\psi_{\Fo, \, \alpha}}$. In the same way,~\refKloos{Kloos bien def} shows that for any $v \in \Pqa$ and any $\beta \in v$, we have $\Kloos{\Fv}{\psi_{\Fv, \, \beta^{1/q}}}{1} = \Kloos{\Fv}{\psi_{\Fv, \, \beta}}{1}$. Thus, as $(j, \alpha)$ runs through $o$, the Gauss sums $\Gauss{\Fo}{\chi_{\Fo, 3}^{j}}{\psi_{\Fo, \, \alpha}}$ remains constant. In the same way, as $\beta$ runs through $v$, $\Kloos{\Fv}{\psi_{\Fv, \, \beta}}{1}$ remains constant too. Thus, it is natural to define:

\begin{definition}\label{def: sums ass orbits}
	Let $v \in \Pqa$, $\beta \in v$. We define $\kapv := \kl{\Fv}{\psi_{\Fv, \beta}}{1}$, $\kapvbis := \klbis{\Fv}{\psi_{\Fv, \beta}}{1}$. \\
	Let $o \in \Oqa$, fix $(j, \alpha) \in o$. We define $\gamo := \Gauss{\Fo}{\chi_{\Fo, 3}^{j}}{\psi_{\Fo, \, \alpha}}$, $\Klo := \Kloos{\Fo}{\psi_{\Fo, \alpha}}{1}$, $\kapo := \kl{\Fo}{\psi_{\Fo, \alpha}}{1}$, and $\kapobis := \klbis{\Fo}{\psi_{\Fo, \alpha}}{1}$.
\end{definition}

With this notation, \refKloos{Kloos conj valuation} states that $\lbrace \valp(\kapo), \valp(\kapobis) \rbrace = \lbrace 0, \abs{o} \rbrace$. We can also compute the $\pfrac$-adic valuations of $\gamo$:

\begin{lemma}\label{lem: Valuations gauss}
	If $q \equiv 1 \mod 3$, then for any $o \in \Oqa$, 
	\begin{equation*}
		\valp (\gamo) =
		\begin{cases}
			2 \abs{o}/3  &\text{ if } \prun(o) = 1,\\
			\abs{o}/3  &\text{ if } \prun(o) = -1.\\
		\end{cases}
	\end{equation*}
	If $q \equiv 2 \mod 3$, then $\valp(\gamo) = \abs{o}/2$.
\end{lemma}

\begin{proof}
	The case $q \equiv 1 \mod 3$ is a consequence of Stickelberger's Theorem (see \cite[Theorem 6.10]{AGTT}). The case $q \equiv 2 \mod 3$ is a consequence of a result of Tate and Shafarevich, which states that, under this hypothesis, $\gamo = \pm \chi_{\Fo, 3}^{-j}(\alpha)(-q)^{\abs{o}/2}$ (see \cite[Lemma 6.7]{AGTT}). 
\end{proof}

\subsection{Cohomology of auxiliary curves}\label{sect: Auxiliary curves} 
In constructing the curve $C / \K$ (Section~\ref{sect: Constructing ab surf}), we introduced two smooth projective curves $X$ and $Y$ defined over $\Fq$ (for simplicity, we drop the subscript $a$). We now describe the étale cohomology groups of these curves (see \cite{Milne80}). After a more refined analysis of the geometry of $\Cmod$ (\cf \ Section~\ref{sect: Geom Cmod}), this will help us compute the $L$-function of the abelian surface $S$. \\

We fix a prime number $\ell \neq p$ and an isomorphism $\Qlbar \simeq \C$. The main reference for this section is \cite{Milne80}. For any variety $V$ defined over $\Fq$, we denote by $\H^{i}(V) := \Het(V_{\Fqbar}, \Ql) \otimes_{\Ql} \Qlbar$ the $i$-th $\ell$-adic cohomology groups of $V$ (see \cite[Chapter III]{Milne80}). These are finite dimensional $\Qlbar$-vector spaces that come with a natural action of the $q$-th power geometric Frobenius $\Fr$ (we write in the same way the endomorphism of $V_{\Fqbar}$ and the endomorphism on cohomology induced by functoriality).\\

Recall that the defining equations of $X$ and $Y$ are
\begin{equation}\label{eq: Hyp model recall}
	X : u^3 = \ArtSch(t_1) \quad \text{and} \quad Y: v + \frac{1}{v}= \ArtSch(t_2).
\end{equation}
Since both $X$ and $Y$ are smooth projective curves, $\H^{0}(X)$ and $\H^{0}(Y)$ are one-dimensional $\Qlbar$-vector spaces on which $\Fr$ acts trivially (see \cite[Chapter V, Theorem 2.5]{Milne80}). By Poincaré duality, $\H^{2}(X)$ and $\H^{2}(Y)$ are also one-dimensional and the eigenvalue of $\Fr$ on these lines is $q$. We are interested in $\H^{1}(X)$ and $\H^{1}(Y)$. Applying the Riemann--Hurwitz formula (see \cite[Theorem 7.4.16]{Liubook}) we can check that both $X$ and $Y$ have genus $q^a -1$, so $\H^{1}(X)$ and $\H^{1}(Y)$ have dimension $2(q^a -1)$. It is clear that the curve $X$ admits a unique point at infinity, denoted by $\infty_{X}$, which is $\Fq$-rational. On the other hand, $Y$ admits two points at infinity, which we denote by $\infty_{Y_{1}}, \infty_{Y_{2}}$, that are $\Fq$-rational. Indeed, $Y$ is hyperelliptic, described in affine coordinates by the hyperelliptic equation $y^2 = \As^2 - 4$, where the change of variables is explicitly given by $y = 2v - \As$. 

We first describe $\H^{1}(Y)$ and will then move on to $\H^{1}(X)$. The group $\Fqa$ acts on $Y$ as follows: on the affine subset of $Y$ defined by equation \eqref{eq: Hyp model recall} we let $\beta \cdot (v, t_2) := (v ,t_2 + \beta)$ for $\beta \in \Fqa$ and $(v, t_2) \in Y$. This action extends to the whole projective curve $Y$. One can check that the only fixed points of this action on $Y$ are $\infty_{Y_{i}}$, $i=1, 2$. Note, however, that this action is not necessarily defined over $\Fq$ since it does not commute with the action of $\Fr$ as soon as $a > 1$: we have $\Fr (\beta \cdot (v, t_2)) = \beta^q \cdot \Fr (v, t_2)$. 

By functoriality, $\H^{1}(Y)$ also carries an action of $\Fqa$. We can decompose this cohomology space into isotypical components indexed by irreducible representations of $\Fqa$. Since $\Fqa$ is an abelian group, its irreducible representations are its characters. The latter are described by the isomorphism $\Fqa \stackrel{\sim}{\longrightarrow} \widehat{\Fqa}, \beta \mapsto \psi_{\Fqa, \, \beta}$. Denoting by $\H^{1}(Y)^{\beta}$ the subset of $\H^{1}(Y)$ where $\Fqa$ acts by multiplication by $\psi_{\Fqa, \, \beta}$, we obtain the decomposition
\begin{equation}\label{eq: Decomp isotyp Y}
	\H^{1}(Y) = \bigoplus_{\beta \in \Fqa} \H^{1}(Y)^{\beta}.
\end{equation}
Let us notice that $\H^{1}(Y)^{0}$ is exactly the subspace of fixed points by $\Fqa$. Recall that for any $v \in \Pqa$, any representative $\beta \in v$ we have $v = \left \lbrace \beta, \beta^{q}, \beta^{q^2}, \, \ldots, \, \beta^{q^{\abs{v} - 1}} \right \rbrace$, and $\abs{v} = \dv = [\Fv : \Fq]$.

\begin{lemma}\label{lem: Eigenvalues Kloos}
	For any $\beta \neq 0$, $\H^{1}(Y)^{\beta}$ is $2$-dimensional. When $\beta = 0$, $\H^{1}(Y)^{0} = 0$. Moreover, for any $v \in \Pqa$, any $\beta \in v$, $\Fr^{\abs{v}}$ stabilizes $\H^{1}(Y)^{\beta}$, and its $2$ eigenvalues are $\kapv$ and $\kapvbis$.
\end{lemma}

\begin{proof}
	The description of how $\Fr$, and the action of $\Fqa$ interact on $Y$ and a short computation show that for any $\beta \in \Fqa$, $\Fr \left(\H^{1}(Y)^{\beta} \right) = \H^{1}(Y)^{\beta^{1/q}}$. Thus, for any $v \in \Pqa$, the Frobenius $\Fr$ cyclically permutes the subspaces $\H^{1}(Y)^{\beta}$ for $\beta \in v$. Therefore, any $\H^{1}(Y)^{\beta}$ is stabilized by $\Fr^{\abs{v}}$. Following \cite[§II]{Kat81}, and using arguments of \cite[§2]{vdGeervdVlugt}, one can show that, for any $\beta \neq 0$
	\begin{equation*}
		\operatorname{det} \left( 1 - \Fr^{\abs{v}} T \left| \, \H^{1}(Y)^{\beta} \right. \right) = 1 + \Kloos{\Fv}{\psi_{\Fv, \, \beta}}{1} T + q^{\dv} \, T^2.
	\end{equation*}
	In particular $\dim \H^{1}(Y)^{\beta} = 2$. Counting dimensions we deduce that $\dim \H^{1}(Y)^{0} = 0$. The two eigenvalues of $\Fr^{\abs{v}}$ on $\H^{1}(Y)^{\beta}$ are algebraic integers whose product equals $q^{\dv}$ and their sum equals $\Kloos{\Fv}{\psi_{\Fv, \, \beta}}{1}$. Thanks to the Riemann hypothesis (see \cite{Deligne80}), both of these algebraic integers have absolute value $\sqrt{q^{\dv}}$ in any complex embedding. These eigenvalues clearly satisfy relation~\refKloos{Kl HasseDav}, so they must be equal to $\kapv$ and $\kapvbis$ by uniqueness of the latter.
\end{proof}

Using \cite[Lemma 2.2]{Ulmer07} one computes, for any $v \in \Pqa$ and any representative $\beta_{v} \in v$:
\begin{equation*}
	\chpolg{\bigoplus_{\beta \in v} \H^{1}(Y)^{\beta}}  = \operatorname{det} \left( 1 - \Fr^{\abs{v}} T^{\abs{v}} \left| \, \H^{1}(y)^{\beta_v} \right. \right) 
	 = \left( 1 - \kapv \, T^{\abs{v}}\right) \left( 1 - \kapvbis \,  T^{\abs{v}}\right).
\end{equation*}
Taking the product over the orbits we deduce
\begin{equation}\label{eq: Charpol H1Y}
	\chpolg{\H^{1}(Y)} = \prod_{v \in P_{q}(a)} \left( 1 - \kapv \, T^{\abs{v}}\right) \left( 1 - \kapvbis \,  T^{\abs{v}}\right).
\end{equation}

We turn now to the curve $X / \Fq$ Denote by $\mu_{3}(\Fqbar)$ the set of third roots of unity in $\Fqbar$. The group $\GXa := \mu_{3}(\Fqbar) \times \Fqa$ acts on $X$: on the affine open subset of $X$ defined by equation \eqref{eq: Hyp model recall}, we set $(\zeta, \alpha) \cdot (u, t_1) = (\zeta u, t_1 + \alpha)$ for any $(\zeta, \alpha) \in \GXa$, and $(u, t_1) \in X$. This action can be extended to the whole projective curve $X$, but it is not necessarily defined over $\Fq$: one checks that $\Fr((\zeta, \alpha) \cdot (u, t_1)) = (\zeta^q, \alpha^q) \cdot \Fr(u, t_1)$ for any $(u, t_1) \in X, (\zeta, \alpha) \in \mu_{3}(\Fqbar) \times \Fqa$. The only fixed point of $X$ under this action is $\infty_{X}$. By functoriality, $\GXa$ acts on $\H^{1}(X)$ and, again, we decompose this cohomology space into isotypical components:
\begin{equation}\label{eq: Decomp isotyp X}
	\H^{1}(X) = \bigoplus_{\chi \in \widehat{\GXa}} \H^{1}(X)^{\chi}.
\end{equation}
Here $\chi$ runs through the characters of $\GXa$ and $\H^{1}(X)^{\chi}$ is the subspace of $\H^{1}(X)$ on which $\GXa$ acts by multiplication by $\chi$. The group $\widehat{\GXa}$ is isomorphic to $\Z / 3 \Z \times \Fqa$. Indeed, fix a non-trivial element $\zeta \in \mu_{3}(\Fqbar)$, the map that sends $(j, \alpha) \mapsto \left(\zeta \mapsto (\e^{2 i \pi /3})^j, \psi_{\Fqa, \, \alpha} \right)$ establishes the desired isomorphism. We denote by $\H^{1}(X)^{(j, \alpha)}$ the isotypical component corresponding to the character associated to $(j, \alpha)$. Recall the notation $\Ia = \left(\Z / 3\Z \right)\st \times \Fqa\st$ introduced above, and note that, for any $(j, \alpha) \in \Ia$, the character associated to $(j, \alpha)$ is non-trivial.

\begin{lemma}\label{lem: Eigenvalues gamma}
	For any $(j, \alpha) \in \Ia$ the subspace $\H^{1}(X)^{(j, \alpha)}$ is one-dimensional, and $\H^{1}(X)^{(j, \alpha)} = 0$ if $(j, \alpha) \in \left(\Z / 3\Z \times \Fqa \right) \setminus \Ia$.  For any $o \in \Oqa$, any $(j, \alpha) \in o$, $\Fr^{\abs{o}}$ stabilizes the line $\H^{1}(X)^{(j, \alpha)}$  and its eigenvalue is the Gauss sum $\gamo$. 
\end{lemma}

\begin{proof}
	This is a generalisation of \cite[Corollary 2.2]{Kat81} where Katz describes, when $q = p$, the cohomology of the curve $X$ and the eigenvalues of the Frobenius.
\end{proof}

Using \cite[Lemma 2.2]{Ulmer07} again, we have, for any $o \in \Oqa$ and any representative $(j_o, \alpha_{o}) \in o$:
\begin{equation*}
	\chpolg{\bigoplus_{(j, \alpha) \in o} \H^{1}(X)^{(j, \alpha)}}  = \operatorname{det} \left( 1 - \Fr^{\abs{o}} T^{\abs{o}} \left| \, \H^{1}(X)^{(j_o, \alpha_{o})} \right. \right)  = \left( 1 - \gamo \, T^{\abs{o}}  \right),
\end{equation*}
Therefore, we conclude that
\begin{equation}\label{eq: Charpol H1X}
	\chpolg{\H^{1}(X)} = \prod_{o \, \in \Oqa} \left( 1 - \gamo \,  T^{\abs{o}}\right).
\end{equation}

\subsection{Geometry of the minimal regular model}\label{sect: Geom Cmod} 

In this section we give a more thorough analysis of the geometry of the minimal regular model $\Cmod$. We now show that $\Cmod$ is dominated by the product of curves $X \times Y$, thus proving the claim in Section \ref{sect: Constructing ab surf}. Recall from Section \ref{sect: Constructing ab surf} that there is a rational map 
\begin{equation*}
	\ratmap{\vartheta_{0}}{X \times Y}{\Zsurf}{(u, t_1, v, t_2)}{(u, v, t_2 - t_1),}
\end{equation*}
to the surface $\Zsurf / \Fq$ defined by \eqref{eq: Model Zsurfa}. 
\begin{proposition}\label{prop: DPC}
	The map $\vartheta_{0}$ induces a birational map $\vartheta_1 : (X \times Y) / \Fqa \stackrel{\sim}\dashrightarrow \Zsurf$. In particular, the minimal regular model $\Cmod$ of $C$, as a surface over $\Fq$, is dominated by $X \times Y$.
\end{proposition}

\begin{proof}
	We saw above that the group $\Fqa$ acts on both of the curves $X$ and $Y$, so that the product $X \times Y$ also carries a ``diagonal" action of $\Fqa$. This action is given in affine coordinates by $\alpha \cdot \left(u, t_1, v, t_2\right)=\left(u, t_1+\alpha, v, t_2+\alpha\right)$ for any $\alpha \in \Fqa$, any $(u, t_1, v, t_2) \in X \times Y$. Its only fixed points are $(\infty_{X}, \infty_{Y_i})$, $i = 1, 2$. If $(u, t_1, v, t_2) \in X \times Y$ and $\alpha \in \Fqa$, then $\vartheta_0(\alpha \cdot (u, t_1, v, t_2)) = \vartheta_0(u, t_1, v, t_2)$. Thus, $\vartheta_{0}$ factors through the quotient $(X \times Y) / \Fqa$ and gives a rational map $\vartheta_1 : (X \times Y) / \Fqa \dashrightarrow \Zsurf$. Counting elements in the fibers, one checks that both $\vartheta_0$ and the quotient map $X \times Y \rightarrow (X \times Y) / \Fqa$ are of degree $q^a$. We deduce that $\vartheta_1$ has degree $1$, hence is a birational map. To prove the last statement, recall that both $\Zsurf$ and $\Cmod$ are models of $C$ over $\Pp$, so there exists a birational map $\Zsurf \stackrel{\sim}\dashrightarrow \Cmod$.  The composition
	\begin{equation*}
		\begin{tikzcd}
			X \times Y \arrow[r, two heads] & (X \times Y) / \Fqa  \arrow[r, "\sim", "\vartheta_1"', dashed] & \Zsurf \arrow[r, "\sim", dashed] & \Cmod
		\end{tikzcd}
	\end{equation*}
	is dominant, finally giving the desired result.
\end{proof}
\vspace{0.75em}
In order to compute the $L$-function of $S$, we need to understand the indeterminacy of $\vartheta_1$ better. Following \cite[Proposition 3.1.5]{PriesUlmer}, one constructs by successive blow-ups a surface $\Xmod / \Fq$ and a proper birational morphism $\res_1 : \Xmod \rightarrow X \times Y$ resolving the indeterminacy of both $\delta_0 : X \times Y \dashrightarrow \Pp$, $(u, t_1, v, t_2) \mapsto t_2 - t_1$ and $\vartheta_0: X \times Y \dashrightarrow \Zsurf$. This construction gives morphisms $\delta : \Xmod \rightarrow \Pp$ and $\vartheta : \Xmod \rightarrow \Zsurf$ such that the diagram below commutes:

\begin{equation*}
\begin{tikzcd}
	\Xmod \arrow[d, "\res_1"] \arrow[rrd, "\vartheta", bend left] \arrow[dd, "\delta"', bend right = 40, shift right=2] &   &  \\
	X \times Y \arrow[d, "\delta_0", dashed] \arrow[rr, "\vartheta_0", bend left, dashed] \arrow[r]  & (X \times Y) / \Fqa \arrow[ld, dashed, bend left = 10] \arrow[r, "\vartheta_1", dashed] & \Zsurf \arrow[lld, bend left=20] \\
	\Pp  &   &  
\end{tikzcd}	
\end{equation*}
The diagonal action of $\Fqa$ on $X \times Y$ lifts uniquely to $\Xmod$ and fixes the exceptional divisor pointwise. Thus the morphism $\vartheta : \Xmod \rightarrow \Zsurf$ factors through the quotient $\Xmod / \Fqa$. The latter surface might be singular. Again, resolving singularities (see \cite[§8.3.4]{Liubook}), we let $\res_2 : \Xhat \rightarrow \Xmod / \Fqa$ be a proper birational morphism such that $\Xhat$ is a smooth projective surface and the composition $\Xmod / \Fqa \rightarrow \Zsurf \dashrightarrow \Cmod$ induces a morphism $\Xhat \rightarrow \Cmod$. We summarize the situation in the commutative diagram below:
\begin{equation}\label{eq: Final diagram}
\begin{tikzcd}
	& \Xhat & \Cmod \\
	\Xmod & { \Xmod / \Fqa} & \Zsurf \\
	X \times Y & {(X \times Y) / \Fqa} \\
	& \Pp
	\arrow[two heads, from=3-1, to=3-2]
	\arrow[dashed, from=3-2, to=4-2]
	\arrow["\delta_0", bend right = 25, dashed, from=3-1, to=4-2]
	\arrow["\sim" {anchor=south, rotate=90}, "\res_1",  from=2-1, to=3-1]
	\arrow[two heads, from=2-1, to=2-2]
	\arrow[from=2-2, to=3-2]
	\arrow[two heads, from=2-2, to=2-3]
	\arrow["\sim" {anchor=south, rotate=90}, "\res_2", from=1-2, to=2-2]
	\arrow[two heads, from=1-2, to=1-3]
	\arrow["\sim" {anchor=south, rotate=90}, dashed, from=2-3, to=1-3]
	\arrow["\vartheta_1" {rotate =35}, "\sim"{anchor=north, rotate=35} , dashed, from=3-2, to=2-3]
	\arrow[bend left = 35, from=2-3, to=4-2]
\end{tikzcd}
\end{equation}

\noindent We now describe how exceptional divisors behave under the maps in diagram \eqref{eq: Final diagram}. Recall that $C$ has two points at infinity, denoted by $\infty_1, \infty_2$, which are $\K$-rational.

\begin{lemma}\label{lem: Exc div mapping} Let $s_{\infty_i} : \Pp \rightarrow \Cmod$, $i \in \lbrace 1, 2 \rbrace$, denote the sections of $\Cmod \rightarrow \Pp$ corresponding to the points at infinity on $C / \K$. 
	\begin{enumerate}
		\item The strict transforms of $\infty_{X} \times Y$ and $X \times \infty_{Y}$ in $\Xhat$ map to $\Cmod_{\infty}$.
		\item Let $\Emod$ be the image in $\Xmod / \Fqa$ of the exceptional divisor of $\Xmod \rightarrow X \times Y$. Then the strict transforms of the components of $\Emod$ in $\Xhat$ map either to $\Cmod_{\infty}$ or to the image of the infinity sections $s_{\infty_{i}}$.
		\item Every component of the exceptional divisor of $\Xhat \rightarrow \Xmod / \Fqa$ maps to a fiber of $\Cmod \rightarrow \Pp$.
	\end{enumerate}
\end{lemma}

\begin{proof}
		The formula defining $\delta_0$ shows that for any $y \in Y(\Fqbar) \setminus \lbrace \infty_{Y} \rbrace$, $\delta_0$ maps $(\infty_{X}, y)$ to the point at infinity in $\Pp$. In the same way, for $x \in X(\Fqbar) \setminus \lbrace \infty_{X} \rbrace$, we have $\delta_0(x, \infty_{Y}) = \infty$. Commutativity of diagram \eqref{eq: Final diagram} gives the desired result. The second statement is simply Remark $3.1.6$ in \cite{PriesUlmer}. Finally, to prove the third point, notice that any component $\Gamma$ of the exceptional divisor of $\Xhat \rightarrow \Xmod / \Fqa$ lies over a unique point $P \in (\Xmod / \Fqa)(\Fqbar)$. If $t \in \Pp(\Fqbar)$ is the image of $P$ by the composition $\Xmod / \Fqa \rightarrow \Zsurf \rightarrow \Pp$, then $\Xhat \rightarrow \Cmod$ maps $\Gamma$ to the fiber $\Cmod_{t}$.
\end{proof}


\section[The L-function of S]{The $L$-function of $S$}\label{sect: L func}
In this section we use the geometric analysis of $\Cmod$ from the previous section to show that $S_a$ satisfies the BSD conjecture. We then provide an explicit formula for the $L$-function of our abelian surface $S$ and conclude by computing the Mordell--Weil rank and the special value $\spval$.

\subsection[Definition of L(S, T) and statement of the BSD conjecture]{Definition of $L(S, T)$ and statement of the BSD conjecture}\label{sect: Def Lfunc and BSD}
Recall that in Section~\ref{sect: Auxiliary curves} we introduced the $\ell$-adic étale cohomology functors $\Het$ with coefficients in $\Qlbar$. Let us write also $\HS:= \H_{\acute{e}t}^{1}(S \times_{\K} \Ksep, \Qlbar)$ and $\H^{1}(C) := \H_{\acute{e}t}^{1}(C \times_{\K} \Ksep, \Qlbar)$. It is well-known that we have an isomorphism $\HS \simeq \H^{1}(C)$ (see for example \cite[Corollary 9.6]{Milne_Jac}). This is a $\Qlbar$-vector space of dimension $4$ on which the absolute Galois group $\Gal(\Ksep / \K)$ acts. For any place $v$ of $\K$ corresponding to a geometric point of $\Pp$ (also denoted by $v$), recall that $I_v$ denotes the inertia subgroup at $v$ and let $\H^{1}(S)^{I_v}$ be the subspace of $\H^{1}(S)$ fixed by inertia. Moreover, we denote by $\Frv$ the geometric Frobenius at $v$, induced by $x \mapsto x^{q^{\dv}}$.

\begin{definition}
	The Hasse-Weil $L$-function of $S$ is defined by the Euler-product
	\begin{equation*}
		L(S, T)=\prod_{v \in \Mpl} \chpolIn{1}{S}^{-1}.
	\end{equation*}
\end{definition}
The $L$-function can be expanded into a formal power series in $T$. The Hasse-Weil bound on the eigenvalues of $\Frv$ acting on $\HS$ implies that this power series converges on the complex disc $\left\{T \in \C:\abs{T}<q^{-3 / 2}\right\}$. The theory of étale cohomology developed by Grothendieck \textit{et. al.} (see \cite{Milne80}) and the work of Deligne on the Weil conjectures (\cite{Deligne80}) led to the following theorem, describing the main properties of the $L$-function of $S$.

\begin{theorem}\label{thm: Weil Conj}
	Let $S/\Fqt$ be as above and $N(S) \in \Div(\Pp)$ be its conductor divisor.
	\begin{enumerate}
		\item (Rationality) The L-function $L(S, T)$ is a polynomial in $T$ with integer coefficients. Writing $b(S)$ for its degree, we have the equality $b(S) = \deg N(S) -8$.
		\item (Functional equation) There exists $w(S) \in\{ \pm 1\}$ such that $L(S, T)$ satisfies
		\begin{equation*}
			L(S, T)=w(S)\,(qT)^{b(S)}\, L\left(S, \left(q^2 T\right)^{-1}\right).
		\end{equation*}
		\item (Riemann Hypothesis) If $z \in \C$ is such that $L(S, z)=0$, then $\abs{z}=q^{-1}$.
	\end{enumerate}
\end{theorem}
\noindent The reader can find a sketch of the proof in \cite[§6.2.1]{Ulmer14}. 

\begin{remark} 
	\begin{enumerate}		
		\item In general, the $L$-function of an abelian variety over $\K$ is a rational function in $T$ with integral coefficients. Since $S$ is non-isotrivial, \cite[Lemma 6.2.2]{Ulmer14} implies that $L(S, T)$ is actually a polynomial in $T$.
		
		\item We can define a function of the complex variable $s \mapsto L(S, q^{-s})$, which converges in the right half-plane $\Re(s) > 3/2$ because of the Hasse-Weil bound. The functional equation $ii)$ shows that it extends to a holomorphic function on the whole complex plane. The third point in Theorem~\ref{thm: Weil Conj} is called the Riemann hypothesis because it states that all the zeros of $s \mapsto L(S, q^{-s})$ lie on the vertical line $\Re(s) = 1$. 
	\end{enumerate}	
\end{remark}

\noindent The special value $\spval$ of the $L$-function of $S$ at $q^{-1}$ is defined as
\begin{equation*}\label{eq: Def special value}
	\spval := \left.\frac{L(S, T)}{(1-q T)^{\rkan}}\right|_{T=q^{-1}} \qquad \text{ where } \rkan = \ord_{T = q^{-1}}L(S, T).
\end{equation*}
The integer $\rkan$ is called the analytic rank of $S$. The next statement is called the Birch and Swinnerton-Dyer conjecture (BSD): it is conjectured to be true for any abelian variety over $\K$, and remains an open question in general. We refer the reader to \cite{Tate66} for a survey. In \cite{PriesUlmer}, the authors give a recipe to construct sequences of abelian varieties over function fields satisfying the BSD conjecture. As already mentioned in Section~\ref{sect: Geom Cmod}, the construction of our curve $C$ and its Jacobian $S$ follows that of Pries and Ulmer, to ensure that $S$ has the desired properties.

\begin{theorem}\label{thm: BSD}
	The abelian surface $S$ satisfies the BSD conjecture, \ie,  
	\begin{itemize}
		\item the algebraic and analytic ranks of $S$ coincide, that is, $\rank S(\K) = \rkan.$
		
		\item The group $\Sh(S / \K)$ is finite. 
		
		\item The BSD formula holds:
		\begin{equation}\label{eq: Strong BSD}
			\spval = \frac{\abs{\ShaS} \Reg(S) \, q^2 \prod_{v} c_v(S)}{\abs{\tor{S(\K)}}^{2} \, H(S)}.
		\end{equation}
	\end{itemize} 
\end{theorem}

\begin{proof}
	Kato and Trihan have proved that the three statements are equivalent (\cite{KatoTrihan}). Thus it suffices to establish the first point, called the ``weak BSD conjecture". We relate this to Tate's conjecture for $\Cmod$. Tate's conjecture for a variety $V$ defined over $\Fq$ relates the rank of the Néron--Severi group $\NS(V)$ to the order of the pole of the zeta function of $V$ at $s=1$ (see \cite[§6]{Ulmer14} for precise statements). Now Tate's conjecture is known to be true for curves, and Tate's theorem on homomorphisms of abelian varieties implies that products of curves satisfy Tate's conjecture. Moreover, the latter is preserved by domination (see \cite[§2.11 \& §$2.12$]{Ulmer11}). 
	
	In Proposition~\ref{prop: DPC} we saw that the minimal regular model $\Cmod / \Fq$ is dominated by the product of curves $X \times Y$, and the discussion above implies that $\Cmod$ satisfies Tate's conjecture. But Shioda--Tate theorem implies that Tate's conjecture for $\Cmod$ is equivalent to the BSD conjecture for $S$, hence the result (the reader can find detailed proofs in \cite[§$4.2$, §$6.2.2$]{Ulmer14}). 
\end{proof}

\subsection[Explicit expression for L(S, T)]{Explicit expression for $L(S, T)$} 

Let $a \geq 1$, recall from Section~\ref{sect: Orbits} that $\Oqa$ denotes the set of orbits of $\left(\Z / 3\Z \right)\st \times \Fqa\st$ under the action of $\langle q \rangle$ . Moreover, to any $o \in \Oqa$, we have attached character sums $\gamo$ and $\Klo$.

\begin{theorem}\label{thm: L func formula}
	The $L$-function of $S$ is given by the formula
	\begin{equation}\label{eq: L func expl}
		L(S, T) = \prod_{o \, \in \Oqa} \left( 1 - \gamo \Klo \,  T^{\abs{o}} + \gamo^2 \, q^{\abs{o}} \, T^{2 \abs{o}} \right).
	\end{equation} 
\end{theorem}	

\noindent We devote the rest of section~\ref{sect: L func} to the proof of theorem~\ref{thm: L func formula}.  

\subsection[Cohomological computation of L(S,T)]{Cohomological computation of $L(S,T)$}\label{sect: Comput L section} 

We begin by relating $L(S, T)$ to the zeta function of $\Cmod / \Fq$. Since we will be manipulating zeta functions, we start by recalling their definition and some of their basic facts (see \cite{Serre70} for details).
	
\begin{definition}
	Let $\Ff / \Fq$ be a finite extension. For a variety $V$ over $\Ff$, define its zeta function as 
	\begin{equation*}
		Z(V, T)=\prod_{P \in|V|}\left(1-T^{\operatorname{deg} P}\right)^{-1},
	\end{equation*}
	where the product runs over the set of closed points of $V$. 
\end{definition}
As an example, one can compute $Z(\Pp, T) = \left( (1-T)(1-qT)\right)^{-1}$. When $V$ is smooth and projective, the Grothendieck-Lefschetz trace formula implies the rationality of $Z(V, T)$ (see \cite[Chapter VI, Theorem 12.4]{Milne80})
\begin{equation}\label{eq: Groth Lef}
	Z(V, T)=\prod_{i=0}^{2 \operatorname{dim} V}(-1)^{i+1} 
	\operatorname{det}\left(1 - \Fr^{[\Ff : \Fq]} T \, \left | \, \H^{i}(V) \right. \right).
\end{equation}
Besides, the Frobenius $\Fr$ acts as the identity in $V$, so we know its action on $\H^{0}(V)$. Using Poincaré duality (\cite[Theorem 11.1]{Milne80}) we can also deduce its action on $H^{2 \operatorname{dim} V}(V)$:
\begin{equation*}
	\operatorname{det}\left(1 - \Fr^{[\Ff : \Fq]} T \, \left | \, \H^{0}(V) \right. \right) = 1-T \quad \text{and} \quad 
	\operatorname{det}\left(1 - \Fr^{[\Ff : \Fq]} T \, \left | \, \H^{2 \operatorname{dim} V}(V) \right. \right) = 1-q^{[\Ff : \Fq] \operatorname{dim} V}T.
\end{equation*}

\begin{proposition}
	We have the equality
	\begin{equation}\label{eq: Lfunc form temp}
		L(S, T) = \frac{\chpol{2}{\Cmod}}{\chpol{1}{\Cmod} \, \chpol{3}{\Cmod} \, \left(1-qT\right)^{6q^a + 8}}.  		
	\end{equation}
\end{proposition}

\begin{proof}
	Formula \eqref{eq: Groth Lef} and the discussion below applied to the smooth surface $\Cmod / \Fq$ give
	\begin{equation}\label{eq: Z Cmod cohom}
		Z(\Cmod, T) = \frac{\chpol{1}{\Cmod} \, \chpol{3}{\Cmod}}{(1-T)\, \chpol{2}{\Cmod} \, (1-q^2T)}.
	\end{equation}
	By construction $\Cmod$, as a model of $C / \K$, is equipped with a surjective morphism $\Cmod \rightarrow \Pp$. We therefore decompose $\Cmod$ as a disjoint union of fibers of closed points of $\Pp$ (see \cite[§VI, 13.8]{Milne80}). Identifying closed points of $\Pp$ and places of $\K$, we rewrite $Z(\Cmod, T) = \prod_{v \in \Mpl} Z(\Cmod_{v}, T)$. Splitting the product into places of good and bad reduction for $C$, and using \eqref{eq: Groth Lef} for the smooth fibers $\Cmod_{v} / \Fv$ when $v$ is a place of good reduction we get: 
	\begin{equation}\label{eq: Z Cmod split}
		Z(\Cmod, T) = \prod_{v \text{ good}} \frac{\chpolv{1}{\Cmod_{v}}}{\left( 1 - T^{\dv}\right)\left( 1 - (qT)^{\dv}\right)} \, \prod_{v \text{ bad}} Z(\Cmod_{v}, T).
	\end{equation}
	Let $v \in \Mpl$ be a place of good reduction. As stated above, we have $\HS \simeq \H^{1}(C)$ and the Néron--Ogg--Shafarevich criterion (\cf \, \cite{SerreTate}) gives $\H^{1}(S)^{I_v} \simeq \H^{1}(S)$. Moreover, the smooth base change theorem (\cite[Chapter VI, Corollary 4.2]{Milne80}) implies the isomorphism $\H^{1}(C) \simeq \H^{1}(\Cmod_v)$. Therefore we have $\H^{1}(\Cmod_v) \simeq \HS^{I_v}$ and the local factor of $L(S, T)$ at $v$ is $\chpolv{1}{\Cmod_{v}}^{-1}$.

	Let now $v$ be a finite place of bad reduction for $C$, we relate $Z(\Cmod_v, T)$ to the zeta function of the normalization of $\Cmod_{v}$. From Proposition \ref{prop: Fibers} we know that this is the elliptic curve $E_v / \Fv$ described in affine coordinates by \eqref{eq: Def ell curve Ev}. \cite{AubryPerret} describes $Z(\Cmod_v, T)$ in terms of $Z(E_v, T)$, the number and degree of preimages of singular points of $\Cmod_v$ by the normalization morphism $\nu_v : E_v \rightarrow \Cmod_{v}$. The only singular point in $\Cmod_{v}$ is $(0,0)$, which is $\Fv$-rational. Moreover, $\nu_v$ is one-to-one as the preimage of the singular point of $\Cmod_{v}$ is $\nu_{v}^{-1}(0, 0) = (0,0)$, which is again $\Fv$-rational. Hence, \cite[Theorem 2.1]{AubryPerret} gives
	\begin{equation*}
		Z(\Cmod_{v}, T) = Z(E_v, T) = \frac{\chpolv{1}{E_{v}}}{\left( 1 - T^{\dv}\right)\left( 1 - (qT)^{\dv}\right)}.
	\end{equation*}
	Using the universal property of Néron models, one can show that $\HS^{I_v} \simeq \H^{1}(\Sner_{v})$. From Theorem \ref{thm: Red Smod} we know that the reduction $\Sner_v$ is the extension of $E_v$ by a unipotent group of dimension $1$, which has trivial $\ell$-adic cohomology, as $\ell \neq p$. Therefore we have $\HS^{I_v} \simeq \H^{1}(E_v)$ and the local factor of $L(S, T)$ at $v$ is $\chpolv{1}{E_{v}}^{-1}$.
	
	Plugging all this into \eqref{eq: Z Cmod split} we have
	\begin{equation*}
		Z(\Cmod, T) = \left( \prod_{v \text{ good}} \frac{\chpolIn{1}{S}}{\left( 1 - T^{\dv}\right)\left( 1 - (qT)^{\dv}\right)}\right) \left( \prod_{v \, \text{finite\,bad}} \frac{\chpolIn{1}{S}}{\left( 1 - T^{\dv}\right)\left( 1 - (qT)^{\dv}\right)} \right)  Z(\Cmod_{\infty}, T).
	\end{equation*}
	At infinity, $\Sner$ has totally unipotent reduction (\cf \ref{thm: Red Smod}), so $\H^{1}(S)^{I_v}$ is trivial and the local factor of $L(S, T)$ at $\infty$ is $1$. Regrouping the numerators of the terms inside the parenthesis above, we recognize the inverse of $L(S, T)$. On the other hand, regrouping the denominators inside the parenthesis, one recognizes two zeta functions of $\Pp$:
	\begin{equation*}
		\prod_{v \neq \infty} \left( (1 - T^{\dv}) (1 - (qT)^{\dv})\right)^{-1} = Z(\Pp, T) \, Z(\Pp, qT) \, (1- T) (1 - qT).
	\end{equation*}
	The last two factors correspond to the place at infinity that was missing from the left-hand side. Upon using the explicit expression for $Z(\Pp, T)$, we can rewrite
	\begin{equation}\label{eq: Z Cmod simplified}
		Z(\Cmod, T) = L(S, T)^{-1} \, \frac{Z(\Cmod_{\infty}, T)}{(1-qT)(1-q^2T)}.
	\end{equation}
	We now only have to compute $Z(\Cmod_{\infty}, T)$. Recall from Proposition \ref{prop: Fibers} (see Figure~\ref{fig: Fiber Cmod infinity}) that the fiber $\Cmod_{\infty}$ is a tree of $\Ppo$'s, which has $6q^a + 7$ irreducible components, all $\Fq$-rational. All singularities on $\Cmod_{\infty}$ are nodes, so by a direct point counting argument, we have
	\begin{equation}\label{eq: Z fiber infinity}
		Z(\Cmod_{\infty}, T) = \frac{Z(\Pp, T)^{6q^a + 7}}{Z(\Spec(\Fq), T)^{6q^a + 6}} = \frac{1}{(1-T)(1-qT)^{6q^a + 7}}.
	\end{equation}
	Plugging \eqref{eq: Z fiber infinity} into \eqref{eq: Z Cmod simplified}, combining it with \eqref{eq: Z Cmod cohom}, and rearranging terms gives the desired result.
\end{proof}

Some of the terms in formula \eqref{eq: Lfunc form temp} are easily described:

\begin{proposition}\label{prop: Vanish cohom}
	The first and third cohomology groups of $\Cmod$ are trivial. Therefore we have
	\begin{equation*}
		\chpol{1}{\Cmod} = \chpol{3}{\Cmod} = 1.
	\end{equation*}
\end{proposition}

\begin{proof}
	Recall from Section~\ref{sect: Geom Cmod} that the morphism $\Xhat \rightarrow \Cmod$ is dominant, so proving that $\H^{1}(\Xhat) = 0$ would imply that $\H^{1}(\Cmod)$ is trivial too. Recall from Section \ref{sect: Geom Cmod} that $\Xhat$ is obtained from the quotient $\Xmod / \Fqa$ by a succession of blow-ups at closed points of $\Xmod / \Fqa$ and at the image of the exceptional divisor $\Emod$ (which has trivial $\H^1$ because it is the image of a union of $\Ppo$'s). Thus $\H^{1}(\Xhat)$ decomposes as the direct sum of $\H^{1}(\Xmod / \Fqa)$ and a certain number of copies of $\H^{1}(\Ppo)$, all of which are trivial. We deduce the isomorphisms
	\begin{equation*}
		\H^{1}(\Xhat) \simeq \H^{1}(\Xmod / \Fqa) \simeq \H^{1}(\Xmod)^{\Fqa},
	\end{equation*}
	where the upper-script denotes the $\Fqa$-invariant subspace. As explained in Section \ref{sect: Geom Cmod}, $\Xmod$ is obtained from $X \times Y$ by a series of blow-ups, so that, as above, $\H^{1}(\Xmod) \simeq \H^{1}(X \times Y)$. Applying the Künneth formula (\cite[hapter VI, Theorem 8.5]{Milne80}), we have
	\begin{align*}
	\H^{1}(\Xhat) \simeq \H^{1}(X \times Y)^{\Fqa} & \simeq \left( \H^{0}(X) \otimes \H^{1}(Y) \oplus \H^{1}(X) \otimes \H^{0}(Y)\right)^{\Fqa} \\
	& \simeq \H^{0}(X)^{\Fqa} \otimes \H^{1}(Y)^{\Fqa} \oplus \H^{1}(X)^{\Fqa} \otimes \H^{0}(Y)^{\Fqa}.
	\end{align*}
	Recall from Lemmas \ref{lem: Eigenvalues Kloos}, \ref{lem: Eigenvalues gamma} that  $\H^{1}(X)^{\Fqa}$ and $\H^{1}(Y)^{\Fqa}$ are both trivial, so $\H^{1}(\Xhat) = 0$. Therefore $\H^{1}(\Cmod) = 0$, and, by Poincaré duality (\cite[Corollary 11.2]{Milne80}), that $\H^{3}(\Cmod) = 0$ as well.
\end{proof}

The vanishing of the first and third cohomology groups of $\Cmod$ clarifies the formula \eqref{eq: Lfunc form temp} of the $L$-function of $S$. We now further simplify the expression of $L(S, T)$. Recall that the Néron-Severi group $\NS(\Cmod)$ is defined as the group of divisors of $\Cmod$ modulo algebraic equivalence (see \cite[Chapter III]{Milne80}). We will denote by $[\Gamma]$ the class in $\NS(\Cmod)$ of a divisor $\Gamma \in \Div(\Cmod)$. Following \cite{Shioda}, we define the trivial lattice $T \subset \NS(\Cmod)$ as the sublattice generated by the class of the infinity section $\left[s_{\infty_1} \right]$ and the classes of the irreducible components of the fibers of $\Cmod \rightarrow \Pp$ (we identify the section $s_{\infty_1}$ and its image in $\Div(\Cmod)$). We write $\gamma_{\Cmod} : \NS(\Cmod) \otimes \Zl \rightarrow \H^{2}(\Cmod)$ for the map induced by the Kummer exact sequence, and we define the subspace $\Lambda := \gamma_{\Cmod}(T) \subset \H^{2}(\Cmod)$.
 
\begin{lemma}\label{lem: L and lambda}
	We have 
	\begin{equation*}
		L(S, T) = \chpolg{\H^{2}(\Cmod) / \Lambda}.
	\end{equation*}
\end{lemma}

\begin{proof}
	Recall from Lemma \eqref{prop: Fibers} that all fibers of $\Cmod \rightarrow \Pp$ at finite places are irreducible, and they are all algebraically equivalent. Fix a finite place $v$. $T$ is generated by $[s_{\infty_1}]$, $[\Cmod_v]$, and $\left([\Gamma_i]\right)_{i}$ where the $\Gamma_{i}$'s are the irreducible components of $\Cmod_{\infty}$ not meeting $s_{\infty_1}$. As Figure \ref{fig: Fiber Cmod infinity} in Proposition \ref{prop: Fibers} shows, there are $6q^a + 6$ of these irreducible components, that are all $\Fq$-rational. Using the intersection pairing, one sees that $T$ is not spanned by any subfamily of these $6q^a + 8$ classes of divisors, and so $\Lambda \subset \H^{2}(\Cmod)$ has dimension $6q^a + 8$. The Frobenius acts by multiplication by $q$ on $\gamma_{\Cmod}([s_{\infty_1}])$, $\gamma_{\Cmod}([\Cmod_v])$ and each of the $\gamma_{\Cmod}([\Gamma_i])$, so we obtain
	\begin{equation*}
		\chpolg{\Lambda} = (1 - qT)^{6q^a + 8}.
	\end{equation*} 
	But we have $\chpolg{\H^{2}(\Cmod) / \Lambda} = \chpolg{\H^2(\Cmod)} / \chpolg{\Lambda}$, so combining the discussion above with \eqref{eq: Lfunc form temp} and Proposition~\ref{prop: Vanish cohom} gives the desired result.
\end{proof}

We have now all the necessary tools to prove Theorem~\ref{thm: L func formula}.

\begin{proof}[Proof of Theorem~\ref{thm: L func formula}]
Given Lemma \ref{lem: L and lambda}, we need to compute the characteristic polynomial of $\Fr$ on $\H^2(\Cmod) / \Lambda$. We start by relating $\H^{2}(\Cmod) / \Lambda$ to the cohomology groups of the auxiliary curves $X$ and $Y$. Since the morphism $\Xhat \rightarrow \Cmod$ is dominant, the map $\H^2(\Xhat) \rightarrow \H^2(\Cmod)$ is surjective. Consider the subspace $B \subset \H^2(\Xhat)$ spanned by the classes of the components of the exceptional divisor of $\Xhat \rightarrow \Xmod / \Fqa$ and the strict transform of $\Emod$ in $\Xmod$. Then 
\begin{align*}
	\H^2(\Xhat) & \simeq \H^2(\Xmod / \Fqa) \oplus B \\
	& \simeq \H^{2}(\Xmod)^{\Fqa} \oplus B \\
	& \simeq \H^2(X \times Y)^{\Fqa} \oplus B \\
	& \simeq \left( \H^{0}(X) \otimes \H^{2}(Y) \right)^{\Fqa} \oplus \left( \H^{1}(X) \otimes \H^{1}(Y) \right)^{\Fqa} \oplus \left( \H^{2}(X) \otimes \H^{0}(Y) \right)^{\Fqa} \oplus B,
\end{align*}
where the last isomorphism follows from the Künneth formula (see \cite{Milne80}). Note that all the isomorphisms above are compatible with the action of $\Gal\left(\Fqbar / \Fq\right)$. Now $\H^{0}(X) \otimes \H^{2}(Y)$ and $H^{2}(X) \otimes \H^{0}(Y)$ are respectively generated by the strict transforms of $X \times \infty_{Y}$ and $\infty_{X} \times Y$ in $\Xhat$. Lemma~\ref{lem: Exc div mapping} implies that $\H^{0}(X) \otimes \H^{2}(Y)$, $\H^{2}(X) \otimes \H^{0}(Y)$, and $B$ all map to the trivial lattice under $\H^{2}(\Xhat) \rightarrow \H^{2}(\Cmod)$. We deduce from this the existence of a surjective map 
\begin{equation*}
	\Upsilon : \left( \H^{1}(X) \otimes \H^{1}(Y) \right)^{\Fqa} 
	\twoheadrightarrow \H^{2}(\Cmod) / \Lambda.
\end{equation*}
Moreover, $\Upsilon$ is $\Gal(\Fqbar / \Fq)$-equivariant because all the isomorphisms above are compatible with the action of $\Gal(\Fqbar / \Fq)$. We further describe the left-hand side. Recall from Section~\ref{sect: Auxiliary curves} that the cohomology groups $\H^{1}(X)$ and $\H^{1}(Y)$ decompose into isotypical components: \eqref{eq: Decomp isotyp Y} and \eqref{eq: Decomp isotyp X} yield
\begin{equation*}
	\H^{1}(X) \otimes \H^{1}(Y) = \bigoplus_{(j, \alpha) \in \Ia} \bigoplus_{\beta \in \Fqa\st} \H^{1}(X)^{(j, \alpha)} \otimes \H^{1}(Y)^{\beta}.
\end{equation*}
The group $\Fqa$ acts on $\H^{1}(X)$ (resp. $\H^{1}(Y)$) by multiplication by the character $\chi$ associated to $(j, \alpha)$ (resp. $\psi_{\Fqa, \beta}$). Thus, among all the planes $\H^{1}(X)^{(j, \alpha)} \otimes \H^{1}(Y)^{\beta}$, the $\Fqa$-invariant ones are those for which $\alpha = \beta$. We get
\begin{equation*}
	\left( \H^{1}(X) \otimes \H^{1}(Y) \right)^{\Fqa} = \bigoplus_{\substack{(j, \alpha) \in \Ia \\ \beta \in \Fqa\st}} \left( \H^{1}(X)^{(j, \alpha)} \otimes \H^{1}(Y)^{\beta}\right)^{\Fqa}  = \bigoplus_{(j, \alpha) \in \Ia} \H^{1}(X)^{(j, \alpha)} \otimes \H^{1}(Y)^{\alpha}.
\end{equation*}
This space has therefore dimension $2 \abs{\Ia} = 4(q^a -1)$. Because of Lemma~\ref{lem: L and lambda}, the dimension of $\H^{2}(\Cmod) / \Lambda$ equals the degree of $L(S, T)$. Now the Grothendieck-Ogg-Shafarevich formula (\cite[Chapter V, Theorem 2.12]{Milne80}) implies that $\deg L(S, T) = \deg N_S - 8 = 4(q^a -1)$. By equality of dimensions, $\Upsilon$ is actually an isomorphism. We deduce that
\begin{equation*}
	L(S, T) = \chpolg{\left( \H^{1}(X) \otimes \H^{1}(Y) \right)^{\Fqa}}.
\end{equation*}
For any $o \in \Oqa$ and any choice of representative $(j_o, \alpha_o) \in o$, the $\abs{o}$-th iterate of $\Fr$ stabilizes the plane $\H^{1}(X)^{(j_o, \alpha_o)} \otimes \H^{1}(Y)^{\alpha_o}$ (see Lemmas \ref{lem: Eigenvalues Kloos} and \ref{lem: Eigenvalues gamma}). The computation of the eigenvalues of $\Fr^{\abs{o}}$ done in Section \ref{sect: Auxiliary curves} gives
\begin{equation*}
	\operatorname{det}\left(1 - \Fr^{\abs{o}} T \, \left | \, \H^{1}(X)^{(j_o, \alpha_o)} \otimes \H^{1}(Y)^{\alpha_o} \right. \right) = (1 - \gamo \kapo \, T)(1 - \gamo \kapobis \, T).
\end{equation*}
Applying \cite[Lemma 2.2]{Ulmer07} again, we deduce
\begin{equation*}
	\chpolg{\bigoplus_{(j, \alpha) \in \, o} \H^{1}(X)^{(j, \alpha)} \otimes \H^{1}(Y)^{\alpha}} = \left(1 - \gamo \kapo T^{\abs{o}} \right) \left(1 - \gamo \kapobis T^{\abs{o}} \right).
\end{equation*}
Finally, developing the right-hand-side and taking the products over all $o \in \Oqa$ we conclude
\begin{equation*}
	L(S, T) = \prod_{o \in \Oqa} \left( 1 - \gamo \Klo \,  T^{\abs{o}} + \gamo^2 \, q^{\abs{o}} \, T^{2 \abs{o}} \right).
\end{equation*}
\end{proof}

\subsection[Special value and rank of S]{Special value and rank of $S$}\label{sect: Special value section} 
From now on we fix a complex embedding $\iota : \Qbar \hookrightarrow \C$. Thanks to~\refGauss{Gauss norm} and \refKloos{Kloos norm} we know that, for any $o \in \Oqa$, $\abs{\iota(\kapo)} = \abs{\iota(\kapobis)}  = \abs{\iota(\gamo)} = q^{\abs{o} / 2}$ and this implies that $\iota(\kapo)$ and $\iota(\kapobis)$ are complex conjugate.

\begin{definition}\label{def: Angles theto gamo}
	Let $o \in \Oqa$, we define the Gauss and Kloosterman angles attached to $o$ as the unique elements $\epso \in \interv$ and $\theto \in \intrv$ such that
	\begin{equation}\label{eq: Def angles}
		\iota(\gamo) = q^{\abs{o}/2} \, \e^{i \, \epso}   \quad \text{and} \quad \iota(\Klo) = 2 \, q^{\abs{o}/2}  \cos(\theto).
	\end{equation}
\end{definition}
\noindent Up to relabelling $\kapo$ and $\kapobis$, we may assume that $\theto$ is the argument of $\kapo$ and $- \theto$ is the argument of $\kapobis$, \ie, \ $\iota( \kapo) = q^{\abs{o}/2} \, \e^{i \, \theto}$ and $\iota( \kapobis) = q^{\abs{o}/2} \, \e^{- i \, \theto}$.

\vspace{0.5em}
We describe now the analytic rank of $S_a$ and the special value of $L(S, T)$ in terms of these angles.

\begin{theorem}\label{thm: Special value thm}
	The order of vanishing of $L(S, T)$ at $T = q^{-1}$ is $0$. The special value of $L(S, T)$ equals
	\begin{equation}\label{eq: Special value}
		\spval = \prod_{o \, \in \Oqa} 4 \, \abs{\sin \left( \frac{\epso + \theto}{2}\right) \, \sin \left( \frac{\epso - \theto}{2}\right)}.
	\end{equation}
\end{theorem}

\begin{proof}
	Using the notation introduced in Definition~\ref{def: Angles theto gamo} we rewrite the result of Theorem~\ref{thm: L func formula}:
	\begin{equation}\label{eq: L func angles}
		L(S, T) = \prod_{o \, \in \Oqa} \left( 1 - \e^{i (\epso + \theto )} (qT)^{\abs{o}}\right) \left( 1 - \e^{i (\epso - \theto )} (qT)^{\abs{o}}\right).
	\end{equation}
	For any $o \in \Oqa$, \refKloos{Kloos conj valuation} implies that $\lbrace \valp(\kapo), \valp(\kapobis) \rbrace = \lbrace 0, \abs{o} \rbrace$ and Lemma~\ref{lem: Valuations gauss} gives $\valp (\gamo) \in \left \lbrace \frac{\abs{o}}{3}, \frac{\abs{o}}{2}, \frac{2 \abs{o}}{3} \right \rbrace$. Thus for $i = 1, 2$ the $\pfrac$-adic valuation of the algebraic integer $\gamo \kappa_{i}(o)$ is different from $\abs{o}$, hence $\gamo \kappa_{i}(o) q^{-\abs{o}} \neq 1$. This means that, when evaluating \eqref{eq: L func angles} at $T = q^{-1}$, no term in the product vanishes. Therefore the order of vanishing of $L(S, T)$ at $T = q^{-1}$ is $0$, and the special value is simply $L(S, q^{-1})$. Substituting $T = q^{-1}$ in \eqref{eq: L func angles} gives
	\begin{equation*}
		\spval = \prod_{o \, \in \Oqa} \left( 1 - \e^{i (\epso + \theto )}\right) \left( 1 - \e^{i (\epso - \theto )}\right).
	\end{equation*}
	Elementary trigonometric manipulations then show that 
	\begin{equation}\label{eq: Spval temp trig}
		\spval = \prod_{o \, \in \Oqa} 4 \, \e^{i \epso} \sin \left( \frac{\epso + \theto}{2}\right) \sin \left( \frac{\epso - \theto}{2}\right).
	\end{equation}
	The Riemann hypothesis for $L(S, T)$ (see Theorem \ref{thm: Weil Conj}) implies that $\spval$ is a positive real number, so $\abs{\spval} = \spval$. Taking absolute values on both sides of \eqref{eq: Spval temp trig} concludes the proof.
\end{proof}

Since $S$ satisfies the BSD conjecture (see Theorem \ref{thm: BSD}), we deduce the rank of its Mordell--Weil group, as well as the value of its Néron--Tate regulator.

\begin{corollary}\label{cor: Alg rank 0}
	We have $\rank(S(\K)) = 0$. In other words, $S(\K) = \tor{S(\K)}$ is finite. Moreover $\Reg(S) =1$.
\end{corollary}


\section{Distribution of angles}\label{sect: Angles}

Formula \eqref{eq: Special value} describing the special value of $L(S, T)$ involves angles of Gauss and Kloosterman sums. In order to study the size of $\spval$ in the following section, we describe here how these angles distribute in the interval $\interv$. 

\subsection{Angles of Gauss sums}  
As discussed in Section~\ref{sect: Orbits}, when $q \equiv 1 \mod 3$, there is a map $\prun : \Oqa \rightarrow \lbrace \pm 1 \rbrace$: for any $o \in \Oqa$, $\prun(o)$ is a lift to $\lbrace \pm 1 \rbrace$ of the ``first coordinate" $j$ of any representative $(j, \alpha) \in \Ia$ of $o$ .

\begin{proposition}\label{prop: Possible eps angles} 
	\begin{itemize}
		\item If $q \equiv 1 \mod 3$, there exists an irrational number $\phi \in \interv$, depending at most on $q$, such that for any $o \in \Oqa$, the angle $\epso$ satisfies
		\begin{equation*}
			\epso - \prun(o)\abs{o} \phi \in \left\lbrace 0, \frac{2 \pi}{3}, \frac{4 \pi}{3} \right \rbrace.
		\end{equation*}
		Moreover, for any element of the form $\eta + j d \phi$, where $\eta \in \left \lbrace 0, \frac{2 \pi}{3}, \frac{4 \pi}{3}  \right \rbrace$, $j \in \lbrace \pm 1 \rbrace$ and $d \divi a$ there is some $o \in \Oqa$ such that $\epso = \eta + j d \phi$.

		\item If $q \equiv 2 \mod 3$ then for any $o \in \Oqa$ we have $\epso \in \left \lbrace  0, \, \frac{\pi}{3}, \, \frac{ 2 \pi}{3}, \, \pi, \, \frac{4 \pi}{3}, \, \frac{5 \pi}{3}  \right \rbrace$.
	\end{itemize}
\end{proposition}

\begin{proof}
	We treat the case $q \equiv 1 \mod 3$ first. Let $g_{q, 1} := \Gauss{\Fq}{\chi_{\Fq, 3}^{ }}{\psi_{\Fq}}$, $g_{q, 2} := \Gauss{\Fq}{\chi_{\Fq, 3}^{2}}{\psi_{\Fq}}$, and $\phi$ to be the unique element of $\interv$ such that $\iota \left(g_{q, 1} \right) = \sqrt{q} \, \e^{i \phi}$. Since $\chi_{\Fq, 3}^{ }$ is a character of order $3$, $\chi_{\Fq, 3}^{2} = \overline{\chi_{\Fq, 3}^{}}$, and so $\iota(g_{q, 1})$ and $\iota(g_{q, 2})$ are complex conjugate, implying that $\iota(g_{q, 2}) = \sqrt{q} \, \e^{-i \phi}$. Let $o \in \Oqa$, and $(j, \alpha) \in \Ia$ be a representative of $o$. By definition we have $\prun(o) \equiv j \mod 3$, $\prun(o) \in \lbrace \pm 1 \rbrace$, and $\gamo = \Gauss{\Fo}{\chi_{\Fo, 3}^{j}}{\psi_{\Fo, \, \alpha}}$. Applying~\refGauss{Extract Gauss}, the Hasse--Davenport relation~\refGauss{Hasse Davenport} and the discussion above, this last equality becomes
	\begin{equation}\label{eq: Extract alpha}
		\gamo = (\chi_{\Fo, 3}^{ }(\alpha))^{-\prun(o)} \, g_{q, 1}^{\prun(o) \, \abs{o}}.	
	\end{equation}
	Now $(\chi_{\Fo, 3}^{ }(\alpha))^{-\prun(o)}$ is a third root of unity, so taking arguments of the associated complex numbers in \eqref{eq: Extract alpha} shows that $\epso - \prun(o)\abs{o} \phi \in \left\lbrace 0, \frac{2 \pi}{3}, \frac{4 \pi}{3} \right \rbrace$. A straightforward application of Lemma \ref{lem: Valuations gauss} when $a = 1$ shows that $\valp(g_{q, 1}) = 2/3$. Therefore $\valp \left(g_{q, 1} / \sqrt{q} \right) = 1/6 \neq 0$, and so $\e^{i \phi}$ is not a root of unity, implying that $\phi$ is irrational.
	
	Finally, fix $d \divi a$, $j \in \lbrace \pm 1 \rbrace$ and $\eta \in \left \lbrace 0, \frac{2 \pi}{3}, \frac{4 \pi}{3}  \right \rbrace$. If $\eta = 0$, then for any $\alpha \in \Ff_{q^d} \setminus \Ff_{q^{d-1}}$, equation \eqref{eq: Extract alpha} implies that the orbit $o \in \Oqa$ of $(j, \alpha^3)$ satisfies $\epso = j d \phi$. If $\eta \neq 0$, by an element-counting argument, one can check that there is always some $\alpha \in \Fqa\st$ such that $[\Fq(\alpha) : \Fq] = d$ and $\alpha$ is not a cube in $\Ff_{q^d}$. Switching $\alpha$ and $\alpha^2$ if necessary, we may assume that $\iota \left(\chi_{\Ff_{q^d}, 3}^{-j}(\alpha)\right) = \e^{i \eta}$ and again, using \eqref{eq: Extract alpha}, the orbit $o \in \Oqa$ of $(j, \alpha)$ satisfies the equality $\epso = \eta + j d \phi$.	
	
	When $q \equiv 2 \mod 3$, recall from the proof of Lemma~\ref{lem: Valuations gauss} that for any $o \in \Oqa$, any $(j, \alpha)\in o$ we have $\gamo = \pm \chi_{\Fo, 3}^{-j}(\alpha)(-q)^{\abs{o}/2}$. Therefore, $\gamo^6 = q^{3 \abs{o}}$ and $(\e^{i \epso})^6 = 1$. Hence, $\epso$ is as stated in the proposition.
\end{proof}

\begin{remark}\label{rmk: Possible gauss angles}
	When $q \equiv 2 \mod 3$, there are only $6$ possible values for the angles $\epso$ as $o$ ranges through $\Oqa$. When $q \equiv 1 \mod 3$, the irrationality of $\phi$ implies that $\epso$ may take many different values as $o$ runs through elements of $\Oqa$. Nevertheless, if one only considers orbits having a fixed size, the angles $\epso$ only take at most $6$ different possible values. This will be used in Section~\ref{sect: Analytic estimates spval} to estimate the size of $\spval$.  
\end{remark}

\subsection{Equidistribution of Kloosterman angles}\label{sect: Equid sect}

In this section, we state an equidistribution result concerning the angles of our Kloosterman sums $\Klo$. In Section~\ref{sect: Bounding Sato-Tate sect} we will manipulate the Kloosterman angles $\theto$ associated to orbits in $\Oqa$. First, we define and discuss how Kloosterman angles associated to orbits in $\Pqa$ distribute. Later, we will explain how these are related to the angles $\theto$. Before stating any result, we introduce the necessary theoretical background to discuss equidistribution of Kloosterman angles. \\

Recall from \cite[§4]{Katz88} the construction of the  Kloosterman sheaf $\Klsh := Kl_{2}(1/2)$. This is an  $\ell$-adic lisse sheaf on $\Gm \, / \, \Fq$ of rank $2$ which is pure of weight $0$. Consider the morphism $[2] : \Gm \rightarrow \Gm, x \mapsto x^2$ and define $\F:= [2]^{*} \Klsh$. $\F$ is also a lisse $\Qlbar$-sheaf of rank $2$ on $\Gm$ which is pure of weight $0$ (see \cite{FuLiu}). More concretelly, $\F$ is characterized by the following condition: at any place $v$ of $\Gm \, / \, \Fq$, we have
\begin{equation*}
	\Trpol{\Frv}{\F_{\overline{v}}} = \abs{\Fv}^{-1/2} \, \Kloos{\Fv}{\psi_{\Fv}}{\beta_{v}^2},
\end{equation*}
where $\beta_v \in \Gm(\Fqbar)$ is any element in the closed point of $\Gm$ corresponding to $v$, and $\overline{v}$ is a geometric point lying over $v$. Note that the change of variables $x \mapsto \beta_{v} x$ gives $\Kloos{\Fv}{\psi_{\Fv}}{\beta_{v}^2} = \Kloos{\Fv}{\psi_{\Fv, \beta_{v}}}{1}$, which corresponds to the Kloosterman sums introduced in Section \ref{sect: Orbits}. The lisse sheaf $\F$ gives rise to a $\Qlbar$-representation $\rho_{\F}$ of the arithmetic fundamental group of $\Gm \, / \, \Fq$. Remark 1 in \cite{FuLiu} shows that the Zariski closure $G$ of $\rho_{\F} ( \pi_{1}(\Gm \otimes_{\Fq} \Fqbar))$ in $\GL(2)$ is $G = \SL(2)$ (same as for the sheaf $\Klsh$) and that $\F$ satisfies the assumptions of \cite[§3.1 - §3.3]{Katz88}, namely that $\image(\rho_{\F}) \subset G(\Qlbar)$, \ie, $\rho_{\F}$ has trivial determinant.

Via the isomorphism $\Qlbar \simeq \C$ fixed above, we may view $G$ as an algebraic group over $\C$. The subgroup $K := \SU(2, \C)$ is then maximal compact in $G(\C)$. The set of conjugacy classes of $K$, denoted by $K\nat$, may be identified with the interval $\intrv$. Indeed any matrix $M \in K$ is conjugate (in $K$) to a diagonal matrix $\Diag(\e^{i \, \theta_{M}}, \e^{-i \, \theta_{M}})$ for some $\theta_{M} \in \intrv$, hence $\Trace(M) = 2 \cos \theta_{M}$. Denote by $\nu\nat$ the measure on $K\nat$ defined as the direct image of the Haar measure on $K$, normalised to have total mass $1$. The pushforward of $\nu\nat$ under the identification $K\nat \stackrel{\sim}{\longrightarrow} \intrv$ is the Sato-Tate measure $\nust:= \frac{2}{\pi}\sin^{2}\theta \, \d\theta$ (see \cite[Chap. 13]{Katz88} for details).

We can now define the Kloosterman angle associated to a place. For any place $v \in \Mpl \setminus \lbrace 0, \infty \rbrace$ corresponding to a closed point of $\Gm$, denoted also by $v$, the semi-simplification of $\rho_{\F}(\Frv) \in G(\C)$ is conjugate to an element of $K$. Therefore, we can define the conjugacy class (in $K$) of $\rho_{\F}(\Frv)\semis$, and associate to $v$ (\textit{via} the identification $K\nat \simeq \intrv$) a unique element $\thetv \in \intrv$. We call $\thetv$ the Kloosterman angle associated to the place $v$. For any $\beta \in v$, we define $\theta(\beta) := \thetv$, which is independent of $\beta$ by \refKloos{Kloos bien def} for the angle associated to $v$ as above. Under the identification $K\nat \simeq \intrv$ we have
\begin{equation}\label{eq: thetv to Klv}
	2 \, \cos \thetv = \Trace(\rho_{\F}(\Frv)\semis) = \Trpol{\Frv}{\F_{\overline{v}}} = \abs{\Fv}^{-1/2} \, \Kloos{\Fv}{\psi_{\Fv}}{\beta^2}.
\end{equation}

\begin{remark}\label{rmk: theto to thetv}
	We can relate the angles $\thetv$ to the angles $\theto$ defined in Section~\ref{sect: Special value section}. Let $o \in \Oqa$ and $v:= \Theta(o) \in \Pqa$. Equation \eqref{eq: thetv to Klv} combined with~\refKloos{Kl HasseDav} gives $\theto \equiv \begin{cases}
		\thetv \mod \pi & \text{ if } \abs{o} = \abs{v}, \\
		2 \, \thetv \mod \pi  & \text{ if } \abs{o} = 2 \abs{v}.
	\end{cases}$
\end{remark} \vspace{1em}

We associate to the family of angles $\left( \thetv \right)_{v \in \Pqamax}$ a measure on $K\nat \simeq \intrv$.

\begin{definition}\label{def: measure nua}
	For any $a \geq 1$, denote by $\Pqamax:= \ensemble{v \in \Pqa}{\dv = a}$ the subset of orbits of maximal size. We define $\nua$ to be the discrete probability measure on $\intrv$ given by
	\begin{equation*}
		\nua := \frac{1}{\Pi_{q}(a)} \sum_{v \in \Pqamax} \dlthetv,
	\end{equation*}
	where $\delta \lbrace x \rbrace$ denotes the Dirac delta measure at $x \in \intrv$ and, according to the notation of Section \ref{sect: Orbits}, we have $\Pi_{q}(a) = \abs{\Pqamax}$. 
\end{definition}

We have the equality\footnote{In \cite[Theorem 3.6]{Katz88}, the measure $\nua$ is denoted by $Y_a$ and defined by the right-hand side of \eqref{eq: Relating nua to Katz's Ya}.}
\begin{equation}\label{eq: Relating nua to Katz's Ya}
	\nua = \frac{1}{\abs{\ensemble{\beta \in \Fqa\st }{[\Fq(\beta) : \Fq] = a }}} \sum_{[\Fq(\beta) : \Fq] = a} \delta \left \lbrace \theta(\beta) \right \rbrace.	
\end{equation}
Indeed, $\langle q \rangle$ acts on $\Fqa\st$ as the Frobenius, so one can partition $\ensemble{\beta \in \Fqa\st }{[\Fq(\beta) : \Fq] = a }$ into $\Pi_{q}(a)$ orbits of size $a$. Regrouping elements by orbits, the right-hand side of \eqref{eq: Relating nua to Katz's Ya} may be rewritten as
\begin{equation*}
	\frac{1}{a \, \Pi_{q}(a)} \sum_{v \in \Pqamax} \sum_{\beta \in v} \delta \left \lbrace \theta(\beta) \right \rbrace = \frac{1}{a \, \Pi_{q}(a)} \sum_{v \in \Pqamax} a \, \delta \left \lbrace \thetv \right \rbrace
\end{equation*}
and the right-hand side is clearly equal to $\nu_a$.

We cite a theorem of Katz establishing the equidistribution of Kloosterman angles with respect to the Sato--Tate measure (see \cite[Theorem $3.6$]{Katz88} and \cite{FuLiu}). Recall that a sequence of Borel measures $\left(\mu_i\right)_{i \geq 1}$ on $\intrv$ is said to converge weak--$*$ to $\nust$ if, for every continuous $\C$-valued function $f$ on $\intrv$, the sequence $\int_{\intrv} f \d \mu_i$ converges to  $\int_{\intrv} f \dnust$ as $i \rightarrow \infty$. 

\begin{theorem}\label{thm: Equid thm}
	The sequence of probability measures $(\nua)_{a \geq 1}$ converges weak--$*$ to the Sato--Tate measure $\nust$. 
\end{theorem} 

\begin{proof}[Sketch of the proof]
	It suffices to treat the case of characters of irreducible finite dimensional representations of $K$, since these functions are dense in the space of continuous central functions on $K$. The integral of such a character with respect to the Sato--Tate measure equals $0$. In \cite{Katz88, FuLiu}, the authors show that for any irreducible finite-dimensional representation $\Psi$ of $K$, 
	\begin{equation}\label{eq: Katz bound}
		\abs{\int_{K\nat} \Trace(\Psi) \, \dnua} \ll_q \frac{a \, \dim \Psi}{q^{a/2}}.
	\end{equation}	
	This follows from the Lefschetz trace formula, Deligne's purity theorem and the Grothendieck--Ogg--Shafarevich formula combined to a detailed analysis of the monodromy of Kloosterman sheaves	(see \cite{Katz88, FuLiu} for details).
\end{proof}

To study the size of the special value of $L(S, T)$, we will need an effective version of Theorem~\ref{thm: Equid thm}, which we state below:

\begin{theorem}\label{thm: Equid eff thm}
	For any continuously differentiable function $g$ on $\intrv$, any $a \geq 1$, 
	\begin{equation*}
		\abs{\int_{\intrv} g \, \dnua - \int_{\intrv} g \, \dnust} \ll_{q} \frac{a^{1/2}}{q^{a/4}} \, \int_{0}^{\pi} \abs{g'(x)} \dx.
	\end{equation*}
\end{theorem}

\begin{proof}
	We follow the proof of Proposition $6.11$ in \cite{Gri19}. Recall the definition of the star discrepancy $\Discrep$:
	\begin{equation*}
		\Discrep := \sup_{x \in \intrv} \abs{\int_{[0, x)} \dnua - \int_{[0, x)} \dnust}.
	\end{equation*}	
	This quantity measures ``how far" the angles $\thetv$ are from being perfectly equidistributed with respect to the Sato--Tate measure $\nust$.	Let $g : \intrv \rightarrow \C$ be a continuously differentiable function. Corollary $2$ in \cite{Niederreiter} establishes the analogue of Koksma's inequality in the context of the Sato--Tate measure:
	\begin{equation}\label{eq: Koksma}
		\abs{\int_{\intrv} g \, \dnua - \int_{\intrv} g \, \dnust} \leq \Discrep \int_{0}^{\pi} \abs{g'(x)} \dx.
	\end{equation}
	Let $n \geq 1$ be an integer and write $\Psi_n = \mathrm{Symm}^{n}(\mathrm{std})$ for the $n$-th symmetric power of the standard representation of $K$. We apply \eqref{eq: Katz bound} to $\Psi_n$ and combine it with Lemma $3$ in \cite{Niederreiter} to get the upper bound
	\begin{equation}\label{eq: Bound discrep}
		\Discrep \ll_q \frac{a^{1/2}}{q^{a/4}}.
	\end{equation}
	The result follows by putting \eqref{eq: Koksma} and \eqref{eq: Bound discrep} together. 
\end{proof}

\subsection{Kloosterman angles ``avoid" Gauss angles}\label{sect: Distance angles}
Now we show that the angles of Gauss and Kloosterman sums cannot be equal nor arbitrarily close to each other. This phenomenon was partly observed in \cite{Gri19} or \cite{GriDWit}, where the authors show that Kloosterman angles avoid certain rational multiples of $\pi$ (which are angles of Gauss sums). Our results follow from Liouville's inequality and hold in any odd characteristic. That is why, in this section, we work in broader generality and temporarily drop the assumption that $p \geq 7$. \textbf{Until the end of the section, we fix a prime number $p \geq 3$ and integers $n, N \geq 1$, with $n$ coprime to $p$.}

Consider a finite extension $\Ff$ of $\Fp$, a non-trivial multiplicative character $\chi$ on $\Ff\st$, a non-trivial additive character $\psi$ on $\Ff$ and $\alpha \in \Ff\st$. Recall from \refGauss{Gauss norm} that $\abs{\iota \left(\Gauss{\Ff}{\chi}{\psi} \right)} = \sqrt{\abs{\Ff}}$. We let $\epspm \in \interv$ be the argument of $\iota \left(\Gauss{\Ff}{\chi}{\psi} \right)$, so that 
\begin{equation*}
	\iota \left(\Gauss{\Ff}{\chi}{\psi} \right)= \sqrt{\abs{\Ff}} \, \e^{i \, \epspm}.
\end{equation*}
Similarly, recall that $\Kloos{\Ff}{\psi}{\alpha}$ decomposes as $\kl{\Ff}{\psi}{\alpha} + \klbis{\Ff}{\psi}{\alpha}$, where $\kl{\Ff}{\psi}{\alpha}$, $\klbis{\Ff}{\psi}{\alpha}$ are algebraic integers satisfying $\abs{\kl{\Ff}{\psi}{\alpha}} = \abs{\klbis{\Ff}{\psi}{\alpha}} = \sqrt{\abs{\Ff}}$, such that $\iota \left( \kl{\Ff}{\psi}{\alpha}\right)$, $\iota \left( \klbis{\Ff}{\psi}{\alpha} \right)$ are conjugate. Switching $\kl{\Ff}{\psi}{\alpha}$, $\klbis{\Ff}{\psi}{\alpha}$ if necessary, we may assume that $\iota \left( \kl{\Ff}{\psi}{\alpha} \right)$ lies in the upper half plane, and we let $\thetpm \in \intrv$ denote its argument, so that
\begin{multline*}
	\iota \left( \kl{\Ff}{\psi}{\alpha} \right)= \sqrt{\abs{\Ff}} \, \e^{i \, \thetpm},  \quad \iota \left( \klbis{\Ff}{\psi}{\alpha} \right)= \sqrt{\abs{\Ff}} \, \e^{- i \, \thetpm}  \\  \text{ and } \quad \iota \left( \Kloos{\Ff}{\psi}{\alpha} \right)= 2 \, \sqrt{\abs{\Ff}} \, \cos(\thetpm).
\end{multline*}
We write $R_{N}$ for the set of angles of $N$-th roots of unity and $\GamFn$ for the set of angles of Gauss sums associated to characters of order $n$ on $\Ff\st$:
\begin{align*}
	R_{N} & := \biggl \lbrace \frac{2j \pi}{N} ; \ j = 0, \, \ldots , \, N-1 \biggr \rbrace, \\
	\GamFn & := \biggl \lbrace \epspm \ ; \  \chi : \Ff\st \rightarrow \Qbar\st \text{ of order } n, \  \psi: \Ff \rightarrow \Qbar\st \text{ non-trivial} \biggr \rbrace.
\end{align*}
For any finite subset $I  \subset [0, 2 \pi]$ and any $x \in [0, 2 \pi]$, we write $\dist\left( x, I \right) = \min_{y \in I} \abs{x - y}$.

\begin{theorem}\label{thm: Dist G K angles}
	Let $p \geq 3$ be a prime number and $N, \, n \geq 1$ be integers with $n$ coprime to $p$. There are explicit constants $\tau_{p, N} > 0$ and $\sigma_{p, n} > 0$ satisfying the following: for any finite extension $\Ff / \Fp$, any non-trivial additive character $\psi$ on $\Ff$, any $\alpha \in \Ff\st$ we have
	\begin{enumerate}
		\item $\min \Bigl \lbrace \dist \left( \thetpm, R_N  \right), \, \dist \left(2 \pi - \thetpm, R_N \right) \Bigr \rbrace \geq \abs{\Ff}^{- \tau_{p, N}}$.
		
		\item $\min \Bigl \lbrace \dist \left( \thetpm, \GamFn  \right), \, \dist \left(2 \pi - \thetpm, \GamFn \right) \Bigr \rbrace \geq \abs{\Ff}^{- \sigma_{p, n}}$.
	\end{enumerate}
\end{theorem}
\vspace{1em}

Before proving this theorem, we show that Kloosterman angles are not in $R_N \, \cup \, \GamFn$, \ie, the distances appearing in Theorem \ref{thm: Dist G K angles} are positive.

\begin{lemma}\label{lem: Thetpm not in GamFn}
	Under the hypotheses of the previous theorem, we have
	\begin{enumerate}
		\item $\thetpm \notin R_N$ and $2\pi - \thetpm \notin R_N$,
		\item $\thetpm \notin \GamFn$ and $2\pi - \thetpm \notin \GamFn$.
	\end{enumerate}
\end{lemma}

\begin{proof}
	Just as in Section~\ref{sect: Character sums}, consider a prime ideal $\pfrac$ of $\Zbar$ lying over $p$ and let $\valp$ be the $\pfrac$-adic valuation normalised so that $\valp(\Ff) = 1$ (in the notation of Section~\ref{sect: Orbits} this would correspond to $\Ff = \Fq$). The angle $\thetpm$ (resp. $2 \pi - \thetpm$) belongs to $\pi \Q$ if and only if $\kl{\Ff}{\psi}{\alpha}/ \sqrt{\abs{\Ff}}$ (resp. $\klbis{\Ff}{\psi}{\alpha}/ \sqrt{\abs{\Ff}}$) is a root of unity.  But~\refKloos{Kloos conj valuation} states that $\valp \left(\kappa_{i, \, \Ff}(\psi, \, \alpha) / \sqrt{\Ff} \right) = \pm 1/2$ for $i = 1, 2$, so $\kappa_{i, \, \Ff} \left(\psi, \, \alpha \right)/ \sqrt{\abs{\Ff}}$ cannot be a root of unity, hence proving the first item. To check item $ii)$, consider a multiplicative character $\chi$ on $\Ff\st$ of order $n$, a non-trivial additive character $\psi'$ on $\Ff$. Stickelberger's theorem implies $\valp(\Gauss{\Ff}{\chi}{\psi'}) \in \left[ \frac{1}{n}, 1 - \frac{1}{n}\right]$ (see \refGauss{Stickelberger Griffon}), and with our normalization of $\valp$, we deduce from~\refKloos{Kloos conj valuation} that $\lbrace \valp(\kl{\Ff}{\psi}{\alpha}), \valp(\klbis{\Ff}{\psi}{\alpha}) \rbrace= \lbrace 0, 1 \rbrace$. This implies that $\Gauss{\Ff}{\chi}{\psi'} \neq \kappa_{i, \, \Ff}(\psi, \, \alpha)$ for $i = 1, 2$, hence $\thetpm, 2\pi - \thetpm \notin \GamFn$. 
\end{proof}

We now turn to the proof of Theorem~\ref{thm: Dist G K angles}, the main tool used being a form of Liouville's inequality.

\begin{proof}[Proof of Theorem~\ref{thm: Dist G K angles}]
	Fix a finite extension $\Ff / \Fp$, a multiplicative character $\chi$ on $\Ff\st$ of order $n$, two non-trivial additive characters $\psi$ and $\psi'$ on $\Ff$ and $\alpha \in \Ff\st$. Define three algebraic numbers $g:= \Gauss{\Ff}{\chi}{\, \psi} / \sqrt{\abs{\Ff}}$ and $k_{i} := \kappa_{i, \, \Ff}^{ }(\psi, \, \alpha)/ \sqrt{\abs{\Ff}}$ for $i= 1, 2$. Let $P_N(T) = T^N - 1 \in \Q[T]$, and $\pig \in \Q[T]$ be the minimal polynomial of $g$ over $\Q$. 
	
	\textit{Proof of $i)$:}  We know from Lemma~\ref{lem: Thetpm not in GamFn} that $P_N(k_i) \neq 0$ for $i = 1, 2$. Using the convexity of $z \mapsto \e^{i z}$ we get, for any $j \in \lbrace 0, \, \ldots, \, N-1 \rbrace$, 
	\begin{equation}\label{eq: Ineq PreLiouv rat}
		\abs{\iota \left( P_N(k_1) \right)} = \abs{\e^{i N \thetpm} - \e^{i N \frac{2j \pi}{N}}} \leq N \abs{\thetpm - \frac{2j \pi}{N}},
	\end{equation}
	and similarly
		\begin{equation}\label{eq: Ineq PreLiouv rat2}
		\abs{\iota \left( P_N(k_2) \right)} \leq N \abs{2 \pi - \thetpm - \frac{2j \pi}{N}}.
	\end{equation}	
	On the other hand, Liouville's inequality (see  \cite[Proposition $3.14$]{Waldschmidt}) gives, for $i =1, 2$,
	\begin{equation}\label{eq: Liouville rat}
		\log \abs{\iota(P_N(k_i))} \geq - [\Q(k_i) : \Q] \left(\log \normun{P_N} + N \, \height(k_i) \right).
	\end{equation}
	Here, we have denoted by $\height : \Qbar \rightarrow \R$ the absolute logarithmic Weil height and by $\normun{P_N}$ the sum of the sizes of the coefficients of $P_N$.  Lemma 5.3 in \cite{Gri19} gives the upper bounds $[\Q(k_i) : \Q] \leq 2(p-1)$ and $\height(k_i) \leq \log \sqrt{\abs{\Ff}}$, and we clearly have $\normun{P_N} = 2$. Combining these estimates we obtain
	\begin{equation}\label{eq: Liouv eff rat}
		\log \abs{\iota(P(k_i))} \geq -2(p-1)\left( \log 2 + N \log \sqrt{\abs{\Ff}} \right) \geq -2(p-1)\left(\frac{\log 2}{\log 3} + \frac{N}{2} \right)  \log \abs{\Ff}, 
	\end{equation}
	where the last inequality uses the fact that $\abs{\Ff} \geq 3$. Combining \eqref{eq: Ineq PreLiouv rat} with \eqref{eq: Liouv eff rat} we get
	\begin{equation*}
		\abs{\thetpm - \frac{j \pi}{N}} \geq \frac{\abs{\Ff}^{-2(p-1)(\log 2 / \log 3 + N/2)}}{N} \geq \abs{\Ff}^{-2(p-1)(\log 2 / \log 3 \, + \, N/2) - \log(N)/ \log(3)},
	\end{equation*}
	(the second inequality uses again $\abs{\Ff} \geq 3$). Setting $\tau_{p, N} := 2(p-1)(\log 2 / \log 3 + N/2) + \log(N)/ \log(3)$ gives the desired inequality. Finally, combining \eqref{eq: Ineq PreLiouv rat2} and \eqref{eq: Liouv eff rat} and following the same line of reasoning, we conclude that $\dist\left(2 \pi - \thetpm, R_N \right) \geq \abs{\Ff}^{- \tau_{p, N}}$ too.
	
	\textit{Proof of $ii)$:} We proceed as above, by first noticing that $\pig(k_i) \neq 0$ for $i = 1, 2$. Indeed, for any $\sigma \in \Gal(\Qbar / \Q)$ we have $\sigma \left( \Gauss{\Ff}{\chi}{\psi}\right) = \Gauss{\Ff}{\sigma \circ \chi}{\sigma \circ \psi}$. The character $\sigma \circ \chi$ on $\Ff\st$ has order $n$, so that the argument of $\sigma \left( \Gauss{\Ff}{\chi}{\psi}\right)$ belongs to $\GamFn$ and Lemma~\ref{lem: Thetpm not in GamFn} shows that $\sigma(g) \neq k_i$.
	
	Let $c \in \Z$ be the lcm of the denominators of coefficients of $\pig$ and define $P := c \, \pig \in \Z[T]$. Write $d = \deg( \pig)$ for the degree of $\pig$ and $g_1 = g, g_2, \, \ldots, \, , g_d$ for the Galois conjugates of $g$. We have
	\begin{equation*}
		\abs{\iota(P(k_1))} = \abs{c} \, \abs{\iota(k_1 - g_1)} \, \prod_{j = 2}^{d} \abs{\iota(k_1 - g_j)}.
	\end{equation*}
	Since $\abs{\iota(k_1)} = \abs{\iota(g_j)} = 1$ for all $j$ (\cf~\refGauss{Gauss norm} and~\refKloos{Kloos norm}), the convexity of $z \mapsto \e^{i z}$ implies
	\begin{equation}\label{eq: Ineq preLiouv}
		\abs{\iota(P(k_1))} \leq 2^{d - 1} \abs{c} \,\abs{\e^{i \, \thetpm} - \e^{i \, \epspmbis}} \leq 2^{d -1} \abs{c} \, \abs{\thetpm - \epspmbis},
	\end{equation}
	and, similarly,
	\begin{equation}\label{eq: Ineq preLiouvbis}
		\abs{\iota(P(k_2))} \leq 2^{d -1} \abs{c} \, \abs{2 \pi -\thetpm - \epspmbis}.
	\end{equation}
	In this setting, Liouville's inequality yields:
	\begin{equation}\label{eq: Liouville}
		\log \abs{\iota(P(k_i))} \geq - [\Q(k_i) : \Q] \left(\log \normun{P} + d \, \height(k_i) \right).
	\end{equation}
	Recall that $\Gauss{\Ff}{\chi}{\psi} \in \Q(\zeta_{n p})$ so $g \in \Q(\zeta_{n p}, \sqrt{p})$ and therefore $d \leq 2np$. Using elementary symmetric polynomials, we have $\normun{P} \leq \abs{c} \, \sum_{r = 0}^{d} \binom{d}{r} = \abs{c} \, 2^d$, since all of the $\iota(g_j)$ have absolute value $1$. In order to bound $\abs{c}$, notice that $\abs{\Ff} g = \Gauss{\Ff}{\chi}{\psi} \sqrt{\abs{\Ff}}$ is an algebraic integer so its minimal polynomial $\pi_{\abs{\Ff} g}$ has integral coefficients. On the other hand $\abs{\Ff}^{-d} \pi_{\abs{\Ff} g}(\abs{\Ff} T) \in \Q[T]$ vanishes at all $g_j$, is unitary, and has the same degree as $\pig$. We deduce that
	\begin{equation*}
		\abs{\Ff}^{d} \pig(T) =  \pi_{\abs{\Ff} g}(\abs{\Ff} T) \in \Z[T],
	\end{equation*}
	implying that $\abs{c} \leq \abs{\Ff}^d$. Plugging the bounds above into \eqref{eq: Liouville} gives
	\begin{align*}
		\log \abs{\iota(P(k_i))} & \geq -2(p-1)\left(\log \left((2 \abs{\Ff})^d \right) +2np \log \sqrt{\abs{\Ff}} \right) \\
		& \geq -2(p-1)\left(2np \log \left(2 \abs{\Ff}  \right) +np \log \abs{\Ff} \right) \\
		& \geq - 2 np (p-1)\left( 3 \log\abs{\Ff} + 2 \log 2 \right) \\
		& \geq - 2 np (p-1)\left( 3 + \frac{2 \log 2}{\log 3} \right) \log \abs{\Ff} \qquad \text{as } \abs{\Ff} \geq 3.
	\end{align*}
	Combining this with \eqref{eq: Ineq preLiouv} gives
	\begin{equation*}
		\abs{\thetpm - \epspmbis} \geq \frac{\abs{\Ff}^{-2np(p-1)(3 + 2\log 2 / \log 3)}}{2^{2np} \abs{\Ff}^{2np}} \geq \abs{\Ff}^{- \sigma_{p, n}},
	\end{equation*}
	where $\sigma_{p, n} := 2np \left( \frac{\log 2}{\log 3} + (p-1) \left(3 + \frac{2 \log 2}{\log 3}\right) +1 \right)$ (again we used the fact that $\abs{\Ff} \geq 3$). In a similar way, we also have $\abs{2 \pi - \thetpm - \epspm} \geq \abs{\Ff}^{- \sigma_{p, n}}$.
\end{proof}

Here is the special case of Theorem \ref{thm: Dist G K angles} that we need in Section \ref{sect: Analytic estimates spval}.

\begin{corollary}\label{cor: Dist G K cor}
	Recall the notation from Section~\ref{sect: Special value section} and define $\sigma_{p} := \sigma_{p, 3}$. For any $o, o' \in \Oqa$ such that $\abs{o} = \abs{o'}$ we have
	\begin{equation*}
		\min \left(\abs{\theto - \epsoprim}, \abs{\theto  + \epsoprim}, \abs{2\pi - \theto - \epsoprim}, \abs{2\pi - \theto + \epsoprim}\right)\geq q^{- \sigma_{p} \abs{o}}.
	\end{equation*}
\end{corollary}

\begin{proof}
	Let $\Fo = \Ff_{o'}$ be the unique extension of $\Fq$ of degree $\abs{o}$, so that $\abs{\Fo} = q^{\abs{o}}$. Since $\epsoprim \in \Gamma_{\Fo, 3}$, we deduce from Theorem \ref{thm: Dist G K angles} that $\abs{\theto - \epsoprim} \geq q^{- \sigma_{p} \abs{o}}$ and $\abs{2\pi - \theto - \epsoprim} \geq q^{- \sigma_{p} \abs{o}}$. In order to show the two other inequalities, we simply need to check that $-\epsoprim \in \Gamma_{\Fo, 3}$. Fix a representative $(j', \alpha') \in o'$. By definition, $\eps_{o'} = \eps_{\Fo}^{ }(\chi_{\Fo,3}^{j'}, \psi_{\Fo, \alpha'})$ and since $\chi_{\Fo, 3}$ is of order $3$, we have $-\epsoprim = \eps_{\Fo}^{ }(\chi_{\Fo,3}^{-j'}, \psi_{\Fo, \alpha'}) \in \Gamma_{\Fo, 3}$, hence the result.
\end{proof}


\section{Analytic estimates on the special value}\label{sect: Analytic estimates spval}

In this section we estimate the asymptotic behaviour of the special value of $L(S, T)$ as the parameter $a$ tends to $\infty$. To make the dependence on $a$ more explicit, we write $\Sa$ with the subscript $a$ instead of just $S$. The goal of this section is to prove the following theorem:

\begin{theorem}\label{thm: Ratio spval height thm}
	Let $\Fq$ be a finite field of characteristic $p \geq 7$, $\K := \Fq(t)$. For any $a \geq 1$, let $C_a$ be the hyperelliptic curve defined by \eqref{eq: Hyp model}, and $S_a := \Jac(C_a)$ be its Jacobian. There are positive constants $c_1, c_2$, depending at most on $q$, such that for any $a \geq 1$, we have 
	\begin{equation}\label{eq: Ratio spval height}
		- \frac{c_1}{a} \leq \frac{\log \spvala}{\log H(\Sa)} \leq \frac{c_2}{a}.
	\end{equation}
\end{theorem}

\begin{remark}\label{rmk: Ratio H into qa}
	By definition we have $\log H(\Sa) = h(\Sa) \log q$. Recall from Proposition~\ref{prop: Height} that $\log H(\Sa)$ and $q^a$ have the same order of magnitude. Thus, Theorem~\ref{thm: Ratio spval height thm} is equivalent to the following statement: there are positive constants $c_{1}', c_{2}'$ depending at most on $q$ such that $- \frac{c_{1}'}{a} \leq \frac{\log \spvala}{q^a} \leq \frac{c_{2}'}{a}$. We will study the ratio $\log \spvala / q^{a}$ rather than $\log \spvala / \log H(\Sa)$.
\end{remark}

Let us estimate the logarithm of the special value. Recall from Theorem~\ref{thm: Special value thm} that we have, for any $a \geq 1$,
\begin{equation}\label{eq: Spvala restated}
	\spvala = \prod_{o \, \in \Oqa} 4 \, \abs{\sin \left( \frac{\epso + \theto}{2}\right) \, \sin \left( \frac{\epso - \theto}{2}\right)}.
\end{equation}
The special value is a positive real number, so applying the logarithm to \eqref{eq: Spvala restated} and dividing by $q^a$ gives
\begin{equation}\label{eq: Spval over q^a}
	\frac{\log \spval}{q^a} =  \log 4 \, \frac{\abs{\Oqa}}{q^a} + \frac{1}{q^a} \, \sum_{o \, \in \Oqa} \log \abs{\sin \left( \frac{\epso + \theto}{2}\right) \, \sin \left( \frac{\epso - \theto}{2}\right)}.
\end{equation}

We can already prove the upper bound in Theorem \ref{thm: Ratio spval height thm}.

\begin{proof}[Proof of the upper bound in Theorem~\ref{thm: Ratio spval height thm}]
	For any orbit $o \in \Oqa$, we have the trivial upper bound $\log \abs{\sin \left( \frac{\epso + \theto}{2}\right) \, \sin \left( \frac{\epso - \theto}{2}\right)} \leq 0$, which implies $\frac{\log \spval}{q^a} \leq \frac{\abs{\Oqa}}{q^a} \, \log 4$. We conclude by using the estimate $\abs{\Oqa} \ll_q q^a/a$ given in Lemma~\ref{lem: Estimates numb orbits} $iii)$.
\end{proof}

We devote the remainder of the section to proving the lower bound in Theorem~\ref{thm: Ratio spval height thm}. Since $\log 4 \abs{\Oqa}/q^a = O(1/a)$, we focus on the sum appearing in the right-hand side of \eqref{eq: Spval over q^a}.

\subsection{Reduction to maximal orbits}

In this section, we introduce the notion of maximal orbits and we show that the contribution of non-maximal orbits to the logarithm of the size of the special value is negligible.

\begin{definition}
	For any $a \geq 1$, define the set of maximal orbits 
	\begin{equation*}
		\Oqamax := \ensemble{o \in \Oqa}{\abs{o} \text{ is maximal among orbits in } \Oqa}.
	\end{equation*}
	Any orbit $o \in \Oqamax$ has length
	\begin{equation*}
		\abs{o} = 
		\begin{cases}
			a  & \text{ if } q \equiv 1 \mod 3 \text{ or } a \text{ is even}, \\
			2a  & \text{ if } q \equiv 2 \mod 3 \text{ and } a \text{ is odd}. 
		\end{cases}
	\end{equation*}
\end{definition}

\begin{remark}\label{rmk : Relating Oqamax & Pqamax}
	We defined in Section \ref{sect: Equid sect} the notion of a maximal orbit in $\Pqa$. Using $\Theta$, we can relate maximal orbits in $\Oqa$ to maximal orbits in $\Pqa$ as follows. Recall from Section \ref{sect: Orbits} that there is a map $\Theta : \Oqa \rightarrow \Pqa$ induced by the natural map $\left( \Z / 3 \Z \right)\st \times \Fqa\st \rightarrow \Fqa\st$. Fix $o \in \Oqa$ and let $v := \Theta(o) \in \Pqa$.
	\begin{enumerate}
		\item When $q \equiv 1 \mod 3$, we have $o \in \Oqamax \Leftrightarrow v \in \Pqamax$, in which case $\abs{o} = \abs{v} = a$ and $\abs{\Theta^{-1}(v)} = 2$.
		\item When $q \equiv 2 \mod 3$ and $a$ is odd, we have $o \in \Oqamax \Leftrightarrow v \in \Pqamax$, in which case $\abs{o} = 2a$, $\abs{v} = a$, and $\abs{\Theta^{-1}(v)} = 1$.
		\item When $q \equiv 2 \mod 3$ and $4 \divi a$, we have $o \in \Oqamax \Leftrightarrow v \in \Pqamax$, in which case $\abs{o} = \abs{v} = a$ and $\abs{\Theta^{-1}(v)} = 2$.
		\item When $q \equiv 2 \mod 3$ and $2$ divides $a$ exactly once, we have $o \in \Oqamax \Leftrightarrow \abs{v} \in \left \lbrace a, \frac{a}{2} \right \rbrace$. If $\abs{v} = a$ then $\abs{\Theta^{-1}(v)} = 2$ and if $\abs{v} = \frac{a}{2}$ then $\abs{\Theta^{-1}(v)} = 1$.
	\end{enumerate}
	All these results are a direct application of \eqref{eq: Size orbit eqt} and \eqref{eq: Preimages by Theta}.
\end{remark}

Using the definition of the orbits given in Section \ref{sect: Orbits}, we have $\sum_{o \in \Oqa} \abs{o} = \abs{\Ia} = 2(q^{a} -1)$. As a consequence of the Prime Number Theorem for $\Fq[t]$, Lemma \ref{lem: Estimates numb orbits} $iii)$ states that $\abs{\Oqa}$ has order of magnitude $q^{a}/a$. Thus, the average length of an orbit is $a$, so that, on average, orbits in $\Oqa$ are as large as they can be.	On the other hand, using the above remark, we estimate the typical size of a non-maximal orbit:

\begin{lemma}\label{lem: Sum nonmax sizes}
	For any $a \geq 1$, 
	\begin{equation*}
		\sum_{o \in \Oqa \setminus \Oqamax} \abs{o} \ll_{q} q^{a/2}.
	\end{equation*}
\end{lemma}

\begin{proof}
	We distinguish cases according to $q \mod 3$ and the parity of $a$.
	
	When $q \equiv 1 \mod 3$, we combine the estimate $\Pi_{q}(k) \leq q^k/k$ from Lemma~\ref{lem: Estimates numb orbits} $i)$ with the first case in Remark~\ref{rmk : Relating Oqamax & Pqamax} to get 	
	\begin{multline*}
		\sum_{o \in \Oqa \setminus \Oqamax} \abs{o} = \sum_{\substack{1 \leq k < a \\ k \divi a}} \sum_{\substack{o \in \Oqa \\ \abs{\Theta(o)} = k}} \abs{o} \, = \sum_{\substack{1 \leq k < a \\ k \divi a}} 2 \sum_{\substack{v \in \Pqa \\ \abs{v} = k}} \abs{v} \\
		= 2 \sum_{\substack{1 \leq k \leq \left \lfloor a/2 \right \rfloor  \\ k \divi a}} k \, \Pi_{q}(k) \,\leq 2 \sum_{1 \leq k \leq a/2} q^k \ll_q \, q^{a/2}.
	\end{multline*}
	The third equality uses the fact that a divisor of $a$ strictly smaller than $a$ is smaller than $\left \lfloor \frac{a}{2} \right \rfloor$.
	
	When $q \equiv 2 \mod 3$ and $a$ is odd then $v \in \Pqa$ is not maximal if and only if $\abs{v} \leq \left \lfloor a/3 \right \rfloor$ (the smallest integer $> 1$ that could divide $a$ is 3). Combining Lemma~\ref{lem: Estimates numb orbits} $i)$ with Remark~\ref{rmk : Relating Oqamax & Pqamax} $ii)$ we get
	\begin{equation*}
		\sum_{o \in \Oqa \setminus \Oqamax} \abs{o}  =  \sum_{\substack{1 \leq k \leq \left \lfloor a/3 \right \rfloor \\ k \divi a}} \sum_{\substack{o \in \Oqa \\ \abs{\Theta(o)} = k}} \abs{o} = \sum_{\substack{1 \leq k \leq \lfloor a/3 \rfloor \\ k \divi a}} \sum_{\substack{v \in \Pqa \\ \abs{v} = k}} 2 \abs{v} \leq \sum_{\substack{1 \leq k \leq a/3 \\ k \divi a}} 2k \, \Pi_{q}(k) \, \ll_q \, q^{a/3}.
	\end{equation*}
	When $q \equiv 2 \mod 3$ and $a$ is even, again, a divisor of $a$ is strictly smaller than $a$ if and only if it is smaller than $a/2$. Following the same line of reasoning as above we have
	\begin{multline*}
		\sum_{o \in \Oqa \setminus \Oqamax} \abs{o} = \sum_{\substack{1 \leq k \leq a/2 \\ k \divi a \\ k \text{ even}}} \sum_{\substack{o \in \Oqa \\ \abs{\Theta(o)} = k}} \abs{o} + \sum_{\substack{1 \leq k \leq a/2 \\ k \divi a \\ k \text{ odd}}} \sum_{\substack{o \in \Oqa \\ \abs{\Theta(o)} = k}} \abs{o} \\
		= \sum_{\substack{1 \leq k \leq a/2 \\ k \divi a \\ k \text{ even}}} 2 \sum_{\substack{v \in \Pqa \\ \abs{v} = k}} \abs{v} + \sum_{\substack{1 \leq k \leq 	a/2 \\ k \divi a \\ k \text{ odd}}} \sum_{\substack{v \in \Pqa \\ \abs{v} = k}} 2 \abs{v} \ll_q q^{a/2}.
	\end{multline*}
	In all three cases, independently of $q \mod 3$ and the parity of $a$, the sum $\sum_{o \in \Oqa \setminus \Oqamax} \abs{o}$ is indeed bounded from above by $O \left( q^{a/2} \right)$.
\end{proof}

The interest of introducing $\Oqamax$ is that restricting to the ``subfamily" of angles $(\epso)_{o \in \Oqamax}$ reduces the amount of possible values these angles can take (\cf \ Remark~\ref{rmk: Possible gauss angles}). As the result below proves, the contribution of the non-maximal orbits to the sum appearing in \eqref{eq: Spval over q^a} is negligible:

\begin{proposition}\label{prop: Maj sum orbs nonmax}
	As $a \rightarrow \infty$, we have 
	\begin{equation*}
		\abs{\frac{1}{q^a} \sum_{o \, \in \Oqa \setminus \Oqamax} \log \abs{\sin \left( \frac{\epso + \theto}{2}\right) \, \sin \left( \frac{\epso - \theto}{2}\right)}} \ll_q \frac{1}{q^{a/2}}.
	\end{equation*}
\end{proposition}

\begin{proof}	
	For any $o \in \Oqa$, we have by Lemma~\ref{lem: Thetpm not in GamFn} $ii)$, $\frac{\epso + \theto}{2}, \frac{\epso - \theto}{2} \notin \lbrace 0, \pi \rbrace$, implying the inequalities $0 < \abs{\sin \left( \frac{\epso + \theto}{2}\right) \sin \left( \frac{\epso - \theto}{2}\right)} \leq 1$ for any $o \in \Oqa$. All the terms in the sum are non-positive, so
	\begin{multline*}
		\abs{\sum_{o \, \in \Oqa \setminus \Oqamax} \log \abs{\sin \left( \frac{\epso + \theto}{2}\right) \, \sin \left( \frac{\epso - \theto}{2}\right)}} \\
		=  \sum_{o \, \in \Oqa \setminus \Oqamax} \log  \left(\abs{\sin \left( \frac{\epso + \theto}{2}\right) \, \sin \left( \frac{\epso - \theto}{2}\right)}^{-1} \right).
	\end{multline*}
	Let $o \in \Oqa$, recall that $\epso$ belongs to the interval $\interv$ and $\theto$ belongs to $\intrv$ so $\frac{\epso + \theto}{2}, \frac{\epso - \theto}{2}$ both belong to $\left[- \frac{\pi}{2}, \frac{3 \pi}{2} \right) \, \setminus \lbrace 0, \pi \rbrace$. For any $x \in \left[- \frac{\pi}{2}, \frac{3 \pi}{2} \right) \,  \setminus \lbrace 0, \pi \rbrace$, we have $\frac{1}{\abs{\sin (x)}} \leq \frac{\pi^2}{\abs{x (\pi - x)}}$. Applied to $x = \frac{\epso + \theto}{2}$ and $x = \frac{\theto - \epso}{2}$, this upper bound gives
	\begin{equation*}
		\abs{\sin \left( \frac{\epso + \theto}{2}\right) \, \sin \left( \frac{\epso - \theto}{2}\right)}^{-1} \leq \frac{(2 \pi)^4}{\abs{(\epso + \theto) \, ( 2\pi - \theto - \epso) \, (\epso - \theto) \, (2 \pi - \theto + \epso)}}.
	\end{equation*}
	Using Corollary~\ref{cor: Dist G K cor} we get
	\begin{equation*}
		\log  \left( \abs{\sin \left( \frac{\epso + \theto}{2}\right) \, \sin \left( \frac{\epso - \theto}{2}\right)}^{-1} \right) \leq \log  \left( \frac{(2 \pi)^4}{ q^{- 4 \, \sigma_{p} \abs{o}}} \right) \ll_q \abs{o}.
	\end{equation*}
	Summing over non-maximal orbits $o \in \Oqa \setminus \Oqamax$, we obtain
	\begin{equation*}
		\abs{\frac{1}{q^{a}} \sum_{o \in \Oqa \setminus \Oqamax} \log \abs{\sin \left( \frac{\epso + \theto}{2}\right) \, \sin \left( \frac{\epso - \theto}{2}\right)}} \leq \frac{1}{q^{a}} \sum_{o \in \Oqa \setminus \Oqamax} \abs{o},
	\end{equation*}
	and Lemma \ref{lem: Sum nonmax sizes} allows to conclude.
\end{proof}

\subsection{Lower bounds on the sum over maximal orbits}

In this section we give lower bounds on $\frac{1}{q^a} \, \sum_{o \, \in \Oqamax} \log \abs{\sin \left( \frac{\epso + \theto}{2}\right) \, \sin \left( \frac{\epso - \theto}{2}\right)}$, which is the term we still need to control in order to establish Theorem \ref{thm: Ratio spval height thm}. 

We begin by bounding the term inside the logarithm. Let us briefly recall the description of Gauss angles given in Proposition \eqref{prop: Possible eps angles}. When $q \equiv 1 \mod 3$, there is a map $\prun: \Oqa \rightarrow \lbrace \pm 1 \rbrace$ lifting the first coordinate of any orbit in $\Oqa$. Moreover, in this case, there is an irrational number $\phi \in \interv$ such that for any $o \in \Oqamax$, the angle $\epso$ is of the form
\begin{equation}
	\epso = 
	\begin{cases}
		\prun(o) a \phi + \frac{2 j \pi}{3} & \text{ with } j \in \lbrace 0, 1, 2 \rbrace \text{ if } q \equiv 1 \mod 3, \\
		\frac{j \pi}{3} & \text{ with } j \in \lbrace 0, \, \ldots, \, 5 \rbrace \text{ if } q \equiv 2 \mod 3.
	\end{cases}
\end{equation}
As a simple exercise in real analysis reveals, we have, for any $x \in \interv$, 
\begin{equation}\label{eq: Minor sinus}
	\left( \frac{\sin \left( 6x \right)}{4} \right)^2 \leq \min_{j = \, 0, \, \ldots, \, 5} \abs{ \sin \left( x + \frac{j \pi}{6}\right)}.
\end{equation}
Define $\lambda_{a} := a \phi$ if $q \equiv 1 \mod 3$ and $\lambda_{a} := 0$ if $q \equiv 2 \mod 3$. We can write, independently of $q \mod 3$ and the parity of $a$, $\epso = \prun(o) \lambda_{a} + j \pi /3$, for some $j \in \lbrace 0, \, \ldots, \, 5 \rbrace$. With this notation
\begin{equation*}
	\frac{\epso + \theto}{2} = \frac{\theto \pm \lambda_{a}}{2} + \frac{j \pi}{6} \quad \text{ and } \quad \frac{\epso - \theto}{2} = \frac{-\theto \pm \lambda_{a}}{2} + \frac{j \pi}{6} 
\end{equation*}
where $j \in \lbrace, 0, \, \ldots, \, 5 \rbrace$ and the sign before $\lambda_{a}$ is $\prun(o)$ for both terms. Since $y \mapsto (\sin y)^2$ is even and $\pi$-periodic, applying inequality \eqref{eq: Minor sinus} to $x = \left(\theto \pm  \lambda_{a} \right)/2$ we get
\begin{equation}\label{eq: Get rid of epso}
	\abs{\sin \left( \frac{\epso + \theto}{2}\right) \, \sin \left( \frac{\epso - \theto}{2}\right)} \geq \left( \frac{\sin \left(3 (\theto + \lambda_{a})\right)}{4} \right)^2 \left( \frac{\sin \left( 3(\theto - \lambda_{a})\right)}{4} \right)^2.
\end{equation}
The right-hand side is positive, because Proposition \ref{prop: Possible eps angles} ensures that any $\pm \lambda_{a} + j \pi /3$ for $j \in \lbrace 0, \, \ldots, \, 5 \rbrace$ is the angle $\epsoprim$ of a Gauss sum $\gamma(o')$ for some $o' \in \Oqamax$. Thus Lemma \ref{lem: Thetpm not in GamFn} guarantees that $3 (\theto \pm \lambda_{a}) \not \equiv 0 \mod \pi$. 
\begin{definition}
	For any $r \in \Z \setminus \lbrace 0 \rbrace$, any $\lambda \in \intrv$, define a function $\Wlam : \intrv \rightarrow (- \infty, 0]$ by
	\begin{equation*}
		\Wlam (x) = \begin{cases}
				\log \abs{\sin(r(x + \lambda))} & \text{ if } r(x + \lambda) \not\equiv 0 \mod \pi, \\
				0 & \text{ if } r(x + \lambda) \equiv 0 \mod \pi.
		\end{cases}
	\end{equation*}
\end{definition}
\noindent Since both terms in \eqref{eq: Get rid of epso} are positive, we use the definition above to write
\begin{multline}\label{eq: Lower bound sum Wlam}
	\sum_{o \in \Oqamax} \log \abs{\sin \left( \frac{\epso + \theto}{2}\right) \, \sin \left( \frac{\epso - \theto}{2}\right)} \\
	\geq 2 \sum_{o \in \Oqamax} \left(W^{3, \lambda_{a}}(\theto) + W^{3, -\lambda_{a}}(\theto) \right) - 4 \log 4 \abs{\Oqamax}.
\end{multline}
The sum over $o \in \Oqamax$ in the right-hand side of \eqref{eq: Lower bound sum Wlam} may be interpreted as an integral over $\intrv$. Recall indeed from Definition \ref{def: measure nua} that for any $a \geq 1$, the measure $\nua$ is given by $\nua := \frac{1}{\Pi_{q}(a)} \sum_{v \in \Pqamax} \dlthetv$.

\begin{proposition}\label{prop: Lower bound by integrals}
	For any $a \geq 1$, let $\lambda_{a}$ be as above and define $r_{a}:= 3$ if $q \equiv 1 \mod 3$ or $a$ is even, and $r_{a} = 6$ if $q \equiv 2 \mod 3$ and $a$ is odd. Define also $\varpi(a):= 1$ if $q \equiv 2 \mod 3$ and $a$ is divisible exactly once by $2$, and $\varpi(a) := 0$ if not. Then we have
	\begin{multline*}
		\frac{1}{q^{a}} \sum_{o \in \Oqamax} \log \abs{\sin \left( \frac{\epso + \theto}{2}\right) \, \sin \left( \frac{\epso - \theto}{2}\right)} \\
		\geq \frac{4 \Pi_{q}(a)}{q^{a}} \int_{\intrv} W^{r_{a}, \lambda_{a}} + W^{r_{a}, - \lambda_{a}} \,  \dnua \ + \ \varpi(a) \, \frac{4 \Pi_{q}\left( \lfloor a / 2  \rfloor \right)}{q^{a}} \int_{\intrv} W^{r_{\lfloor a/2 \rfloor} , \lambda_{ \lfloor a/2 \rfloor} } \, \d \nu_{\lfloor a / 2 \rfloor } + O \left(\frac{1}{a}\right).
	\end{multline*}
\end{proposition}

\begin{proof}
	We start by rewriting the sum in the right-hand side of \eqref{eq: Lower bound sum Wlam} over $o \in \Oqamax$ as a sum (or sums) over $v \in \Pqamax$. To do so, we use the case-by-case analysis done in Remark \ref{rmk : Relating Oqamax & Pqamax}. Remark \ref{rmk: theto to thetv} states that for any $o \in \Oqamax$ letting $v := \Theta(v)$ we have $\theto \equiv \frac{\abs{o}}{\abs{v}} \, \thetv \mod \pi$. Treating separately each of the different cases in Remark \ref{rmk : Relating Oqamax & Pqamax} gives:
	\begin{enumerate}
		\item if $q \equiv 1 \mod 3$, then $\sum_{o \, \in \Oqamax} W^{3, \lambda_{a}}(\theto) = 2 \, \sum_{v \, \in \Pqamax} W^{3, \lambda_{a}}(\thetv)$, 
		\item if $q \equiv 2 \mod 3$ and $a$ is odd, then $\sum_{o \, \in \Oqamax} W^{3, 0}(\theto) = \sum_{v \, \in \Pqamax} W^{6, 0}(\thetv)$, 
		\item if $q \equiv 2 \mod 3$ and $4 \divi a$, then $\sum_{o \, \in \Oqamax} W^{3,0}(\theto) = 2 \, \sum_{v \, \in \Pqamax} W^{3, 0}(\thetv)$,
		\item if $q \equiv 2 \mod 3$ and $2$ divides $a$ exactly once, then 
		\begin{equation*}
			\sum_{o \, \in \Oqamax} W^{3, 0}(\theto) = 2 \, \sum_{v \, \in \Pqamax} W^{3, 0}(\thetv) + \sum_{\substack{v \in \Pqa \\ \abs{v} = a/2}} W^{6, 0}(\thetv).
		\end{equation*}
	\end{enumerate}
	In the first case, the same formula holds if one replaces $W^{3, \lambda_{a}}$ by $W^{3, -\lambda_{a}}$. In the last three cases, $\lambda_{a} = 0$, so we simply have $W^{3, \lambda_{a}} + W^{3, - \lambda_{a}} = 2 \, W^{3, 0}$. Using the definition of $r_a$ given in the proposition, we obtain in cases $i), \, ii)$ and $iii)$
	\begin{equation}\label{eq: Sumo to sumv Wlam}
		\sum_{o \in \Oqamax} \left(W^{3, \lambda_{a}}(\theto) + W^{3, -\lambda_{a}}(\theto) \right) \geq 2 \sum_{v \in \Pqamax} \left(\Wlama(\thetv) + W^{r_{a}, -\lambda_{a}}(\thetv) \right).
	\end{equation}
	In the fourth case we have the extra term $\sum_{\abs{v} = a/2} W^{6, 0}(\thetv) = \sum_{\abs{v} = a/2} W^{r_{a/2}, \lambda_{a/2}}(\thetv)$. Using the term $\varpi(a)$ defined in the proposition and dividing by $q^{a}$, we get for any value of $a$
	\begin{multline*}
		\frac{1}{q^{a}} \sum_{o \in \Oqamax} \left(W^{3, \lambda_{a}}(\theto) + W^{3, -\lambda_{a}}(\theto) \right) \\
		\geq \frac{2}{q^{a}} \sum_{v \in \Pqamax} \left(\Wlama(\thetv) + W^{r_{a}, -\lambda_{a}}(\thetv) \right) + \varpi(a) \, \frac{2}{q^{a}} \sum_{\abs{v} = \lfloor a/2 \rfloor } W^{6, 0}(\thetv) \\
		= \frac{2 \Pi_{q}(a)}{q^{a}} \int_{\intrv} \Wlama + W^{r_{a}, \lambda_{a}} \, \dnua \ + \varpi(a) \, \frac{2 \Pi_{q}\left(\lfloor a / 2 \rfloor \right)}{q^{a}} \int_{\intrv} W^{6, 0} \, \d \nu_{\lfloor a / 2 \rfloor },
	\end{multline*}
	where the equality uses the definition of $\nua$. Combining this with \eqref{eq: Lower bound sum Wlam} and using the estimate $4 \log 4 \abs{\Oqamax} / q^{a} = O(1/a)$ gives the desired result.
\end{proof}

\subsection{Bounding a sequence of integrals}\label{sect: Bounding Sato-Tate sect}

In this section, we bound the integrals appearing in Proposition \ref{prop: Lower bound by integrals}.

\begin{proposition}\label{prop: Seq intg bounded}
	Keep notation as in the previous section and choose $\epsilon \in \lbrace \pm 1 \rbrace$. Then the sequence of integrals $\left( \int_{\intrv} W^{r_{a}, \, \epsilon \lambda_{a}} \, \dnua \right)_{a \geq 1}$ is bounded.
\end{proposition}

If we were able to apply the equidistribution results from Section \ref{sect: Equid sect} to $\Wlama$, then Proposition~\ref{prop: Seq intg bounded} would easily follow. Indeed, Theorem~\ref{thm: Equid thm} would imply that the sequence of integrals appearing in Proposition~\ref{prop: Seq intg bounded} converges and is therefore bounded. However, the function $\Wlama$ is not continuous on $\intrv$, so we cannot directly apply Theorem \ref{thm: Equid thm}. Instead, we approximate $\Wlama$ by continuously differentiable functions on $\intrv$ whose integral with respect to $\nua$ we control. 

\vspace{0.75em}

For $r \in \Z \setminus \lbrace 0 \rbrace$ and $\lambda \in \intrv$, let $\Disc := \Bigl \lbrace \frac{j \pi}{r} - \lambda \, ; \, j \in \Z \Bigr \rbrace  \cap \intrv$. The function $\Wlam$ is $\frac{\pi}{r}$-periodic and has $r$ points of discontinuity, namely the elements of $\Disc$. Moreover, $\Wlam$ is symmetric around each of the elements of $\Disc$, meaning that if $x_{0} \in \Disc$, then
\begin{equation}\label{eq: Symmetry around discont pts}
	 \Wlam(x_0 - x) = \Wlam(x_{0} + x) \qquad \text{ for any } x \in \left[ 0, \pi \right].
\end{equation}
We first give an explicit bound on the integral of $\Wlam$ with respect to the Sato--Tate measure:

\begin{lemma}\label{lem: Bound intg Wlam lemma}
	For any $r \in \Z \setminus \lbrace 0 \rbrace$ and any $\lambda \in \intrv$, we have 
	\begin{equation*}
		\abs{\int_{\intrv}\Wlam \, \dnust} \leq \log 4.
	\end{equation*}
\end{lemma}

\begin{proof}
	Using the trivial upper bound $\sin^{2} \leq 1$, the fact that $\Wlam$ is non-positive on $\intrv$ and doing the change of variables $u = r(x + \lambda)$ gives
	\begin{multline*}
		\abs{\int_{\intrv}\Wlam \, \dnust} = \abs{\frac{2}{\pi} \, \int_{0}^{\pi} \log \abs{\sin (r(x+\lambda)) \sin^{2}(x)}\d x} \\ 
		\leq - \frac{2}{r \pi} \int_{r \lambda}^{r \pi + r \lambda} \log \abs{\sin(u)} \d u = - \frac{2}{\pi} \int_{0}^{\pi} \log \abs{\sin(u)} \d u,
	\end{multline*}
	where the last equality follows by $\pi$-periodicity of $u \mapsto \log \abs{\sin(u)}$. Now it is well-known that $\int_{0}^{\pi} \log (\sin(u)) \d u = -\pi \log 2$ (see \cite[$4.224.3$]{GradRyz}), hence the result. 
\end{proof}

We now construct continuously differentiable functions on $\intrv$ approximating $\Wlam$ and vanishing around the points of discontinuity of $\Wlam$. Let $\beta_{0} : [0, 1] \rightarrow [0, 1]$ be an infinitely differentiable function such that $\beta_{0}(x) = 0$ for $x \in [0, \frac{1}{3}]$ and $\beta_{0}(x) = 1$ for $x \in [\frac{2}{3}, 1]$. We explicitly choose $\beta_0$ as follows: for any $x \in (0, 1)$, let $h(x) := \exp(-1/x)$ and $k(x) := h(x)/(h(1-x) + h(x))$. For any $x \in (1/3, 2/3)$ we let $\beta_{0}(x) := k(3x -1)$. For any $\eta \in \left(0, \frac{\pi}{2} \right)$, consider the function $\beta_{\eta} : \intrv \rightarrow [0, 1]$ defined as follows. For any $x \in \intrv$, let
\begin{equation*}
	\beta_{\eta}(x) := 
	\begin{cases}
		\beta_{0} \left(\frac{x}{\eta} \right) & \text{ if } x \in [0, \eta], \\
		1  & \text{ if } x \in [\eta, \, \pi - \eta ], \\
		\beta_{0} \left( \frac{\pi - x}{\eta} \right) & \text{ if } x \in [\pi - \eta, \, \pi ].
	\end{cases}
\end{equation*}
Consider now $r \in \Z \setminus \lbrace 0 \rbrace$, $\lambda \in \intrv$, and $\eta \in \left(0, \frac{\pi}{2r}\right)$. Let $\beteta : \intrv \rightarrow [0, 1]$ by $\beteta(x) := \beta_{r \eta}(r(x + \lambda))$ and $\Delteta := \ensemble{x \in \intrv}{\dist \left( x, \Disc\right) \geq \eta}$, we have $\beteta(x) = 1$ for any $x \in \Delteta$. 

Define $\Weta := \Wlam \, \beteta$, this is a continuously differentiable function on $\intrv$ which coincides with $\Wlam$ on $\Delteta$ and which vanishes in the neighbourhood of elements of $\Disc$. More precisely, for any $x_{0} \in \Disc$ and any $x \in \left[ x_{0} - \eta /3, x_{0} + \eta /3\right]$, we have $\Weta(x) = 0$. In addition to that, $\Weta$ is also $\frac{\pi}{r}$-periodic and symmetric around each of the points of discontinuity of $\Wlam$, meaning that equation \eqref{eq: Symmetry around discont pts} holds if one replaces $\Wlam$ by $\Weta$.

\begin{lemma}\label{lem: Analytic bounds intg}
	Let $r \in \Z \setminus \lbrace 0 \rbrace$, and $\lambda \in \intrv$. There is a positive constant $c(r)$ depending only on $r$ such that, for any $\eta \in \left(0, \frac{\pi}{2r} \right)$, we have
	\begin{equation*}
		\displaystyle \int_{0}^{\pi} \abs{\Wprim(x)} \dx \leq c(r) \abs{\log \eta}.
	\end{equation*}
\end{lemma}

\begin{proof}
	First of all, notice that $\Wprim$ is $\pi$-periodic (it is even $\frac{\pi}{r}$-periodic), so that it suffices to treat the special case $\lambda = 0$. When $\lambda = 0$, the points of discontinuity of $W_{r, 0}$ are those in  $\mathcal{D}_{r, 0} = \Bigl \lbrace \frac{j \pi}{r} \, ; \, j \in \lbrace 0, \, \ldots \, , r-1 \rbrace \Bigr \rbrace$. To lighten notation, we write $W$ for $W^{r, 0}$, $\beta_{\eta}$ for $\beta_{\eta}^{r, 0}$ and $W_{\eta}$ for $W_{\eta}^{r, 0}$. By $\frac{\pi}{r}$-periodicity and symmetry of $W'$ around each $\frac{j \pi}{r}$, it is enough to prove that 
	\begin{equation}\label{eq: Intg Wlamprim}
		\int_{0}^{\pi/2r} \abs{W_{\eta}'(x)}\dx \ll_{r} \abs{\log \eta}.
	\end{equation}
	For any $x \in \intv$, we have $\beta_{\eta}(x) \leq 1$ so $\abs{\Wpri(x)} \leq \abs{W'(x)} + \abs{\beta_{\eta}'(x)} \abs{W(x)}$. Write $\nrmbet$ for the sup norm of $\beta_{0}'$ on $[0, 1]$, then we have $\beta_{\eta}'(x) = 0$ for any $x \in \left[ \eta, \frac{\pi}{2r} \right]$, and $\abs{\beta_{\eta}'} \leq \frac{\nrmbet}{\eta}$ for $x \in [0, \eta]$. With the specific choice for $\beta_0$ made above, one can check that $\nrmbet = 6$. Secondly, we bound $W$. It is well-known that $\sin(u) \geq \frac{2u}{\pi}$ for any $u \in [0, \pi/2]$. This implies that, for any $x \in \left(0, \frac{\pi}{2r}\right]$, we have $\abs{W(x)} = \log \left( \frac{1}{\sin (rx)} \right) \leq \log \left( \frac{\pi}{2rx}\right)$. Finally, to control $W'$, we use the same lower bound on sin to deduce, for any $x \in \left(0, \frac{\pi}{2r}\right]$,
	\begin{equation*}
		\abs{W'(x)} = \abs{\frac{r \cos(rx)}{\sin(rx)}} \leq \frac{r}{\sin(rx)} \leq \frac{\pi}{2x}.
	\end{equation*}
	With these inequalities at hand, we can estimate the integral displayed in \eqref{eq: Intg Wlamprim}, which we split into integrals over $(0, \eta]$ and $\left[\eta, \frac{\pi}{2r}\right]$. $W_{\eta}$ is constant equal to $0$ over $\left[0, \frac{\eta}{3} \right]$, hence $\Wpri = 0$ on this interval, and
	\begin{align*}
		\int_{0}^{\eta} \abs{\Wpri(x)} \dx = \int_{\eta / 3}^{\eta} \abs{\Wpri(x)} \leq & \int_{\eta / 3}^{\eta} \abs{W'(x)} \dx \, + \, \int_{\eta / 3}^{\eta} \abs{\beta_{\eta}'(x)} \, \abs{W(x)} \dx \\
		\leq & \int_{\eta / 3}^{\eta} \frac{\pi}{2x}\, \dx \, +  \, \frac{6}{\eta} \int_{\eta / 3}^{\eta} \log \left( \frac{\pi}{2rx} \right) \dx \ll_{r} \abs{\log \eta \,}.
	\end{align*}
	To estimate the integral of $\abs{\Wpri}$ on $\left[ \eta, \frac{\pi}{2r}\right]$ we combine the bounds established above with the fact that $\beta_{\eta}'(x) = 0$ for any $x \in \left[ \eta, \frac{\pi}{2r}\right]$:
	\begin{equation*}
		\int_{\eta}^{\pi / 2r} \abs{\Wpri(x)} \dx \leq \int_{\eta}^{\pi/2r} \abs{W'(x)} \, \dx \leq \int_{\eta}^{\pi/2r} \frac{\pi}{2x} \,  \dx \ll_{r} \abs{\log \eta \, }.
	\end{equation*}
	Both implicit constants depend only on $r$, so adding the two bounds gives the desired result.	
\end{proof}
\vspace{0.5em}

We have now all the necessary tools to prove Proposition~\ref{prop: Seq intg bounded}.

\begin{proof}[\textbf{Proof of Proposition~\ref{prop: Seq intg bounded}}]
	We show that the sequence of integrals of $\Wlama$ is bounded, the proof would be the same if one replaced $\Wlama$ by $W^{r_{a}, - \lambda_{a}}$ (this corresponds to $\epsilon = -1$ in the notation of the proposition). We use the triangle inequality to obtain
	\begin{multline*}
		\abs{\int_{\intrv} \Wlama \dnua} \leq \abs{ \int_{\intrv} \Wlama - \Wetaa \, \dnua } \\ 
		+ \abs{\int_{\intrv} \Wetaa \, \dnua - \int_{\intrv}\Wetaa \, \dnust} + \abs{\int_{\intrv} \Wetaa \, \dnust}.
	\end{multline*}
	
	We claim that the first integral of the right-hand side vanishes if $\eta$ is small enough. Recall from Section~\ref{sect: Distance angles} the definition of the constants $\sigma_{p}$ and $\tau_{p, N}$. Fix $\eta := (\min \left( (q^{a})^{- \sigma_{p}}, (q^{a})^{- \tau_{p, 3}}, (q^{a})^{- \tau_{p, 6}} \right))^{1/2}$ (using the assumption $p \geq 7$, one can check that $\eta < \pi / 12$). The points of discontinuity of $\Wlama$ are either angles of Gauss sums or of the form $\frac{j \pi}{6}$ for some $j \in \lbrace 0, \, \ldots \, , 5 \rbrace$. Theorem~\ref{thm: Dist G K angles} implies that $\dist(\thetv, \Disca) \geq \eta$  for any $v \in \Pqamax$. This means that the support of the measure $\nua$ is contained in $\Delta_{\eta_{a}}^{r_a, \lambda_{a}}$. By construction of $\Wetaa$ we have $\Wlama(x) = \Wetaa(x)$ for any $x \in \Delta_{\eta_{a}}^{r_a, \lambda_a}$, so we deduce that $\int_{\intrv} \abs{\Wlama - \Wetaa}\dnua =0$.
	
	We show that the second term is small when $a \rightarrow \infty$. Indeed, the function $\Wetaa$ is continuously differentiable on $\intrv$, and we may apply Theorem~\ref{thm: Equid eff thm}. Combining it with Lemma~\ref{lem: Analytic bounds intg} gives
	\begin{equation*}
		\abs{\int_{\intrv} \Wetaa \, \dnua - \int_{\intrv}\Wetaa \, \dnust} \ll_q \frac{a^{1/2}}{q^{a/4}} \, \int_{0}^{\pi} \abs{(\Wetaa)'(x)} \dx \ll_{q} c(r_a) \,  \frac{a^{1/2}}{q^{a/4}} \, \abs{\log \eta}.
	\end{equation*}
	Having chosen $\eta = (\min \left( (q^{a})^{- \sigma_{p}}, (q^{a})^{- \tau_{p, 3}}, (q^{a})^{- \tau_{p, 6}} \right))^{1/2}$, we have $\abs{ \log \eta} \ll_{q} a$, so
	\begin{equation}
		\abs{\int_{\intrv} \Wetaa \, \dnua - \int_{\intrv}\Wetaa \, \dnust} \ll_{q} c(r_a) \frac{a^{3/2}}{q^{a/4}} \ll_{q} \max(c(3), c(6)) \, \frac{a^{3/2}}{q^{a/4}} \, .
	\end{equation}
	
	Finally, the third term is bounded independently of $a$. Indeed, $\abs{\Wetaa(x)} \leq \abs{\Wlama(x)}$ for any $x \in \intrv$, so Lemma~\ref{lem: Bound intg Wlam lemma} implies $\abs{\int_{\intrv} \Wetaa \, \dnust} \leq \log 4$ for any $a \geq 1$. Adding up the three contributions:
	\begin{equation*}
		\abs{\int_{\intrv} \Wlama \dnua} \leq 0 + O \left(\frac{a^{3/2}}{q^{a/4}} \right) + \log 4 \ll_{q} 1.
	\end{equation*}
\end{proof}

\subsection{Conclusion of the proof of Theorem \ref{thm: Ratio spval height thm}}
In this section we conclude the proof of Theorem \ref{thm: Ratio spval height thm} by establishing the lower bound appearing in it.

\begin{proof}[Proof of the lower bound in \ref{thm: Ratio spval height thm}]
	As discussed in Remark \ref{rmk: Ratio H into qa}, it is sufficient to prove that we have $\log \spvala / q^{a} \gg_{q} 1/a$. Recall from \eqref{eq: Spval over q^a} that
	\begin{equation}\label{eq: Spval over q^a bis}
		\frac{\log \spval}{q^a} = \log 4 \, \frac{\abs{\Oqa}}{q^a}  + \frac{1}{q^a} \, \sum_{o \, \in \Oqa} \log \abs{\sin \left( \frac{\epso + \theto}{2}\right) \, \sin \left( \frac{\epso - \theto}{2}\right)}.
	\end{equation}
	Lemma \ref{lem: Estimates numb orbits} $iii)$ implies that $\log 4 \abs{\Oqa}/q^a = O(1/a)$, and Proposition \ref{prop: Maj sum orbs nonmax} gives
	\begin{equation*}
		\frac{1}{q^a} \sum_{o \, \in \Oqa \setminus \Oqamax} \log \abs{\sin \left( \frac{\epso + \theto}{2}\right) \, \sin \left( \frac{\epso - \theto}{2}\right)} = O \left(\frac{1}{q^{a/2}} \right).
	\end{equation*}
	Plugging these estimates into \eqref{eq: Spval over q^a bis} we get
	\begin{equation*}
		\frac{\log \spvala}{q^{a}} = \frac{1}{q^{a}} \sum_{o \in \Oqamax} \log \abs{\sin \left( \frac{\epso + \theto}{2}\right) \, \sin \left( \frac{\epso - \theto}{2}\right)} + O \left(\frac{1}{a} \right) . 
	\end{equation*}
	Recall from Lemma \ref{lem: Estimates numb orbits} $i)$ that $\Pi_{q}(a) / q^{a} \geq  1/a$. Combining this with Proposition \ref{prop: Lower bound by integrals} implies that 
	\begin{equation*}
		\frac{\log \spvala}{q^{a}} \geq \frac{4}{a} \int_{\intrv} \left( W^{r_{a}, \lambda_{a}} + W^{r_{a}, - \lambda_{a}} \right)  \dnua \ + \  \frac{8 \, \varpi(a)}{a} \int_{\intrv} W^{r_{\lfloor a/2 \rfloor}, \lambda_{{\lfloor a/2 \rfloor}}} \, \d \nu_{ \lfloor a / 2 \rfloor} + O \left(\frac{1}{a}\right).
	\end{equation*}
	Proposition \ref{prop: Seq intg bounded} states that for any $\epsilon \in \lbrace \pm 1 \rbrace$, the integrals $\int_{\intrv} W^{r_{a}, \, \epsilon \lambda_{a}}\, \dnua$ form a bounded sequence, so the right-hand side in the inequality above is a $O(1/a)$, finally giving the desired result.
\end{proof}


\section{Large Tate--Shafarevich groups and the Brauer--Siegel ratio}

In this section we conclude the proof of Theorem~\ref{thm: Main thm introd}, and deduce various other properties on the Tate--Shafarevich groups of the surfaces $S_a$. We finally reinterpret Theorem~\ref{thm: Main thm introd} as an analogue of the Brauer--Siegel theorem.  

\subsection{Large Tate--Shafarevich groups}

From Section~\ref{sect: L func} we know that the surface $\Sa$ satisfies the BSD conjecture for all $a \geq 1$ (see Theorem \ref{thm: BSD}): this implies in particular that its Tate--Shafarevich group $\Sh(\Sa)$ is finite. We now deduce from Theorem~\ref{thm: Ratio spval height thm} the main result of our article, which describes the rate of growth of the Tate--Shafarevich groups of $S_a$.

\begin{theorem}\label{thm: Sha compared to H}
	As $a \rightarrow \infty$, we have the estimate
	\begin{equation}\label{eq: Sha compared to H}
		\abs{\Sh \left( \Sa \right)} = H(\Sa)^{1 + O(1/a)},
	\end{equation}
	The implicit constant in \eqref{eq: Sha compared to H} is effective and depends only on $q$.
\end{theorem}

\begin{proof}
	By Corollary~\ref{cor: Alg rank 0}, we have $\Reg(\Sa) = 1$. Taking the logarithm of both sides of the BSD formula \eqref{eq: Strong BSD} and rearranging terms gives
	\begin{equation}\label{eq: BSD logarithm}
		\frac{\log \abs{\Sh (\Sa)}}{\log H(\Sa)} = 1 + \frac{\log \spval}{\log H(\Sa)} - \frac{\log \left( q^2 \, \prod_{v \in \Mpl} c_{v}(\Sa) \right) - 2\, \log  \abs{\tor{S_{a}(\K)}}}{\log H(\Sa)}.			
	\end{equation}
	From Section~\ref{sect: Invariants} we know that $\prod_{v \in \Mpl} c_{v}(\Sa) = c_{\infty}(\Sa) \divi 9$. On the other hand, Theorem $3.8$ in \cite{HindryPacheco} states that $\abs{\tor{S_{a}(\K)}}^{2} \ll_{q} h(S_a)^{4}$. Proposition \ref{prop: Height} yields $h(S_a) = q^{a} +1$, so we deduce that $\frac{2 \log \abs{\tor{S_{a}(\K)}}}{\log H(S_a)} \ll_q \frac{a}{q^{a}}$. On the other hand, Theorem~\ref{thm: Ratio spval height thm} states $\log \spvala / \log H(S_a) = O(1/a)$, so substituting all these estimates in \eqref{eq: BSD logarithm} gives
	\begin{equation}\label{eq: Sha over height}
		\frac{\log \abs{\Sh (\Sa)}}{\log H(\Sa)} = 1 + O \left( \frac{1}{a} \right).
	\end{equation}
	The statement follows by exponentiating \eqref{eq: Sha over height}.
\end{proof}

\subsection{An analogue of the Brauer--Siegel theorem}

In this section we view Theorem \ref{thm: Sha compared to H} as an analogue of the Brauer--Siegel theorem. In \cite{HindryPacheco}, Hindry and Pacheco defined the Brauer--Siegel ratio of an abelian variety $A / \K$ whose Tate--Shafarevich group is finite. It is given by
\begin{equation*}
	\BS(A) := \frac{\log \left( \abs{\Sh(A)} \, \Reg(A) \right)}{\log H(A)}.
\end{equation*}
In our setting, for any $a \geq 1$, the abelian surface $S_a$ has Mordell--Weil rank $0$, so its Néron--Tate regulator is trivial and the Brauer--Siegel ratio is $\BS(S_a) = \log \abs{\Sh(S_a)} / \log H(S_a)$. With this notation, Theorem \ref{thm: Sha compared to H} may be rewritten in the following way:

\begin{corollary}\label{cor: Brauer Siegel tends to 1}
	As $a \rightarrow \infty$, we have
	\begin{equation*}
		\BS(S_a) = 1 + O(1/a).
	\end{equation*}
\end{corollary}

We refer the reader to \cite{Hindry07, HindryPacheco} for a detailed discussion on the definition and properties of the Brauer--Siegel ratio and more on the analogy with the original Brauer--Siegel theorem.
\vspace{1em}

In \cite{Ulmer19}, Ulmer associates a new invariant, called the \textit{dimension} of $\Sh$, to an abelian variety $A / \K$ with finite Tate--Shafarevich group. For any $n \geq 1$, let $\K_{n} := \Ff_{q^{n}}(t)$ and define
\begin{equation*}
	\dim \Sh (S_a) := \lim_{n \rightarrow \infty} \frac{\log \abs{\Sh \left( S_a \times_{\K} \K_{n} \right)[p^{\infty}]}}{\log \left(q^{n} \right)}.
\end{equation*}  
\noindent Ulmer has proved (see \cite[Proposition $4.1$]{Ulmer19}) that this limit indeed exists and is a non-negative integer. One of the interests of $\dim \Sh(S_a)$ is that it provides a lower bound on $\BS(S_a)$ (see Proposition 4.6 in \cite{Ulmer19}):
\begin{equation}\label{eq: Dim Sh vs BS}
	\frac{\dim \Sh(S_a)}{h(S_a)} \leq \BS(S_a) + o(1) \qquad \text{ as } a \rightarrow \infty.
\end{equation}

\begin{proposition}\label{prop: Dimension Sh}
	For any $a \geq 1$, we have 
	\begin{equation*}
		\dim \Sh(S_a) = 0.
	\end{equation*}
\end{proposition}

\begin{proof}
	We use Proposition $4.2$ from \cite{Ulmer19}, which gives a closed formula for $\dim \Sh(S_a)$ in terms of $h(S_a)$ and the $p$-adic slopes of $L(S_a, T)$. For any $o \in \Oqa$ and any $i \in \lbrace 1, 2 \rbrace$, the Newton polygon of $1 - \gamo \kappa_{i}(o) T^{\abs{o}} \in \Zbar[T]$ is a segment of length $\abs{o}$ and slope $\valp(\gamo \kappa_{i}(o))$. With our notation, the proposition mentioned above implies 
	\begin{equation*}
		\dim \Sh(S_a) = h(S_a) -2 - \sum_{o \in \Oqa} \sum_{i = 1}^{2} \max \bigl \lbrace 0, \abs{o} - \valp(\gamo \kappa_{i}(o)) \bigr \rbrace.
	\end{equation*}
	Using \refKloos{Kloos conj valuation}, and switching $\kapo$ and $\kapobis$ if necessary, we may assume that $\valp(\kapo) = 0$ and $\valp(\kapobis) = \abs{o}$. 
	
	When $q \equiv 1 \mod 3$, Lemma \ref{lem: Valuations gauss} implies that for any $o \in \Oqa$, we have
	\begin{equation*}
		\valp(\gamo \kapo) = \begin{cases}
			\frac{2\abs{o}}{3} & \text{ if } \prun(o) = 1, \\
			\frac{\abs{o}}{3} & \text{ if } \prun(o) = -1, 
		\end{cases}
	\quad \text{ and } \quad \valp(\gamo \kapobis) = \begin{cases}
		\frac{5\abs{o}}{3} & \text{ if } \prun(o) = 1, \\
		\frac{4\abs{o}}{3} & \text{ if } \prun(o) = -1.
	\end{cases}
	\end{equation*}
	Therefore, we get
	\begin{equation*}
		\sum_{o \in \Oqa} \sum_{i = 1}^{2} \max \bigl \lbrace 0, \abs{o} - \valp(\gamo \kappa_{i}(o)) \bigr \rbrace = \frac{1}{3} \sum_{\substack{o \in \Oqa \\ \prun(o) = 1}} \abs{o} \  + \ \frac{2}{3} \sum_{\substack{o \in \Oqa \\ \prun(o) = -1}} \abs{o}.
	\end{equation*}
	Since $q \equiv 1 \mod 3$, the map $\Theta$ establishes a bijection between $\ensemble{o \in \Oqa}{\prun(o) = 1}$ and $\Pqa$ (in the same way $\ensemble{o \in \Oqa}{\prun(o) = -1}$ and $\Pqa$ are also in bijection). Moreover, we have $\abs{o} = \abs{\Theta(o)}$ for any $o \in \Oqa$, so we get
	\begin{equation*}
		\sum_{o \in \Oqa} \sum_{i = 1}^{2} \max \bigl \lbrace 0, \abs{o} - \valp(\gamo \kappa_{i}(o)) \bigr \rbrace = \left(\frac{1}{3} + \frac{2}{3} \right) \sum_{v \in \Pqa} \abs{v} = \abs{\Fqa\st} = q^{a}-1.
	\end{equation*}

	When $q \equiv 2 \mod 3$, we deduce from Lemma \ref{lem: Valuations gauss} that
	\begin{equation*}
		\valp(\gamo \kapo) = \abs{o}/2 \quad \text{ and } \quad \valp(\gamo \kapobis) = 3\abs{o}/2.
	\end{equation*}
	Following the same line of reasoning as above:
	\begin{multline*}
		\sum_{o \in \Oqa} \sum_{i = 1}^{2} \max \bigl \lbrace 0, \abs{o} - \valp(\gamo \kappa_{i}(o)) \bigr \rbrace = \sum_{o \in \Oqa} \left( \abs{o} - \valp (\gamo \kapo) \right) \\ 
		= \sum_{o \in \Oqa} \frac{\abs{o}}{2} = \frac{1}{2} \abs{\left(\Z / 3 \Z\right)\st \times \Fqa\st} = q^{a}-1.
	\end{multline*}
	In both cases we obtain
	\begin{equation*}
		\dim \Sh (S_a) = q^{a} + 1 -2 -(q^{a} - 1) = 0.
	\end{equation*}
\end{proof}

\noindent In \cite{Ulmer19}, Ulmer shows that different families of abelian varieties have Brauer--Siegel ratios that tend to $1$ by proving that their dimension of $\Sh$ are large enough. Notice that, in the case of the sequence $(S_a)_{a \geq 1}$, Proposition \ref{prop: Dimension Sh} implies that $\dim \Sh(S_a) / h(S_a) = 0$ for any $a \geq 1$, meaning that the algebraic approach presented in \cite{Ulmer19} is not sufficient to prove Corollary \ref{cor: Brauer Siegel tends to 1}.
\vspace{0.75em}

The Brauer--Siegel ratio has been studied for different sequences of elliptic curves over function fields (see \cite{HindryPacheco,  Gri19, GriUlm, GriDWit}) for which the Brauer--Siegel ratio is almost always near $1$. In \cite{AGTT, Ulmer19}, the authors study some families of abelian varieties of dimension greater than $1$ whose Brauer--Siegel ratios tend to $1$. The family $\left( S_a\right)_{a \geq 1}$ studied in this article constitutes a new example of higher-dimensional abelian varieties over $\K$ that are non-isotrivial and whose Brauer--Siegel ratios tend to $1$. This provides some more evidence that the Brauer--Siegel ratio of an abelian variety over $\K$ should ``often" be close to $1$ when the exponential height is large, a phenomenon already observed in the articles mentioned above (see \cite{Hindry07, HindryPacheco} for a discussion on why this is the case).


\noindent\hfill\rule{7cm}{0.5pt}\hfill\phantom{.}

\paragraph{Acknowledgements}
It is a pleasure to thank Richard Griffon for his encouragements and useful comments on earlier versions of this work. Special thanks go to Qing Liu and Céline Maistret for taking the time to explain details about the geometry of curves. The author would also like to thank Marc Hindry and Douglas Ulmer for fruitful discussions about various parts of the paper. This work was funded by the École Doctorale des Sciences Fondamentales, Université Clermont Auvergne.


\normalsize\vfill
\noindent\rule{7cm}{0.5pt}

\smallskip
\noindent
\textsc{Martin Azon} (\textit{martin.azon@uca.fr}) --
{\sc Laboratoire de Mathématiques B. Pascal, Université Clermont Auvergne,} 
Campus des Cézeaux. 3 place Vasarely, TSA 60026 CS 60026, 63178 Aubière Cedex (France).

\end{document}